\documentclass{alggeom} 
\usepackage{latexsym,xspace,enumerate,amsfonts,amsmath,amssymb,amscd}
\usepackage{xypic,color,hyperref,syntonly,slashed,paralist,stmaryrd,graphicx}
\usepackage{bbm}
\usepackage{thmtools}
\usepackage{tikz-cd}
\usepackage{etex}
\usepackage{amsfonts}
\usepackage{amscd}
\usepackage{amsbsy}
\usepackage{mathtools}
\usepackage{amsxtra}
\usepackage{amssymb}
\usepackage{epsfig}
\usepackage{epic,eepic}
\usepackage{graphicx}
\usepackage[all]{xy}
\usepackage{psfrag}
\usepackage{amsmath}
\usepackage{amsthm}
\usepackage{setspace}
\usepackage{hyperref}
\usepackage{url}
\usepackage{here}
\usepackage{todonotes}
\usepackage{xparse}
\usepackage{xspace}
\usepackage{enumitem}

\PassOptionsToPackage{hyphens}{url}

\DeclareMathAlphabet\mathbfcal{OMS}{cmsy}{b}{n}
\usepackage{adjustbox}

\DeclareRobustCommand{\SkipTocEntry}[9]{}

\newcommand*{\calA}{\mathcal A}
\newcommand*{\calB}{\mathcal B}

\newcommand*{\calE}{\mathcal E}
\newcommand*{\calF}{\mathcal F}
\newcommand*{\calG}{\mathcal G}
\newcommand*{\calH}{\mathcal H}
\newcommand*{\calI}{\mathcal I}

\newcommand*{\calL}{\mathcal L}
\newcommand{\calM}{{\mathcal M}}
\newcommand{\calO}{{\mathcal O}}

\newcommand{\calQ}{{\mathcal Q}}
\newcommand{\calN}{{\mathcal N}}
\newcommand*{\calR}{\mathcal R}
\newcommand*{\calS}{\mathcal S}

\newcommand*{\calX}{\mathcal X}

\newcommand{\calY}{{\mathcal Y}}

\renewcommand{\AA}{{\mathbb{A}}}
\newcommand{\BB}{{\mathbb{B}}}
\newcommand{\CC}{{\mathbb{C}}}

\newcommand{\GG}{{\mathbb{G}}}

\newcommand{\PP}{{\mathbb{P}}}
\newcommand{\QQ}{{\mathbb{Q}}}

\newcommand{\ZZ}{{\mathbb{Z}}}

\newcommand{\st}{\mathrm{st}}

\newcommand{\GL}{\mr{GL}}
\newcommand{\SL}{\mr{SL}}

\newcommand{\mr}{\mathrm}

\newcommand{\wt}{\widetilde}
\newcommand{\End}{\mathop{\rm End}\nolimits}

\newcommand{\Jac}{\mathop{\rm Jac}\nolimits}
\newcommand{\Id}{\mathop{\rm Id}\nolimits}

\renewcommand{\ker}{\mathop{\rm ker}\nolimits}

\renewcommand{\deg}{\mathop{\rm deg}\nolimits}
\newcommand{\rk}{\mathop{\rm rk}\nolimits}
\newcommand{\Tr}{\mathop{\rm Tr}\nolimits}

\newcommand{\Hit}{\mathop{\rm Hit}\nolimits}
\newcommand{\Pic}{\mathop{\rm Pic}\nolimits}

\newcommand{\SD}{\mathop{\rm SD}\nolimits}

\renewcommand{\mod}{\mathop{\rm mod}\nolimits}

\renewcommand{\div}{\mr{div}\,}

\newcommand{\note}[1]{\marginpar{\raggedright\if@twoside\ifodd\c@page\raggedleft\fi\fi\sf\scriptsize \red{RMK: #1}}}

\newcommand\red[1]{\textcolor{red}{#1}}

\newcommand{\be}{\begin{equation}}
\newcommand{\ben}{\begin{equation}\nonumber}
\newcommand{\ee}{\end{equation}}
\newcommand{\bp}{\begin{para}}
\newcommand{\ep}{\end{para}}
\newcommand{\bps}{\begin{paras}}
\newcommand{\eps}{\end{paras}}

\newcommand{\frakm}{\mathfrak{m}}

\newcommand{\bflambda}{{\boldsymbol{\lambda}}}

\newcommand{\bftau}{{\boldsymbol{\tau}}}

\newcommand{\bfPhi}{{\boldsymbol{\Phi}}}

\newcommand{\bfz}{{\mathbf{z}}}

\newcommand{\bfq}{{\mathbf{q}}}
\newcommand{\bfd}{{\mathbf{d}}}

\newcommand{\IVC}{\mathrm{IVC}}
\newcommand{\Higgs}{\mathrm{Higgs}}
\newcommand{\Nm}{\mathsf{Nm}}

\newcommand{\Spec}{\mathop{\rm Spec}\nolimits}
 
\newtheorem{proposition}{\textbf{Proposition}} [section]
\newtheorem{lemma}[proposition]{\textbf{Lemma}} 

\newtheorem{remark}[proposition]{\textbf{Remark}}
\newtheorem{theorem}[proposition]{\textbf{Theorem}}

\theoremstyle{definition}
\newtheorem{definition}[proposition]{\textbf{Definition}}
\newtheorem{example}[proposition]{\textbf{Example}}
\newtheorem*{example*}{\textbf{Example}}
\newtheorem*{theorem*}{\textbf{Theorem}}

\def\be{\begin{equation}}   \def\ee{\end{equation}}     \def\bes{\begin{equation*}}    \def\ees{\end{equation*}}
\def\ba{\be\begin{aligned}} \def\ea{\end{aligned}\ee}   \def\bas{\bes\begin{aligned}}  \def\eas{\end{aligned}\ees}
\def\={\;=\;}  \def\+{\,+\,}          

\newcounter{para}[section]
\newenvironment{para}[2][]{\refstepcounter{para}\noindent\ignorespaces{\bf #1\thepara. #2.} \rmfamily}{\noindent\ignorespacesafterend\bigskip}
\newenvironment{paras}[1]{\noindent\ignorespaces{\bf #1.} \rmfamily}{\noindent\ignorespacesafterend\bigskip}

\numberwithin{proposition}{section}
\numberwithin{definition}{section}

\newcommand{\can}{\mr{can}}

\newcommand{\ol}{\overline}

\def\wR{\widehat R}  
\def\wX{\widehat X}

\newcommand{\whmu}{\widehat{\mu}}

\def\wh#1{\widehat{#1}}

\newcommand{\moduli}[1][g]{{\mathcal M}_{#1}}

\newcommand{\barmoduli}[1][g]{{\overline{\mathcal M}}_{#1}}

\DeclareDocumentCommand{\barPr}{O{\wh{\pi}} O{ }}{{\overline{\mathcal{P}}}_{#2}({#1})}

\DeclareDocumentCommand{\tildecJ}{O{g,n} O{d}}{{\widetilde{\mathcal J}}_{#1}^{#2}}
\DeclareDocumentCommand{\barcJ}{O{g,n} O{d}}{{\overline{\mathcal J}}_{#1}^{#2}}

\DeclareDocumentCommand{\LMS}{ O{\mu} O{g,n}} {\Xi\overline{\mathcal{M}}_{#2}(#1)}

\DeclareDocumentCommand{\kLMS}{ O{\mu} O{g,n}} {\Xi^k\!\overline{\mathcal{M}}_{#2}(#1)}

\DeclareDocumentCommand{\QMS}{ O{\mu} O{g,n}} {\overline{\mathcal{Q}}_{#2}(#1)}

\newcommand*{\msds}{multi-scale differentials\xspace}

\newcommand{\iomoduli}[1][g,n]{{\Omega^i\!\mathcal M}_{#1}}

\usetikzlibrary{intersections} 
\tikzset{
	every loop/.style={very thick},
	comp/.style={circle,black,draw,thin,inner sep=0pt,minimum size=5pt,font=\tiny},
	order bottom left/.style={pos=.05,left,font=\tiny},
	order top left/.style={pos=.9,left,font=\tiny},
	order bottom right/.style={pos=.05,right,font=\tiny},
	order top right/.style={pos=.9,right,font=\tiny},
	order middle right/.style={draw,thin,shape=rectangle,inner sep=2pt,outer sep=5pt,pos=.5,right,font=\tiny},
	order middle left/.style={draw,thin,shape=rectangle,inner sep=2pt,outer sep=5pt,pos=.5,left,font=\tiny},
	order node dis/.style={text width=.75cm},
	circled number/.style={circle, draw, inner sep=0pt, minimum size=12pt},
	below left with distance/.style={below left,text height=10pt},
	below right with distance/.style={below right,text height=10pt}
}

\DeclareRobustCommand{\SkipTocEntry}[9]{}

\usepackage[style=alphabetic,
            isbn=false,
            doi=false,
            url=false,
            backend=bibtex,
            maxnames =5, 
            maxalphanames = 6,
            block=space,
           ]{biblatex}
\defbibheading{bibliography}[\bibname]{%
  \section*{References}%
}

\DeclareFieldFormat[article]{title}{{\it #1}}
\DeclareFieldFormat{journaltitle}{{\rm #1}}
\renewbibmacro{in:}{%
  \ifentrytype{article}{}{\printtext{\bibstring{in}\intitlepunct}}}
\DeclareFieldFormat[incollection]{title}{{\it #1}}
\DeclareFieldFormat{journaltitle}{{\rm #1}}

\bibliography{BibHiggs}

\title[Spectral data on semistable curves]
    {Compactifying the rank two Hitchin system via spectral data on semistable curves}

\date{\today}

\begin{document}

\author{Johannes Horn}
\email{horn@math.uni-frankfurt.de}
\address{Goethe-Universität Frankfurt, Robert-Mayer-Str. 6-8, 60325 Frankfurt am Main}

\author{Martin M\"oller}
\email{moeller@math.uni-frankfurt.de}
\address{Goethe-Universität Frankfurt, Robert-Mayer-Str. 6-8, 60325 Frankfurt am Main}

\classification{14D06, 14H40, 14H70, 14H10}
\keywords{Hitchin systems, quadratic differentials}

\begin{abstract}
We study resolutions of the rational map to the moduli space of stable curves that associates with a point in the Hitchin base the spectral curve. In the rank two case the answer
is given in terms of the space of quadratic multi-scale differentials
introduced in [BCGGM3]. This space defines a compactification (of the
projectivization) of the regular locus of the $\GL(2,\CC)$-Hitchin base
and provides a compactification of the Hitchin system by compactified
Jacobians of pointed stable curves.
\par
We show how the classical $\GL(2,\CC)$- and $\SL(2,\CC)$-spectral
correspondence extend to the compactified Hitchin system by a
correspondence along an admissible cover between torsion-free rank $1$ sheaves
and (multi-scale) Higgs pairs of rank $2$.
\end{abstract}

\maketitle

\tableofcontents

\section{Introduction}\label{sec:intro}

The complexity of the Hitchin fibration on the moduli space of Higgs bundles
stems from the variety of singularities of the spectral curve. Already in the
case of $\GL(2,\CC)$-Higgs bundles non-reduced curves and all singularities
locally given by the equation $\lambda^2-z^\ell$ have to be considered. In
this paper
we investigate the possible modifications of the Hitchin base over which the
family of smooth spectral curves on the regular locus can be extended as
a family of \emph{semistable curves}. One might ask for a minimal modification,
and for a modification that allows for a nice modular interpretation and
a spectral correspondence. We give an answer to both questions in the rank
two case.
\par
The construction of the modified Hitchin base is motivated
by the compactification of strata of quadratic differentials in \cite{LMS}
and works in families as the Riemann surface~$X$ varies. In particular, it
allows to study the (modified) Higgs moduli space on any stable degeneration
of~$X$. Our
definition of multi-scale Higgs pairs works uniformly for smooth and stable~$X$.
We compare this viewpoint with existing constructions on degenerating families.
\par
\medskip
\paragraph{\bf Singular fibers} We start with a summary of the knowledge about
the singular Hitchin fibers. In general the spectral curve is a complex
algebraic curve, that can be reducible and non-reduced, with planar singularities. The corresponding Hitchin fibers were put in  bijection with torsion-free sheaves on the spectral curve by \cite{BNR} and \cite{Schaub}.
For $G=\SL(2,\CC)$ their geometry was studied in \cite{GotOli} using parabolic modules. From the point of view of algebraically completely integrable systems the singular fibers were studied in \cite{HiCritical} establishing lower-dimensional integrable systems supported on the critical locus. The fibers of these
sub-integrable systems are abelian varieties and were augmented in \cite{Horn} by non-abelian coordinates to a full description of singular Hitchin fibers with integral spectral curve for $G=\SL(2,\CC)$. 
The non-abelian parameters can be easily understood as \emph{stratified space}. The global description is more delicate revealing a list of moduli of non-abelian (Hecke) parameters indexed by the singularities of type $A$. For $A_k$ with $k=1,2$ one obtains a $\PP^1$ for $k=3,4$ a $\PP(1,1,2)$ and in general a complex projective variety of dimension $\lfloor\frac{k}{2} \rfloor$. Another complication is a relation between the abelian and non-abelian parameters in the case of multi-branched singularities. One motivation for the present work is to reduce the complexity by only allowing nodes for singularities. 
\par
\medskip
\paragraph{\bf The spectral map and its resolution} Let $\AA_k^\circ$ be the
Hitchin base for $\SL(k,\CC)$, the sum of the vector spaces of global $k$-differentials on a fixed complex projective curve $X$
for $k \geq 2$.  Associating with~$X$ and a point in the Hitchin base the
spectral curve defines a rational map to the moduli space of stable curves
that we call \emph{spectral map}. Our first result in Section~\ref{sec:spectralmap}
determines the resolution of this rational map in the case $k=2$, and equivalently
for $\GL(2,\CC)$-Higgs bundles. There are several possibilities to interpret
this problem. A minimal resolution of a rational map is always defined as the
closure of the graph, in our special case inside $H^0(X,K_X^2) \times
\overline{\calM}_{4g-3}$, since $\AA^\circ_2 = H^0(X,K_X^2)$.  We characterize the
points in this closure in Theorem~\ref{thm:charrank2}. Our characterization
relates the
position of the marked points on each component of the spectral curve to the
existence of a twisted quadratic differential, in spirit close to \cite{strata}.
\par
Working with this minimal resolution as a replacement of the moduli space of
Higgs bundles has several flaws. One major issue being the absence of well-defined coverings
of the limiting spectral curves to (a modification) of the original curve. Hence there
is no way to define a spectral correspondence in this setting. From this perspective
it is more natural to instead look at the admissible cover compactification of the
spectral map, i.e., consider the closure $\widetilde{\calH}_g$ of the image of the
map $H^0(X,K^2) \rightarrow \calH_{g}$, where $\calH_{g}$ is the admissible cover
compactification of simply branched covers of nodal curves of (arithmetic) genus $g$
by nodal curves of genus $\hat{g}=4g-3$. Now the push-forward of a coherent sheaf
on the covering curve will yield a sheaf on the base curve. However, since
$\widetilde{\calH}_{g}$ does not store scaling information for the differentials
we do not quite get Higgs bundles there.
\par
A 'compactification' of strata of differentials with prescribed type of zeros
is the space of multi-scale differentials that was introduced in \cite{LMS} for
abelian differentials and extended to $k$-differentials in \cite{CoMoZa}.
(Strictly speaking, only the $\CC^*$-quotient by the action of rescaling the
differentials is a compact space.) To obtain a space resolving the spectral map we 
look at the compactification of the open strata of quadratic differentials of
type $\mu=(1^{4g-4})$  and denote the moduli space of
quadratic multi-scale differentials of type $\mu$ by $\QMS$.
\par
A point in $\QMS$ contains the datum of a $4g-4$-pointed stable curve~$(X,\bfz)$
together with an order on the dual graph, that partitions~$X$ into the subcurves
$X_i$ at various levels, see Section~\ref{sec:comp} for details. It also
contains a collection~$\bfq = (q_i)$ of non-zero quadratic differentials
indexed by the levels. These differentials are required to vanish at the marked
points and satisfy compatibility
conditions at the nodes. Finally it contains the datum of a double covering
$\wh{X} \to X$ such that the pullback of $\bfq$ is square of a
collection $\bflambda = (\lambda_i)$ of abelian differentials.
In particular there is a birational forgetful map
\[ \QMS \rightarrow \calH_{g}.
\]
and, by composition, also to the minimal resolution of the spectral map.
Obviously the space of multi-scale differentials also comes with a
natural forgetful map to $\barmoduli[g]$.
\par
We first study the fiber $\BB_{X_\st} \subset \QMS$ over a fixed Riemann
surface~$X_\st$ and return below to global aspects. By definition a quadratic multi-scale differential is contained in $\BB_{X_\st}$, if the stable unpointed curve underlying $(X,\bfz)$ is $X_{\st}$. We call $\BB_{X_\st}$ the \emph{modified Hitchin base} since it comes with a natural forgetful birational surjective map
$b: \BB_{X_\st} \to \AA_2^\circ \setminus \{0\}$. In Section~\ref{sec:Hitchinbase} we show
that this map is a bijection over union of the regular locus and the locus where
merely two points collide. If $g(X_\st)=2$ the modified Hitchin base does not depend
on~$X_\st$ and we exhibit
its complete geometry in Example~\ref{ex:g2Hitbase}.
However, if $g(X_\st) \geq 3$, then $\BB_{X_\st}$ does depend on~$X_\st$.
\par
\medskip
\paragraph{\bf Multi-scale spectral data} We propose to turn the classical
spectral correspondence upside down and \emph{define} a modification of the space
of Higgs bundles to be the universal compactified Jacobian over the modified
Hitchin base. The details are most transparently stated in the case of
trace-free $\GL(2,\CC)$-Higgs bundles. We address the additional twists in
the $\SL(2,\CC)$-case in Section~\ref{sec:fixeddet}.
The general $\GL(2,\CC)$-case simply requires to record an abelian differential
in addition to the data used in the trace-free $\GL(2,\CC)$-case.
\par
A trace-free $\GL(2,\CC)$-\emph{multi-scale spectral datum on~$X_\st$} is a quadratic
multi-scale differential $(X,\bfq) \in \BB_{X_\st}$ together with a torsion-free
rank-one sheaf~$\calF$ on the associated double cover~$\wh{X}$ of~$X$. The
stack $\SD_{X_\st}$ of multi-scale spectral data on~$X_\st$ is a universally closed smooth Artin stack and comes with natural
forgetful maps $h: \SD_{X_\st} \to \BB_{X_\st}$ and $\SD_{X_\st} \to \barmoduli[4g-3]$. The latter
map justifies our claim to work in a setting of semistable curves.
The next statement summarizes the geometry of the space~$\SD_{X_\st}$. 

\par
\begin{proposition} \label{intro:modHFvsHF}
  There is a rational map $S: \SD_{X_\st} \dashrightarrow
\Higgs^\circ_{\GL(2,\CC)}(X_\st) \setminus \calN$ to the space of trace-free
$\GL(2,\CC)$-Higgs bundles on~$X_\st$ with image in the complement of the
nilpotent cone~$\calN$, such that the diagram
\be \label{eq:modHFvsHF}
\begin{tikzcd}
\SD_{X_\st}  \arrow[dashed]{r}{S} 
\arrow{d}{h}  & \Higgs^\circ_{\GL(2,\CC)}(X_\st) \setminus \calN  \arrow{d}{\Hit}\\
\BB_{X_\st} \arrow{r}{b}& \AA_2 \setminus \{0\}
\end{tikzcd}
\ee        
commutes. The rational map~$S$ is defined and an isomorphism over the locus
where the quadratic differential $q\in \AA_2$ has simple zeros or at most
one double zero.
\end{proposition}
\par
The definition of $\SD_{X_\st}$ requires the choice of a numerical polarization on the double cover $\wh{X}$. In our setup a natural choice is the pointed canonical polarization. With respect to this polarization the degree $d$ component of $\SD_{X_\st}$ is a proper Deligne-Mumford stack if and only if $\gcd(d-2g+2,6g-6)=1$. 
\par
The $P=W$-conjecture has recently motivated a lot of research on the
perverse filtration of the Hitchin map. E.g.\ in \cite{CHMsupport} the authors
determine support loci with the help of compactified Jacobians on versal
deformations of the singular fiber. It would be interesting to see if
the (well-understood) map~$h$ could lead to additional insights.
\par
The space $\SD_{X_\st}$ has two more features that we address in the next
paragraphs.
\par
\medskip
\paragraph{\bf Multi-scale spectral correspondence} While defined in
terms of spectral data, i.e.\ by data on the 'spectral' curve $\wh{X}$,
the stack $\SD_{X_\st}$ has an interpretation as moduli stack of Higgs pairs
via a spectral correspondence. To begin with we specify certain
Higgs pairs that appear. As in \cite{LMS} and as in the definition of~$\SD_{X_\st}$
these Higgs pairs will be defined level by level on some pointed stable
curve $(X,\bfz)$ that stabilizes after forgetting the marked points
to its top level curve, which is $X_\st$. We define a \emph{trace-free
multi-scale $\GL(2,\CC)$-Higgs pair} to be a tuple consisting of the
following objects. First, a \emph{special torsion-free rank two
sheaf $\calE$} on $X$, i.e.\ $\calE$ is required to be locally free, except for a
special local form at all nodes. Second, it contains an equivalence class
of a collection of trace-free Higgs fields $\bfPhi = (\phi_i)$ on each level
of~$X$. These levelwise Higgs fields are meromorphic with poles and zeros of
higher order at the preimages of the nodes. We refer the reader to Definition~\ref{def:Higgspair}
for full details. We emphasize that the setup allows~$X_\st$ to vary. The following is the special case
of Theorem~\ref{thm:Higgspairs}, in which~$X_\st$ is fixed but may itself be a stable curve.
It should be viewed as stable curve generalization of the classical
BNR-correspondence recalled for comparison in
Section~\ref{sec:HitchinBackground}. 
\par
\begin{theorem} \label{intro:BNR} Given a quadratic multi-scale differential
$(\wh{X} \to X, \bfz, \bfq)$ there is a bijective correspondence between
\begin{itemize}    
\item[i)] torsion-free rank 1 sheaves on $\wh{X}$, and
\item[ii)] special pairs of a torsion-free rank two sheaf $\calE$ on~$X$ and a trace-free Higgs
field $\bfPhi$ with determinant $\bfq$.
\end{itemize}
Equivalently, there is bijective correspondence between trace-free
$\GL(2,\CC)$-multi-scale spectral data and trace-free multi-scale $\GL(2,\CC)$-Higgs
pairs.
\par
Moreover, the correspondence respects natural notions of (semi)stability
on spectral data and Higgs pairs. 
\end{theorem}
\par
The precise definition of stability is given in Section~\ref{sec:spectralcorr}. In Section~\ref{sec:fixeddet} we define a substack of $\SD_{X_\st}$ of spectral data satisfying a Prym condition. We show that the above theorem specializes to a correspondence between torsion-free sheaves satisfying the Prym condition and multi-scale $\SL(2,\CC)$-Higgs
pairs generalizing the classical $\SL(2,\CC)$-spectral correspondence (Theorem \ref{thm:sl_spectral_corr}). 
\par In the case of a quadratic differential on $X_\st$ with one double zero and all other zeros simple the original spectral curve is nodal and $\wh{X}$ is the semistable model obtained from the normalization of the spectral curve by putting a rational bridge at the preimages of the node. In this way one obtains an admissible cover $\wh{X} \to X$, where $X$ is the smooth curve $X_\st$ augmented by a rational tail with a node at the double zero. A multi-scale Higgs pair on $X$ restricts to a meromorphic Higgs bundle on the rational tail, such that at the preimage of the node the Higgs field is diagonal with a pole of order 3. Such meromorphic Higgs bundles appeared recently in the work of Ivan Tulli \cite{Tulli}, who showed that a certain moduli space of Higgs bundles of this kind on a rational curve realizes the Ooguri-Vafa space. The Ooguri-Vafa hyperkähler metric was conjectured in \cite{Neitzke} to be part of the local model for the approximate description of the Hitchin hyperkähler metric at a generic point of the discriminant locus.
\par
\medskip
\paragraph{\textbf{Comparison to original Hitchin fiber}}
In Section \ref{sec:Compare}, we will compare the original Hitchin fibers to the $\phi$-compactified Jacobians over the modified Hitchin base with respect to the pointed canonical polarization $\phi$. There are similarities in special cases, mostly the cases of quadratic differentials with at most double zeros, but in general
the fibers look quite different. For example, the classical Hitchin fibers are stratified by fiber bundles over the Jacobian of the normalized spectral curve. In the present work, the singular spectral curves are replaced by pointed stable curves. Their compactified Jacobians instead are stratified by fiber bundles over the Jacobian of the normalized pointed stable curve with the normalized spectral curve being only one of several connected components. Moreover, the compactified Jacobians of pointed stable curves have irreducible components indicated by the multi-degrees compatible with the stability condition. In contrast in the classical setting, the $\GL(2,\CC)^\circ$-Hitchin fibers are irreducible, whenever there is at least one zero of the quadratic differential of odd order. We will study some special cases in more detail to give references for these differences.
\par
\medskip
\paragraph{\bf Degeneration of the curve~$X$} Applying degeneration techniques
is one of several reasons to consider the moduli space of Higgs bundles on a family
of curves $\calX$ degenerating to a stable curve. Our construction of the stack
of multi-scale spectral data automatically works in families. For comparison
we summarize the known constructions.
\par
A natural approach is to compactify the family of moduli spaces on the generic fiber by the moduli space of Hitchin pairs $(\calE,\Phi)$ on the special fiber $X$ made up of a \emph{torsion-free sheaf} $\calE$
and a morphism $\calE \to \calE \otimes \omega_{X}$. Using Simpson's method for constructing moduli spaces, one can define a moduli space of semistable Hitchin pairs on $X$. However, this moduli space is missing some desirable properties. For example there is no well-defined Hitchin map on the moduli space of Hitchin pairs on the stable curve~$X$.
\par
A resolution to these drawbacks is suggested in the work of Balaji, Barik and Nagaraj \cite{BBN}, however only in the case
of a family degenerating to a stable curve with a single node. Building on the work of Nagaraj-Seshadri \cite{NagSesh} the idea is to consider a modification $\calX^{mod}$ of the family of curves obtained by blowing up $\calX$ repeatedly at the node. One can show that for any family of torsion-free sheaves there exists a modification $\calX^{mod}$, such that it corresponds to a family of locally free sheaves on $\calX^{mod}$. Hence, a family of Higgs pairs on $\calX$ corresponds to a family of Higgs bundles on some modification $\calX^{mod}$. This data of $(\calX^{mod},E,\Phi)$ is referred to as Gieseker-Hitchin pair. Using locally free sheaves the moduli space of Gieseker-Hitchin pairs allows the definition of a Hitchin morphism extending the classical Hitchin map for the generic fiber of~$\calX$. 
\par 
Our approach is opposite. Our degeneration of the Hitchin system on a family of curves $\calX$ as above features a morphism to the (modified) Hitchin base by definition. Then a spectral correspondence dictates a moduli space of multi-scale $\GL(2,\CC)$-Higgs pairs on the pointed stable curve $X$ associated with it. When the zeros of the quadratic differential do not collide, we encounter the Hitchin pairs of the naive approach mentioned before with a special local form at all nodes. However, our definition of semistability is different. Taking into account the zeros $\bfz=(z_1, \dots z_{4g-4})$ of the quadratic differential we consider a polarization associated to the family of stable pointed curves $(\calX, \bfz)$. This is natural and necessary in our setting as the families of unpointed curves used to define the modified Hitchin base are not necessarily stable. It would be interesting to see how this compares to the moduli space of Hitchin pairs on stable curves and to the approach of Balaji, Barik, Nagaraj. We hope to come back to this question in future work.
\par
In contrast we emphasize that limits obtained by rescaling the Higgs fields
(for which e.g.\ the asymptotics of the hyperkähler metric was studied in
a series of works Mazzeo, Swoboda, Weiß, Witt and Fredrickson)
are not covered by the scope of this paper, since the starting point is a
compactification of the \emph{projectivization} of a space of quadratic differentials. 
\par
\medskip
\textbf{Acknowledgments.} This project started when both authors were in residence
at the MSRI, Berkeley, during the Fall 2019 semester, supported by NSF Grant
DMS-1440140. The authors acknowledge support by the DFG under Grant MO1884/2-1
and the Collaborative Research Centre TRR 326 {\it Geometry and
Arithmetic of Uniformized Structures}, project number 444845124.
We thank Martin Ulirsch for inspiring discussions and the referee for detailed
comments. We thank Lucia Caporaso, Sebastian Casalaina-Martin and Magarida Melo for kindly answering our questions.

\section{The \texorpdfstring{$\SL(2,\CC)$}{SL(2,C)}- and the
  \texorpdfstring{$\GL(2,\CC)$}{GL(2,C)}-Hitchin system}
\label{sec:HitchinBackground}

In this section we recall basic terminology for the Hitchin system
in the rank two case, that is for $\SL(2,\CC)$ for $\GL(2,\CC)$
and the trace-free variant of the latter. The results here, the description
of the Hitchin fibers and the BNR-correspondence (after \cite{BNR}), are
well-known and stated for comparison with our results over the modified
Hitchin base.
\par
\medskip
\paragraph{\bf $\GL(2,\CC)$-Higgs bundles}
Let $X$ be a Riemann surface of genus $g \geq 2$ and $K$ its canonical bundle.
A \emph{$\GL(2,\CC)$-Higgs bundle} is a pair $(E,\Phi)$ of a locally free
sheaf of rank~$2$  and a Higgs field $\Phi \in H^0(X,\End(E) \otimes K)$. A
$\GL(2,\CC)$-Higgs bundle is called \emph{(semi-)stable}, if for all
$\Phi$-invariant subbundles $L \subset E$
\be \label{eq:defstab}
\deg(L)< \deg(E)/2 \qquad (\deg(L) \leq \deg(E)/2)\,. 
\ee
Let $\mathrm{Higgs}^{\GL_2}(X)$ denote the moduli space of semistable
$\GL(2,\CC)$-Higgs bundles on $X$. This is an algebraic variety with the
moduli of stable Higgs bundles as a smooth subvariety of dimension $8g-6$
\cite{Nitsure}. The $\GL(2,\CC)$-Hitchin map
\be \label{eq:GLHitchinmap}
\Hit: \Higgs^{\GL_2}(X) \rightarrow \AA_2:=H^0(X,K) \oplus H^0(X,K^2), \quad (E,\Phi)
\mapsto (-\Tr(\phi), \det(\Phi))
\ee
is proper and surjective \cite{Hitorigin, Nitsure} and determines the
characteristic polynomial of the Higgs field.
\par
\medskip
\paragraph{\bf $\SL(2,\CC)$-Higgs bundles}
A \emph{$\SL(2,\CC)$-Higgs bundle} is a pair $(E,\Phi)$ as above, but
now we require that the determinant of~$E$ and
the trace of $\Phi$ are trivial. The moduli space $\mathrm{Higgs}^{\SL_2}(X)$ of
semistable $\SL(2,\CC)$-Higgs bundles on~$X$ is an algebraic variety.
It contains the moduli of stable Higgs bundles as a smooth subvariety of
dimension $6g-6$ \cite{Nitsure}. The $\SL(2,\CC)$-Hitchin map
\be \label{eq:Hitchinmap}
\Hit_{\SL_2}: \Higgs^{\SL_2}(X) \rightarrow \AA_2^\circ=H^0(X,K^2), \quad (E,\Phi) \mapsto
\det(\Phi) \ee
is again proper and surjective.
\par
\medskip
\paragraph{\bf Variants: Twisted Higgs bundles and/or fixed determinant}
We relax the conditions for Higgs bundles in two respects. An \emph{$M$-twisted
rank 2-Higgs bundle} is a pair $(E,\Phi)$ of a locally free sheaf $E$ of rank~$2$ and
a Higgs field $\Phi \in H^0(X, \End(E)\otimes M)$. Let $\Lambda$ be a line bundle on $X$. An \emph{$M$-twisted Higgs bundle
with fixed determinant~$\Lambda$} is a pair $(E,\Phi)$ of a locally free sheaf $E$ of
rank~$2$, such that $\det(E) = \Lambda$ and a Higgs field $\Phi
\in H^0(X, \End(E)\otimes M)$, such that $\mathrm{tr}(\Phi)=0$. That is, as in
the $\SL(2,\CC)$-case, we impose simultaneously the restriction on the determinant
of~$E$ and on the trace of~$\Phi$. The notion of (semi-)stability is still
defined using~\eqref{eq:defstab}.
\par
Let $\mathrm{Higgs}_{\Lambda}(X,M)$ denote the moduli space of semistable
$M$-twisted rank 2-Higgs bundle with fixed determinant $\Lambda$ on $X$.
This is an algebraic variety with the moduli of stable Higgs bundles as a
smooth subvariety \cite{Nitsure}. When $\deg(M) \geq 2g-2$ its dimension
is given by $3\deg(M)$. When $\deg(M) > \frac{g+1}{2}$ this holds for
generic $(X,M)$. Now the Hitchin map reads as
\[ \Hit_{\Lambda,M}: \Higgs_{\Lambda}(X,M) \rightarrow H^0(X,M^2), \quad (E,\Phi) \mapsto \det(\Phi)\,.
\] Here $\dim H^0(X,M^2)=2\deg M -g +1$ as long as $\deg M > g-1$. More generally, for generic $(X,M)$, we have 
\[ \dim H^0(X,M^2)=\left\{ \begin{array}{cc} 2\deg M -g +1 & \text{ for }\deg M  > \frac{g-1}{2} \\
0 & \text{ for } \deg M  \leq \frac{g-1}{2} 
\end{array} \right.
\] 
\par
\medskip
\paragraph{\bf The Hitchin system}
The cotangent bundle to the moduli space of stable rank 2 locally free sheaves
(with trivial determinant) is a dense subset of $\mathrm{Higgs}^{\GL_2}(X)$
(of $\mathrm{Higgs}^{\SL_2}(X)$) and its holomorphic symplectic structure
extends to the whole space. The Hitchin map restricted to pairs
$(q_1,q_2)$ with discriminant $\Delta = q_2-\frac{1}{4}q_1^2$ (resp.\ with~$q_2$) having only
simple zeros is an algebraically completely integrable system, called the
$\GL(2,\CC)$-Hitchin system (resp.\ the $\SL(2,\CC)$-Hitchin system).
In the twisted case: if $MK^{-1}$ has a section, the space $\Higgs_{\Lambda}(X,M)$
carries the structure of a Poisson manifold \cite{Botta,Markman}.
\par
\medskip
\paragraph{\bf The spectral curve}
Let $p_K: K \rightarrow X$ be the total space of the canonical bundle
and let $\eta: K \rightarrow p_K^*K$ be the tautological section. For
$q_i \in H^0(X,K^i)$, the \emph{spectral curve} $\Sigma$ is the zero divisor of
the characteristic section 
\be \label{eq:defeta}
\eta^2 + \eta \cdot p_K^*q_1 +  p_K^*q_2 \in H^0(K,p_K^*K^2)\,.
\ee
The spectral cover $\pi:=p_K|_{\Sigma}: \Sigma \rightarrow X$ is an analytic covering
factoring through the involution $\sigma: \Sigma \rightarrow \Sigma$ interchanging the
sheets. We also write abusively $\eta \in H^0(\Sigma, \pi^*K)$ for the restriction
to the spectral cover. We also frequently view the tautological section as a
section $\lambda \in H^0(\Sigma, K_\Sigma)$. The set of branch points is
$B = \mathrm{div}(\Delta) \subset X$. The section $\eta^\circ:=\eta + \frac12 \pi^*q_1 \in H^0(\Sigma,\pi^*K)$ vanishes on the ramification divisor $\wh{B}=\mathrm{div}(\eta^\circ)$ and is odd with respect to the involution, i.e.\ $\sigma^*\eta^\circ=-\eta^\circ$.
\par
The spectral curve is smooth if and only if $\Delta$ has simple zeros.  The
corresponding Hitchin fibers are torsors over an abelian variety. The following rank
two version of the BNR-correspon\-dence is the model for our correspondence statement.
See \cite{Schaub} for the general version of integral spectral curves,
including stability considerations.
\par
\begin{theorem}[\cite{Hitorigin, BNR}]\label{thm:Hitchinfiber}
Let $(q_1,q_2) \in H^0(X,K) \oplus H^0(X,K^2)$ be differentials such that
$\Delta$ has simple zeros. Then the fiber $\Hit^{-1}(q_1,q_2)$ of the Hitchin
map is a torsor over the Jacobian of the spectral cover.
\par
In the $\SL(2,\CC)$-case the fiber $\Hit_{\SL_2}^{-1}(q)$ over a quadratic differential
$q$ with simple zeros is a torsor over the Prym variety 
\be \label{eq:defPrym}
\ker(\Nm_{\pi}) \=\{ L \in \Pic^0(\Sigma) \mid L \otimes \sigma^*L =  \calO_\Sigma \}.
\ee
More generally, let $q \in H^0(X,M^2)$ be a section with simple zeros. Then
the fiber of the Hitchin map in the $M$-twisted case with determinant~$\Lambda$
is
\be \label{eq:defHitfiber}
\Hit_{\Lambda,M}^{-1}(q) \cong \{ L \in \Pic^{\deg M + \deg\Lambda}(\Sigma) \mid L \otimes \sigma^*L = \pi^*(M \otimes \Lambda)\,. \}
\ee
This is a torsor over the Prym variety $\ker(\Nm_{\pi})$ having dimension $g-1+\deg(M)$.
\end{theorem}
\par
\begin{proof} If we start with any line bundle~$L$ on $\Sigma$, then
$E = \pi_*L$ together with the $\pi$-pushforward
\[ \Phi = \pi_* \eta: \pi_*L \rightarrow \pi_*(L \otimes \pi^*K)=\pi_*L \otimes K
\]
of the multiplication map by~$\eta$ gives a Higgs bundle $(E,\Phi)$. By
Cayley-Hamilton and the definition of~$\eta$ in~\eqref{eq:defeta} we find 
that $(E,\Phi) \in \Hit^{-1}(q_1,q_2)$, see \cite[Proposition~3.6]{BNR}.
This Higgs bundle is indeed stable since a destabilizing subbundle would yield
a factorization of the characteristic polynomial (\cite[Lemma~3.2]{Schaub}),
contradicting the hypothesis that $\Delta$ resp.~$q$ have simple zeros.
\par
Conversely, let $(E,\Phi)$ be a Higgs bundle. The Higgs field, considered
as a map $E \otimes K^{-1} \to E$ makes $E$ into a $\pi_* \calO_\Sigma$-module,
since $\pi_* \calO_\Sigma = \calO_X \oplus K^{-1}$, and this is the same datum
as a line bundle on~$\Sigma$. Clearly the two constructions are inverse to
each other. This completes the first claim.
\par
We need another viewpoint that realizes the converse construction by computing the eigensheaf of the Higgs field pulled back to $\Sigma$. In a holomorphic chart~$(U,z)$ we write
$\eta^\circ = f dz$ with~$f$ holomorphic or with a simple zero. Then $E=\pi_*L$
is locally generated (as $\calO_X$-module) by $1$ and~$f$. 
The map
\bes
L|_U \to (\pi^* \pi_* L) \otimes \calO(\wh{B})|_U, \qquad
1 \mapsto 1 \otimes 1 +  f \otimes \frac{1}{f} 
\ees
is independent of the local frame chosen, injective and thus defines
an embedding. It exhibits
\be
L = \ker((\Phi - \eta \mathrm{Id}) \otimes \mathrm{Id}_{\calO(\wh{B})})
\subset (\pi^* \pi_* L) \otimes \calO(\wh{B})\,. \label{eq:eigensheaf}
\ee
\par
We now turn to the $\SL(2,\CC)$-case.
To prove the most general claim we will construct a morphism from the torsor
over the Prym variety to the Hitchin fiber.  Let $L \in \Pic(\Sigma)$ with
$L \otimes \sigma^*L = \pi^*(M \otimes \Lambda)$. As above let $E:=\pi_*L$.  
The equation (see \cite{hartshorne} Exercise IV.2.6)
\be \label{eq:detEstandard}
\det(E) = \Nm_{\pi}(L) \otimes \det(\pi_* \calO_{\Sigma}), 
\ee  pulls back to 
\be \label{eq:detEpullback}
\pi^* \det(E) \= L \otimes \sigma^* L \otimes \calO(\wh{B})^{-1}
\=\pi^*(M \otimes \Lambda)\otimes \calO(\wh{B})^{-1}\,, 
\ee
where $\wh{B} = \div(\eta) $ is the ramification divisor of $\pi$. Therefore in the
twisted context now $\calO(\wh{B})=\pi^*M$. The pullback $\pi^*: \Pic(X) \rightarrow \Pic(\Sigma)$ is injective as $\pi: \Sigma \rightarrow X$ is not unbranched.
Hence,  $\det(E)=\Lambda$. By the very definition of $\eta=\eta^\circ$, we have $\Tr(\Phi)=0$
and $\det(\Phi)=q$. Summing up, $(E,\Phi)$ defines a Higgs bundle
in $\Hit_{\Lambda,M}^{-1}(q)$. 
\par
Starting  with $(E,\Phi$) we first compute the two eigen-sheaves as in \eqref{eq:eigensheaf}
\[ L \= \ker((\pi^*\Phi-\eta \Id) \otimes \mathrm{Id}_{\calO(\wh{B})}), \quad
L'\=\ker((\pi^*\Phi+\eta \Id)\otimes \mathrm{Id}_{\calO(\wh{B})})\,.\]
As $\sigma^*\eta=-\eta$, we have $L'= \sigma^*L$. The inclusions
into $\pi^*E\otimes \calO(\wh{B})$ define an exact sequence
\[ 0 \rightarrow L \oplus \sigma^*L \rightarrow \pi^*E \otimes \calO(\wh{B}) \rightarrow \calO_{\wh{B}} \rightarrow 0
\] (see \cite[Theorem 5.5]{Horn}). Computing the determinant bundles yields
\[ L \otimes \sigma^*L = \pi^*(\Lambda)(\wh{B})= \pi^*(\Lambda \otimes M). \]
This is the Prym condition. 
\end{proof}
\par
\medskip
\paragraph{\bf Reducing to trace-free $\GL(2,\CC)$} In the following we will restrict our study to trace-free $\GL(2,\CC)$-Higgs bundles and $\SL(2,\CC)$-Higgs bundles. We will refer to a trace free $\GL(2,\CC)$-Higgs bundles as \emph{$\GL(2,\CC)^\circ$-Higgs bundles}. The smooth Hitchin fibers depend only on the discriminant $\Delta \in H^0(X,K^2)$ up to isomorphism. In general, we have an isomorphism 
\[ \Hit^{-1}(q_1,q_2) \rightarrow \Hit^{-1}(0,\Delta), \qquad (E,\Phi) \mapsto (E, \Phi - \tfrac{1}{2}\Id_E \otimes q_1 ).
\] In particular, the rational map $\AA_2 \to \overline{\calM}_{4g-3}$ to the moduli space of stable curves that associates to a point in the regular locus of $\AA_2$ the spectral curve is constant when only varying the abelian differential in $\AA_2$. This justifies the restriction to $\AA_2^\circ$ in terms of the spectral map.  
\par

\section{The compactification of strata of quadratic differentials}
\label{sec:comp}

We recall the construction of a smooth compactification of strata
of $k$-differentials from \cite{CMZ} and \cite{LMS}, specialized
to the case of quadratic differentials and later moreover to the principal
stratum where all zeros are simple. The fiber of this compactification
over a smooth curve is a modification on the space of quadratic differentials,
which we analyze in the subsequent Section~\ref{sec:Hitchinbase}. Finally
we compare this compactification with the incidence variety compactification
from \cite{kdiff}, the naive closure of strata in the Hodge bundle
over~$\barmoduli[g,n]$.
\par
Let  $\mu=(m_1,\ldots,m_n)$ be a \emph{type} of a quadratic differential, i.e.,
integers $m_i$ with $\sum m_i = 2(2g-2)$.
In the sequel we will be mainly interested in the case $\mu=(1^{4g-4})$
and we drop the argument~$\mu$ in that case. Let $\calQ_{g,n}(\mu)$ be
the moduli space of quadratic differentials~$(X,q)$ with~$n$ labeled
special points where~$q$ has a zero or pole of order~$m_i$ for
$i=1,\ldots,n$. We state the goal of the construction
and explain the missing notation subsequently.
\par
\begin{theorem} \label{thm:kLMS}
There exists a smooth Deligne-Mumford stack $\overline{\calQ}_{g,n}(\mu)$,
{\em the moduli space of quadratic \msds}, with the following properties.
\begin{itemize}
\item[i)] The space $\calQ_{g,n}(\mu)$ is dense in $\overline{\calQ}_{g,n}(\mu)$
\item[ii)] The boundary $D = \overline{\calQ}_{g,n}(\mu)
\smallsetminus \calQ_{g,n}(\mu)$ is a normal crossing divisor.
\item[iii)] The rescaling action of $\CC^*$ on $\calQ_{g,n}(\mu)$ extends
to $\overline{\calQ}_{g,n}(\mu)$ and the resulting projectivization
$\PP \overline{\calQ}_{g,n}(\mu)$ is a proper smooth stack.
\item[iv)] The space $\overline{\calQ}_{g,n}(\mu)$ is immersed 
in the compactification $\LMS[\wh{\mu}][\wh{g},\{\wh{n}\}]$ 
of a stratum of abelian differentials with partially labeled points.
\end{itemize}
\end{theorem}
\par
\paragraph{\textbf{Canonical double cover}}
The notion of {\em canonical double cover $\wh\pi: \wh{X} \to X$}
associated with a non-zero quadratic differential~$q$ on a smooth curve~$X$
is ubiquitous to the literature on half-translation surfaces. See
e.g.\ \cite[Section~2.1]{kdiff} for various methods of construction. One
possible definition of the canonical double cover is the normalization of
the spectral curve, i.e.\ $\wh{X} = \Sigma^\nu$. In particular the canonical
double cover is always smooth. It is irreducible if and only if~$q$ is
not the square of an abelian differential on~$X$. However, the genus
depends on the profile of the zeros on~$q$. To specify the type of the double
cover we let
\be \label{eq:typemuhat}
\whmu \,:=\, \Bigl(\underbrace{\wh m_1, \ldots, \wh m_1}_{\gcd(2,m_{1})},\,
\underbrace{\wh m_2,  \ldots, \wh m_2}_{\gcd(2,m_{2})} ,\ldots,\,
  \underbrace{\wh m_n, \ldots, \wh m_n}_{\gcd(2,m_{n})} \Bigr)\,,
\ee
where $\wh m_i := \tfrac{2+m_{i}}{\gcd(2,m_{i})}-1$. We let
$\wh{g} = g(\wX)$ and $\wh{n} = \sum_i \gcd(2,m_{i})$. Suppose that
there are $s_1$ points of odd order and $s_2$ points of even order,
$s_1+s_2=n$. Then by the cyclic covering constructions
(\cite[Section~3]{esnaultvielog}
or \cite[Section~2.1]{kdiff}) these quantities are related by
\be
\wh{g} \= 2g-1 + \frac{s_1}{2} \quad \text{and} \quad 
\wh{n} \= s_1 + 2s_2\,.
\ee
\par
\medskip
\paragraph{\textbf{Ordered vs.\ unordered points}} The point of departure
the space of quadratic differentials $\calQ_{g,n}(\mu)$ and the compactification
$\overline{\calQ}_{g,n}(\mu)$ come with~$n$ ordered marked points. The same
holds for the compactification
$\LMS[\wh{\mu}][\wh{g},\wh{n}]$ of the moduli space of abelian differentials
of type~$\wh{\mu}$. Given~$\mu$ (or~$\wh{\mu}$) there is a natural subgroup
of the symmetric group $S_n$ (or~$S_{\wh{n}}$) that acts by permuting points
with the same order. We indicate the corresponding quotients by the action of the symmetric group by brackets~$[n]$, e.g., $\overline{\calQ}_{g,[n]}(\mu)$. In the case of a signature
$\wh{\mu}$ arising from a double covering as in~\eqref{eq:typemuhat},
there is the subgroup of~$S_{\wh{n}}$ generated by the commuting involutions
swapping the pair of points stemming from the same even~$m_i$. The quotient
of $\LMS[\wh{\mu}][\wh{g},\wh{n}]$ by this involution is the space
$\LMS[\wh{\mu}][\wh{g},\{\wh{n}\}]$ appearing in item~iv) of
Theorem~\ref{thm:kLMS}.
\par
\medskip
\paragraph{\textbf{Multi-scale differentials}}
We now recall the definition of quadratic multi-scale differentials to
the extent we need it. The first piece of datum is a pointed stable
curve $(\wh{X},\wh{\bfz})$ of genus~$\wh{g}$, where the tuple $\wh{\bfz}$
consists of $\wh{n}$ non-singular points on~$\wh{X}$. Second, there is an
involution $\sigma: \wh{X} \to \wh{X}$ that fixes~$s_1$ points
and has $s_2$ orbits of length two on~$\wh{\bfz}$. We let 
$\pi: \wh{X} \to X$ be the corresponding quotient map,
and let $\bfz$ be the (ordered) $n=s_1+s_2$-tuple of images of~$\wh{\bfz}$.
\par
Third, we need a \emph{covering $\pi: \wh{\Gamma} \to \Gamma$ of enhanced
level graphs}. That is, $\wh{\Gamma}$ and $\Gamma$ are the dual graphs
of~$\wh{X}$ and~$X$ respectively, both are provided with a level structure,
and each edge carries an enhancement that we denote by $\widehat{\kappa}_{\wh{e}}
\in \ZZ_{\geq 0}$ for $\wh{e} \in \wh{\Gamma}$ and by $\kappa_e \in \ZZ_{\geq 0}$
for $e \in \Gamma$, subject to the constraints given below.
\par
\begin{figure}[ht]
	\bes	
	\begin{tikzpicture}[
	baseline={([yshift=-.5ex]current bounding box.center)},
	scale=2,very thick,
	bend angle=30,
	every loop/.style={very thick}]
	\node[circle, draw, inner sep=2pt, minimum size=12pt] (T) [] {$g$}; 
	\node[comp, fill] (T-1) [below=of T] {};
	\path (T) edge [shorten >=4pt] (T-1.center);	
	\node[comp,fill] (B) [below=of T] {}
	edge 
	node [order top right] {$2$}
	node [order middle right] {$4$}
	node [order bottom right] {$-6$} (T);
	\node [text height=12pt,below right] (B-2) at (B.south east) {$1$};
	\path (B) edge [shorten >=4pt] (B-2.center);
	\node [text height=12pt,below left] (B-2) at (B.south west) {$1$};
	\path (B) edge [shorten >=4pt] (B-2.center);
	\end{tikzpicture}
	\quad\,\,\,\, {~\overset{\widehat{\pi}}{\longleftarrow}~}
        \quad	\,\,\,\,  
	\begin{tikzpicture}[	 
	baseline={([yshift=-.5ex]current bounding box.center)},
	scale=2,very thick,
	bend angle=30,
	every loop/.style={very thick}]
	\node[circle, draw, inner sep=2pt, minimum size=12pt,
	      label=0:{$\widehat{g}-1$}] (T) [] {}; 
	\node[comp, fill] (T-1) [below=of T] {};
	\path (T) edge [white,shorten >=4pt] (T-1.center);	
	\node[comp,fill] (B) [below=of T] {}
	edge [bend left] 
	node [order top left] {$1$} 
	node [order middle left] {$2$}
	node [order bottom left] {$-3$} (T) 
	edge [bend right]		
	node [order top right] {$1$} 
	node [order middle right] {$2$}
	node [order bottom right] {$-3$} 
	(T);
	\node [text height=12pt,below right] (B-2) at (B.south east) {$2$};
	\path (B) edge [shorten >=4pt] (B-2.center);
	\node [text height=12pt,below left] (B-2) at (B.south west) {$2$};
	\path (B) edge [shorten >=4pt] (B-2.center);
	\end{tikzpicture}
        \qquad \qquad
        	\begin{tikzpicture}[
	baseline={([yshift=-.5ex]current bounding box.center)},
	scale=2,very thick,
	bend angle=30,
	every loop/.style={very thick}]
	\node[circle, draw, inner sep=2pt, minimum size=12pt] (T) [] {$g$}; 
	\node[comp, fill] (T-1) [below=of T] {};
	\path (T) edge [shorten >=4pt] (T-1.center);	
	\node[comp,fill] (B) [below=of T] {}
	edge 
	node [order top right] {$3$}
	node [order middle right] {$5$}
	node [order bottom right] {$-7$} (T);
	\node [minimum width=.8cm,text height=10pt,below left] (B-2) at (B.south west) {$1$};
	\path (B) edge [shorten >=6pt] (B-2.center);
	\node [text height=14pt,below] (B-2) at (B.south) {$1$};
	\path (B) edge [shorten >=2pt] (B-2.center);
	\node [minimum width=.8cm,text height=10pt,below right] (B-3) at (B.south east) {$1$};
	\path (B) edge [shorten >=6pt] (B-3.center);	
	\end{tikzpicture}
	\quad {~\overset{\widehat{\pi}}{\longleftarrow}~} \quad	 
	\begin{tikzpicture}[	 
	baseline={([yshift=-.5ex]current bounding box.center)},
	scale=2,very thick,
	bend angle=30,
	every loop/.style={very thick}]
	\node[circle, draw, inner sep=2pt, minimum size=12pt,
          label=0:{$\widehat{g}-1$}] (T) [] {}; 
	\path (T) edge [shorten >=4pt] (T-1.center);	
	\node[circle, draw, inner sep=0pt, minimum size=12pt] (B) [below=of T] []{$1$}
	edge 		
	node [order top right] {$4$} 
	node [order middle right] {$5$}
	node [order bottom right] {$-6$} 
	(T);
	\node [minimum width=.8cm,text height=10pt,below left] (B-2) at (B.south west) {$2$};
	\path (B) edge [shorten >=6pt] (B-2.center);
	\node [text height=14pt,below] (B-2) at (B.south) {$2$};
	\path (B) edge [shorten >=2pt] (B-2.center);
	\node [minimum width=.8cm,text height=10pt,below right] (B-3) at (B.south east) {$2$};
	\path (B) edge [shorten >=6pt] (B-3.center);
	\end{tikzpicture}	
	\ees 
        \caption{Level graphs of two (left) resp. three (right) zeros coming together and all other zeros simple} \label{cap:23together}
\end{figure}
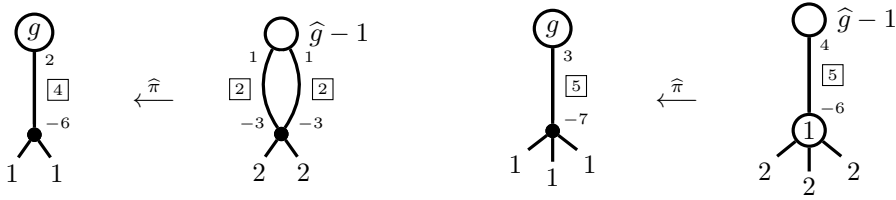
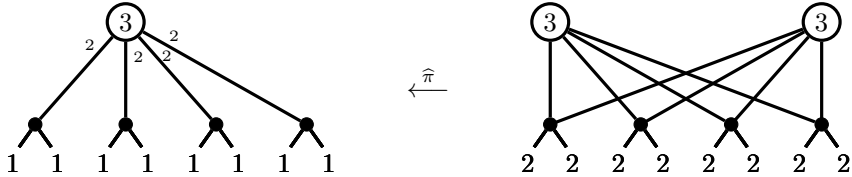
\begin{figure}[ht]
	\bes		
	\begin{tikzpicture}[
	baseline={([yshift=-.5ex]current bounding box.center)},
	scale=2,very thick,
	bend angle=30,
	every loop/.style={very thick}]
	\node[circle, draw, inner sep=2pt, minimum size=12pt] (T) [] {$3$}; 
	\node[comp,fill] (B) [below=of T] {} edge 
		node [order top right,xshift=-1pt,yshift=-3pt] {$2$} (T); 
	\node[comp,fill] (B-1) [left=of B]  {} edge 
		node [order top left] {$2$} (T);
	\node[comp,fill] (B-2) [right=of B] {} edge 
		node [order top right,xshift=2pt,yshift=-4pt] {$2$} (T);	
	\node[comp,fill] (B-3) [right=of B-2] {} edge node [order top right,yshift=2pt] {$2$} (T);			
	\foreach \x in {B,B-1,B-2,B-3} 
	\foreach \y in {B1,B2,B3,B4} 
	  { 
	  	\node [text height=12pt,below right] (\y) at (\x.south east) {$1$};
	  	\path (\x) edge [shorten >=4pt] (\y.center);
	  	\node [text height=12pt,below left] (\y) at (\x.south west) {$1$};
	  	\path (\x) edge [shorten >=4pt] (\y.center);
	  } 	
	\end{tikzpicture}
	\qquad {~\overset{\widehat{\pi}}{\longleftarrow}~} \qquad	 
	\begin{tikzpicture}[	 
	baseline={([yshift=-.5ex]current bounding box.center)},
	scale=2,very thick,
	bend angle=30,
	every loop/.style={very thick}]
	\node[circle, draw, inner sep=2pt, minimum size=12pt] (T) [] {$3$}; 	
	\node[comp,fill] (B) [below=of T] {} edge node [] {} (T); 
	\node[comp,fill] (B-1) [right=of B]  {} edge node [] {} (T);
	\node[comp,fill] (B-2) [right=of B-1] {} edge node [] {} (T);	
	\node[comp,fill] (B-3) [right=of B-2] {} edge node [] {} (T);	
	\node[circle, draw, inner sep=2pt, minimum size=12pt] (T-1) [above=of B-3] {$3$};			
	\foreach \x in {B,B-1,B-2,B-3} 
	\foreach \y in {B1,B2,B3,B4} 
	{ 
		\node [text height=12pt,below right] (\y) at (\x.south east) {$2$};
		\path (\x) edge [shorten >=4pt] (\y.center);
		\node [text height=12pt,below left] (\y) at (\x.south west) {$2$};
		\path (\x) edge [shorten >=4pt] (\y.center);
	} 
	\path (T-1) edge [] (B.center) 
	            edge [] (B-1.center)
	            edge [] (B-2.center)	            
	            edge [] (B-3.center);
	\end{tikzpicture}	
	\ees
        \caption{Level graphs of the locus where the quadratic differential
   in $\calQ_{3,4}(1^8)$  becomes a global square with zeros coming together pairwise} \label{cap:qtosquare}
\end{figure}

A level structure is by definition a weak total order on the vertices,
any two vertices can be compared and equality is permitted.
The top level is usually normalized to
be level zero and depicted on top, see the examples in
Figure~\ref{cap:23together} and Figure~\ref{cap:qtosquare}. The lower levels
are then referred to by~$-1$, $-2$ etc. The covering of graphs, abusively
also called~$\pi$, is level-preserving and the quotient map induced by
the involution~$\sigma$. The level structure on the graph allows to call
edges \emph{horizontal} if
they start and end at the same level, and \emph{vertical} otherwise.
\par
The enhancement~$\kappa_e$ of an edge~$e$ of~$\Gamma$ is even if and only
if~$e$ has two preimages in~$\wh{\Gamma}$. In this case
$\wh{\kappa}_{\wh{e}_j} = \kappa_e/2$ for either of the two preimages~$\wh{e}_j$
of~$e$. If $\kappa_e$ is odd, then $\wh{\kappa}_{\wh{e}} = \kappa_e$.
Moreover, the enhancement~$\kappa_e$ is zero if and
only if the edge~$e$ is horizontal. These conditions reflects
the branching rule of canonical covers, see below.
Boundary strata of $\overline{\calQ}_{g,[n]}(\mu)$ are labeled by the
coverings $\pi: \wh{\Gamma} \to \Gamma$, or usually just by~$\wh{\Gamma}$.
Within the next item we clarify which finite number of enhanced level
graphs~$\wh{\Gamma}$ actually occur as labels of a boundary stratum.
\par
Forth, the core datum of a multi-scale differential is a \emph{a twisted
differential}. This is collection of differentials $\bflambda =
(\lambda_v)_{v \in V(\wh{\Gamma})}$, subject to the following conditions. At the
marked legs, the differential has a zero  of order~$\wh{m}_i$ (or a pole,
depending on the sign of $\wh{m}_i$). At the upper end of a vertical edge~$e$, there
is a zero of order $\wh{\kappa}_e -1$. At the lower end of~$e$
or at horizontal edges there is a pole of order $-\wh{\kappa}_e-1$, with
matching residues along horizontal edges. Finally the residues at the
lower ends of the edges are constrained by a global residue condition (GRC).
See \cite[Section~3]{strata} for a discussion and examples of those conditions.
\par
Fifth, there is a collection of \emph{prong-matchings~$\bftau$},
one for each vertical edge, denoted by~$\tau_e$. A prong-matchings is a
bijection of the outgoing horizontal prongs of the differential at the upper
end of~$e$ with the incoming horizontal prongs of the differential at the lower
end of~$e$, reversing the cyclic order on prongs from the planar structure.
Prong matchings determine the branches when smoothing a stable curve with
a twisted differential. Prong matchings play virtually no role in this
paper, but recording them is necessary to construct the smooth stack in
Theorem~\ref{thm:kLMS}.
See \cite[Section~5]{LMS} for details of the definition and \cite{strataEC}
for examples of the combinatorics of these objects.
\par
So far the tuple of objects $(\wh{X}, \wh{\bfz}, \wh{\Gamma}, \bflambda,
\bftau)$ without the conditions on the existence of involutions or the
map~$\pi$ defines an (ordinary) multi-scale differential, i.e., points
in the space $\LMS[\wh{\mu}][\wh{g},\wh{n}]$. The last point, and the
reason for the name, is to define the equivalence relation.
For this purpose, we group the differentials $\lambda_v$ for all vertices~$v$
on the same level~$-i$ to a tuple $\lambda_i$. Consequently we
write $\bflambda = (\lambda_i)_{i \in L(\wh{\Gamma})}$, where
$L(\wh{\Gamma}) = L(\Gamma)$ is the set of levels of these two graphs.
The group~$\CC^*$ rescales such
a tuple $\lambda_i$ simultaneously. We apply this for all levels except for
the top level (which is rescaled in projectivization, see item~iii) of the
theorem). Suppose~$\Gamma$ has~$L$ levels below zero. A cover of the resulting
torus~$(\CC^*)^L$ acts by rescaling $\lambda_i$ and rotating the prong-matching
simultaneously. The product is called the level rotation torus in
\cite[Section~6]{LMS} and defines the equivalence relation we need.
\par
In the last step, we characterize $\overline{\calQ}_{g,n}(\mu)$
locally (in the domain) as subspace of $\LMS[\wh{\mu}][\wh{g},\{\wh{n}\}]$.
For this, we need in addition to the
(ordinary) multi-scale differential the involution~$\sigma$ and the quotient
map~$\pi$, as above. Finally, we require that the collection~$\bflambda$ is
anti-invariant under~$\sigma$, i.e. $\sigma^* \bflambda = - \bflambda$.
This implies that there is a collection $\bfq = (q_i)_{i \in L(\Gamma)}$
of quadratic multi-scale differentials on the subcurves~$X_i$ of level~$-i$ with simple zeroes at the marked points $\bfz:=\wh{\pi}(\wh{\bfz})$ such
that $\lambda_i^2 = \wh{\pi}^* q_i$. To sum up, on the pointed stable curve $(X,\bfz)$ the data of a quadratic multi-scale differential is the tuple $(\Gamma,\bfq,\bftau)$.
\par
\medskip
\paragraph{\textbf{Universal family}} So far, we have given multi-scale
differentials pointwise. The complete definition for families is lengthy as
it requires to encode the passage from the smooth situation (where no
prong-matching is present) to boundary points (where prong-matchings are
part of the datum). The notion of rescaling ensemble in \cite{LMS} serves
this purpose. We only give the consequences here.
\par
By construction, $\overline{\calQ}_{g,n}(\mu)$ comes with a universal family
$\wh{f} \colon \wh{\calX} \to \overline{\calQ}_{g,n}(\mu)$ (all to be interpreted
in the sense of stacks, i.e.\ on an \'etale chart), a universal family
$f: \calX \to \overline{\calQ}_{g,n}(\mu)$ and the universal double covering
$\wh{\pi}: \wh{\calX} \to \calX$.  There are also the universal families
of multi-scale differentials that we continue to denote $\bflambda$ and~$q$,
as we did for single fibers. 
\par
\medskip
\paragraph{\textbf{Examples}} We give three examples of graphs
with two levels and without horizontal edges. These correspond thus to
boundary divisors in $\overline{\calQ}_{g,n}(\mu)$, i.e., we moreover specialize
to $\mu=(1^{4g-4})$. In these graphs a dot is a vertex of genus zero,
otherwise the genus is written inside or next to the vertex. The number in the
square is the enhancement. The other numbers written next to half-edges
are the orders of zeros of poles that the (abelian or quadratic) multi-scale
differential is required to have. In all the cases the graphs~$\Gamma$ are given together with the
double cover $\wh{\pi}: \wh{\Gamma} \to \Gamma$.
\par
\par
\par
\medskip
\paragraph{\textbf{The covering viewpoint}} Since the total space~$\calX$
of the universal family over $\overline{\calQ}_{g,n}(1^{4g-4})$ is not smooth,
the double cover~$\wh{\pi}$ is not given ``as usual'' in the context of
spectral curves by a section of square of a line bundle. However, for
a fixed point of $\overline{\calQ}_{g,n}(1^{4g-4})$, corresponding to a boundary
stratum labeled by $\wh{\Gamma}$, or more generally for an equisingular
deformation in such a stratum, the cover restricted to each level
is given in this usual way. We determine here the line bundle and
the section in terms of graph data.
\par
Consider the bundle on the level-$i$-subcurve given by
\be \label{eq:Micoverbundle}
M_i \= \omega_X|_{X_i}(K_i), \qquad
K_i \= \Bigl(-\sum\limits_{e^- \in X_i} \lfloor \frac{\kappa_{e}}{2}
\rfloor e^-
+ \sum\limits_{e^+ \in X_i} \lceil \frac{\kappa_{e}}{2}\rceil e^+\Bigr)\,,
\ee
where we sum over all preimages $e^\pm \in X_i$ of the nodes~$e$ that connect
to upper/lower levels with respect to level~$i$. 
Then $q_i \in H^0(X_i,M_i^{\otimes 2})$ with only simple zeros,
namely at the marked points in~$X_i$ and at the nodes where~$\kappa_e$
is odd, by definition of compatibility with the level structure.
Consequently, the covering
\be \label{eq:Miqicover}
\wh{\pi}_i: \wh{X}_i \to X_i
\ee
is the spectral curve for the pair $(M_i,q_i)$.
\par
\medskip
\paragraph{\textbf{The incidence variety compactification (IVC)}}
Instead of the space of multi-scale differentials $\overline{\calQ}_{g,n}(\mu)$
one can consider the compactification $\overline{\calQ}^\IVC_{g,n}(\mu)$
of the stratum ${\calQ}_{g,n}(\mu)$ inside the total space of the $2$nd-Hodge bundle $\Omega^2 \barmoduli[g,n]$,
the space of pointed stable quadratic differentials. There is a forgetful map
\bes
f_\IVC: \overline{\calQ}_{g,n}(\mu) \to \overline{\calQ}^\IVC_{g,n}(\mu)
\ees
that takes a multi-scale differential $(\pi: \wh{X} \to X, \wh{\bfz},
\wh{\Gamma}, \bflambda, \bftau)$ and
\begin{itemize}
\item forgets the prong-matching~$\bftau$
\item forgets the covering surface~$\wh{X}$,  its marked points~$\wh{\bfz}$,
the abelian differentials~$\bflambda$ the graph covering, just
retaining $(X,\bfz, \Gamma,\bfq)$, 
\item forgets the enhancements on~$\Gamma$, but keeps the level structure
\item allows rescaling the lower level components of~$\bfq$ individually on
  each irreducible component of the stable curve.
\end{itemize}
Said differently, points in $\overline{\calQ}^\IVC_{g,n}(\mu)$ are given
by a pointed stable curve $(X,\bfz,\Gamma,\bfq)$ with a non-enhanced level
graph and a twisted quadratic differential compatible with~$\Gamma$.
The latter is defined by the existence of a double covering such that
the pullback is an (abelian) differential in the above sense.
Obviously the map $f_{\IVC}$ is equivariant under the rescaling actions
and defines a map on the projectivizations.
\par
Passing to $\overline{\calQ}^\IVC_{g,n}(\mu)$ we lose smoothness
(singularities are e.g.\ not $\QQ$-factorial), the normal crossing
boundary divisor, the interpretation as a functor of multi-scale differentials,
but retain an object whose points are more easily described.

\section{Resolutions of the spectral map}
\label{sec:spectralmap}

In this section we outline the definitions of the spectral map and its variants
for the $\GL(k,\CC)$-Hitchin base. For the case $k=2$ we give a characterization
of the image and the points in the minimal resolution. 
\par
Let~$X_\st$ be a smooth curve and $\AA_k = \oplus_{i=1}^k H^0(X_\st,K_{X_\st}^{\otimes i})$
be the $\GL(k,\CC)$-Hitchin base. To a collection $\bfq = (q_1,\ldots,q_k)$
of differentials, i.e.\ to a point in~$\AA_k$ there is the corresponding
spectral curve $\Sigma_\bfq$. We let $\Delta_\bfq$ be the
discriminant of~$\bfq$ and $\Delta_\AA = \{\Delta_\bfq = 0\} \subset \AA_k$
the discriminant locus. 
\par
For a generic choice of~$\bfq$ the $k(k-1)$-differential $\Delta_\bfq$ has simple
zeros, thus in fact $n=k(k-1)(2g-2)$ of them. In this case the spectral curve
$\Sigma_\bfq$ is a smooth curve of genus~$\wh{g}$. Associating 
to $(X_\st,\bfq)$ the spectral curve $\Sigma_\bfq$ thus defines a rational map
\bes
\psi_{X_\st,k}: \AA_k \dashrightarrow \barmoduli[\wh{g}]
\ees
that we call the \emph{spectral map}. As every rational map, this map can
be resolved by taking the closure~$\calR^{X_\st}_{g,k}$ of the graph of $\psi_{X_\st,k}$
inside $\AA_k \times \barmoduli[\wh{g}]$. We denote the resolution by
\be
\widetilde{\psi}_{X_\st,k}: \calR^{X_\st}_{g,k} \to \barmoduli[\wh{g}]\,.
\ee
The general goal is a characterization of the points in $\calR^{X_\st}_{g,k}$
or at least of the stable curves that arise as their
$\widetilde{\psi}_{X_\st,k}$-images.
\par
All these definitions have analogs in the trace-free case, where the
abelian differentials are zero. That is, on the $\SL(k,\CC)$- (or trace-free)
Hitchin base $\AA_k^\circ = \oplus_{i=2}^k H^0(X_\st,K_{X_\st}^{\otimes i})$ there is
the rational map $\psi^\circ_{X_\st,k}: \AA_k^\circ \dashrightarrow \barmoduli[\wh{g}]$,
whose resolution we denote by $\calR^{X_\st,\circ}_{g,k}$.
\par
We give an explicit description of the points in resolutions~$\calR^{X_\st}_{g,k}$
and~$\calR^{X_\st,\circ}_{g,k}$ for the case $k=2$. The answer is very close to the
compactification of \cite{strata} and \cite{kdiff} and requires mainly
to keep track of rational tails correctly.
\par
\begin{theorem} \label{thm:charrank2}
A stable curve~$\widetilde{X}$ of genus $\wh{g} = 4g-3$ is in the image of
the $\SL(2,\CC)$-spectral map~$\widetilde{\psi}^\circ_{X_\st}$ if and only if
\begin{itemize}
\item there exists a level structure on the dual graph~$\Gamma$ of~$\widetilde{X}$,
\item there is an involution~$\sigma$ on~$\widetilde{X}$ such that the quotient
of the top level curve by~$\sigma$ is isomorphic to~$X_\st$,
\item each $\sigma$-fixed self-node of~$\widetilde{X}$ can be replaced by a rational bridge so that the resulting semistable curve~$\wh{X}$ has
a $\sigma$-anti-invariant twisted differential of type $\mu=2^{4g-4}$ compatible with
one (equivalently: any) level structure on the dual graph of~$\wh{X}$ extending
the one on~$\Gamma$.
\end{itemize}
In particular, the top level of~$\Gamma$ has a unique vertex or two vertices
exchanged by~$\sigma$.
\par
The $\GL(2,\CC)$-spectral map~$\widetilde{\psi}_{X_\st}$ has the same image.
\end{theorem}
\par
We use two variants of the spectral map for the proof. First, the spectral
maps~$\psi_{X_\st,k}$ can be assembled to a rational map~$\psi_k$ from the
fiber product of bundles $\iomoduli[g]$ of $i$-differentials
($1 \leq i \leq k$) over $\moduli[g]$. Second, the target point
of the map $\psi_{X_\st,k}$ only depends on the class of the point $\AA_k$ up to
the action of $\CC^*$ rescaling the $i$-th component with exponent
$i$. We denote this map by
\bes
\psi_k^{\PP} : \Bigl(\times_{\moduli[g]} \iomoduli[g] \Bigr) / \CC^* \dashrightarrow
\barmoduli[\wh{g}]
\ees
and its resolution by $\calR_{g,k}^{\PP}$. All these objects obviously have their
trace-free variants, decorated by a superscript~$\circ$.
\par
\medskip
\paragraph{\bf Interpolation by Hurwitz spaces} Let $k=2$ from now on. By definition
the space~$\calR^X_{g,2}$ is closely related to the incidence variety compactification, or rather
the quotient $\overline{\calQ}_g^\IVC := \overline{\calQ}^\IVC_{g,4g-4}(1^{4g-4})/
S_{4g-4}$ by the permutation group of the marked points. The following
space of admissible covers interpolates between the two, dominating both birationally
but with neither map being the identity.
\par
We take $\ol{\calH}_g$ to be the admissible cover compactification of the Hurwitz space
of degree~$2$ covers $\pi: \wh{X} \to X$ with simple branching over~$4g-4$ (unmarked)
points with $g(X) =g$ (and thus $g(\wh{X}) = \wh{g}$). The Hurwitz space comes with source and target maps
\bes \ol{\calH}_g \to \barmoduli[\wh{g}], \qquad \ol{\calH}_g \to \barmoduli[g].
\ees We consider the subspace
\bes
\widetilde{\calH}^\circ_g  \,\subset \, \ol{\calH}_g\,,
\ees
defined as closure of the locus  where~$\pi$ is a covering between smooth curves
and the branch divisor in~$X$ is the zero locus of a quadratic differential. By definition there are rational maps, 
\be \widetilde{\calH}^\circ_g \dashrightarrow \Omega^2 \moduli[g]/ \CC^*, \qquad \widetilde{\calH}^\circ_g \dashrightarrow \Omega^2 \barmoduli[g,4g-4] / \CC^*. \label{equ:augm_target}
\ee
These define forgetful maps
\be
f_1 : \wt{\calH}^\circ_g \to \calR^{\PP,\circ}_{g,2} \quad \text{and} \quad
f_2 : \wt{\calH}^\circ_g \to \PP \overline{\calQ}_g^\IVC\,,
\ee
Here the map~$f_1$ is induced by taking product of the first map in \eqref{equ:augm_target} and the source map of the Hurwitz space. Recall that we can think of $\calR^{\PP,\circ}_{g,2}$ as the closure of the graph of $\psi_k^{\PP}$. 
The map~$f_2$ is induced by the second map in \eqref{equ:augm_target}.
They are well-defined since both domain and range are defined as closures
of loci on which the map is obviously well-defined. The theorem will follow
from the following statements, where the first is a direct consequence of the
definitions.
\par
\begin{lemma} The map $f_2$ is quasi-finite, surjective and generically bijective.
The cardinality of the $f_2$-fiber over $(X,\bfz,\Gamma,\bfq)$ is precisely the
number of double coverings $\pi: \wh{\Gamma} \to \Gamma$ of enhanced level
graphs with given target graph.
\end{lemma}
\par
The passage between the two level graphs in the $f_2$-fibers is called
'criss-cross' in \cite{kdiff}, see Example~4.3.
\par

\begin{lemma}\label{lem:admissible_reduction}
  Consider a point in $\calR^{\PP,\circ}_{g,2}$ represented by a triple $(\widetilde{X}, X_\st, q_2)$. There exists
a unique semistable model $\wh{X} \to \widetilde{X}$, which has an admissible
cover $\pi: \wh{X} \to X$, such that the stabilization of $X$ is $X_\st$.
\end{lemma}
\begin{proof}
  Using admissible reduction applied to the germ of a family with special fiber $(\widetilde{X}, X_\st, q_2)$ we obtain a semi-stable model $\wh{X} \to \widetilde{X}$ and an admissible cover $\pi: \wh{X} \to X$. By the process of admissible reduction $X$ is given by $X_\st$ augmented by rational tails. On all irreducible components of $\wh{X}$ of genus larger than one the
('hyperelliptic') involution~$\sigma$ is uniquely determined. Moreover, there is at least
one irreducible component ('top level') where~$q_2$ is non-zero by definition of
where we analyze~$f_1$. On this component the differential~$q_2$ determines
the double cover and the involution~$\sigma$. The argument of uniqueness for
elliptic tails of $\wh{X}$ is now similar, using the point of attachment (or exchange pair of
points of attachment) to a higher genus curve or 'higher level curve'. 
\end{proof}

\begin{lemma}
The map $f_1$ is surjective, and injective over $\{q_2 \neq 0\}$ if~$\wh{X}$ defined by Lemma \ref{lem:admissible_reduction} has no
rational tails. However, the fibers of~$f_1$ are positive-dimensional in general
even over $\{q_2 \neq 0\}$.
\end{lemma}
\par
\begin{proof}
The surjectivity follows from properness of the space of admissible covers.
\par
By Lemma \ref{lem:admissible_reduction} the non-uniqueness stems precisely from a chain of rational tails
that is contracted under the stabilization $\wh{X} \to \widetilde{X}$.
For example take $X$ to be a genus two curve~$X_0=X_\st$, attach a rational
tail $T_1$ with two branch points and to that another tail $T_2$ with
two branch points. Then the double cover $\wh{X} \to X$ is unramified
over the nodes and the stable model~$\widetilde{X}$ of~$\wh{X}$ has two irreducible components: A unramified double cover of $X_0$ augmented by a singular elliptic curve meeting the top level component in two nodes. (One might refer to this singular elliptic curve as fish bridge.) However the admissible cover records the
cross-ratio of the four special points on~$T_1$, while the image in $\calR_{2,2}^{\PP}$
loses that information.
\end{proof}
\par
\begin{proof}[Proof of Theorem~\ref{thm:charrank2}]
This is a restatement of the surjectivity of the map~$f_1$ together
with the 'almost bijectivity' of~$f_2$ and the main theorem of \cite[Theorem 1.5]{kdiff}
that characterizes boundary points in $\overline{\calQ}_2^{\mathrm{IVC}}$,
hence in $\calR^{X_\st,\circ}_{g,2}$.
\par
The statement about the $\GL(2,\CC)$-spectral map~$\widetilde{\psi}_{X,2}$ follows
from the observation that for a dense open subset of $\AA_2$ the spectral double
cover is given by a quadratic differential, the discriminant $\Delta_\bfq$.
Since on this locus the spectral maps have the same image, the same holds for
the closure.
\end{proof}
\par

\section{The modified Hitchin base}
\label{sec:Hitchinbase}

The target of the usual $\SL(2,\CC)$-Hitchin map given in~\eqref{eq:Hitchinmap}
is the space $\AA_2^\circ = H^0(X_\st,K^2)$ of quadratic differentials on a fixed smooth
curve~$X_\st$. We refer to this vector space as the \emph{Hitchin base}. Both
in the compactification and the resolution of the spectral map the fibers of
the forgetful map $\psi: \overline{\calQ}_{g,[4g-4]}(1^{4g-4}) \to \barmoduli[g]$
play an important role. We call this fiber $\BB_{X_\st} = \psi^{-1}([X_\st])$ the
\emph{modified Hitchin base} and describe the geometry of this DM-stack
in more detail.
\par
We describe the enhanced level graphs~$\Gamma$ that appear in $\BB_{X_\st}$
if $X_\st$ is a smooth curve: 
\begin{itemize} \item[i)] The dual graph $\Gamma$ is a tree.
\item[ii)] The top level is has a unique vertex corresponding to
the irreducible curve $X_\st$ of genus~$g$.
\item[iii)] There are no horizontal edges in $\Gamma$
(and thus in $\wh{\Gamma}$).
\end{itemize}
\par
\begin{proposition} \label{prop:MHB1}
There is a (``forgetful'') morphism $b: \BB_{X_\st} \to H^0(X_\st,K^2)$. This map is
birational and it is an isomorphism over the locus where there is at
most one non-simple zero which is of order two.
\end{proposition}
\par
\begin{proof} The morphism $b$ is defined by assigning to a quadratic multi-scale differential in $\BB_{X_\st}$ the quadratic differential on top level $X_0=X_\st$. If the zeros are pairwise disjoint, then the construction of
$\overline{\calQ}_{g,[4g-4]}(1^{4g-4})$ does not modify the differential.
\par
Whenever there is only one double zero, the construction of 
$\overline{\calQ}_{g,[4g-4]}(1^{4g-4})$ replaces the smooth curve with double
zeros by a rational tail with two zeros as in Figure~\ref{cap:qtosquare}
left. Since these rational tails have three
(unordered) marked points and thus have no moduli, the map $b$ is still an isomorphism over
this locus. 
\end{proof}
\par
\begin{remark} \textrm{
Since the level graphs~$\Gamma$ parameterizing boundary strata of the
modified Hitchin base arise from
the collision of points (only), they are always of compact type, i.e.
$h^1(\Gamma)=0$.}
\end{remark}
\par
\begin{example} \label{ex:g2Hitbase} 
  The case $g=2$ is particularly simple since the (second order)
Weierstrass sequences are the same for all curves, and yet this case shows that
$b$ is not birational beyond the locus claimed in
Proposition~\ref{prop:MHB1}. 
\par
Consider $\PP^2 = \PP(H^0(X_\st,K^{\otimes 2}))$ for any~$X_\st$ with $g(X_\st)=2$.
A basis of one-forms is $\{\omega_1 = dx/y, \omega_2 =xdx/y\}$, and
a basis of two-forms is $\{\omega_1^2, \omega_1 \omega_2, \omega_2^2\}$. In particular, all quadratic differentials are invariant under the hyperelliptic involution. 
This implies that the locus of differentials of type~$(2,1,1)$ consists of
six lines (called~$W_i$ in Figure~\ref{fig:MHB} that shows only $i=1,2$),
one for each of the 6~Weierstrass points. Note that the Weierstrass points of the first and the second order agree. The locus where the
zeros are of type $(2,2)$ are the pairwise intersection points together
with the reducible locus ($g=2$--version of Figure~\ref{cap:qtosquare}).
This reducible locus is a curve~$V \cong \PP^1$, embedded via the
Veronese embedding
into the $\PP^2$. The double zeros along the reducible locus are never
Weierstrass points, except when they collide to form a four-fold
zero. This gives a special point on each of the Weierstrass lines.
There are five more special points on each of them, the intersection with the
other Weierstrass lines. This gives $6\cdot 5/2+6$ special points in total.
\par
We claim that $b: \PP\BB_{X_\st} \to \PP^2 $ is the blowup of these 15 points $W_i \cap W_j$
and a $\PP^1$ with an orbifold point of order two blown into the intersection
points $W_i \cap V$ for each $i=1,\ldots,6$. For the latter note that
this intersection is tangent of order two. A first blowup makes the
intersection transversal, with exceptional divisor~$E_i$. A second blowup
creates the curve $F_i$ as exceptional divisors, which is geometrically
explained below. After this blowup, $E_i$ becomes a $(-2)$-curve, whose
contraction lead to an ordinary double point, equivalently the orbifold
point of order two.
\par
We now justify the decoration in Figure~\ref{cap:qtosquare}. Each double zero
produces an $\moduli[0,3]$--rational tail. Two
double zeros thus give a ``cherry''-type level graph. We now have to
distinguish cases. Suppose the double zeros are at two Weierstrass
points. This implies that the quadratic differential is the product
of two abelian double zero differentials, hence not a square.
The relative scale of the one-forms on the two rational tails is the extra
parameter that induces a $\PP^1$ instead of a point in the original
Hitchin base. This gives the curves $D_{ij}$ in Figure~\ref{fig:MHB}
for $1 \leq i < j \leq 6$. These $D_{ij}$ intersect~$W_i$ and~$W_j$
at points where the cherry is slanted left or right.
\par
On the other hand, if the double zeros are not Weierstrass points and
hence hyperelliptic conjugates, then we are on the Veronese curve, the
reducible locus. There, the two differentials on the lower level have to be
of the same scale in order to allow a deformation to proper quadratic
differentials. This can be seen on the double cover ($g=2$--version of
Figure~\ref{cap:qtosquare}, right), where the global residue condition
forces all residues to be the same. Also a global dimension count shows
that scaling differentials on lower level independently would give a
boundary stratum of the same dimension as the stratum 
${\calQ}_{g,[4g-4]}(1^{4g-4})$.
\par
\begin{figure}[ht]
	\bes	
	\begin{tikzpicture}[transform shape,
	baseline={([yshift=-.5ex]current bounding box.center)},
	scale=1, very thick,
	node distance=.4cm, 
	bend angle=30,
	every loop/.style={very thick}]
	\draw[name path=LV] (-.8,1.55) -- +(0:9.5)  node (V) [right] {$V$};
	\draw[name path=LF1] (-.3,2) -- +(270:6.9) node (F1) [right,midway] {$F_1$};
	\draw[name path=LF2] (2,2) -- +(270:3) node (F2) [right,midway] {$F_2$};
	\draw[name path=LW2] (1.5,-.45) -- +(0:5.25) node (W2) [right] {$W_2$};
	\draw[name path=LD12] (5.5,0) -- +(270:4.9) node (D12) [left,midway] {$D_{12}$};
	\draw[name path=LW1] (-.8,-4.4) -- +(0:9)  node (W1) [right] {$W_1$};
	\fill[name intersections={of=LV and LF1}] (intersection-1) circle (2.8pt);
	\fill[name intersections={of=LF1 and LW1}] (intersection-1) circle (2.8pt);
	\fill[name intersections={of=LW1 and LD12}] (intersection-1) circle (2.8pt);
	\fill[name intersections={of=LV and LF2}] (intersection-1) circle (2.8pt);
	\fill[name intersections={of=LF2 and LW2}] (intersection-1) circle (2.8pt);
	\fill[name intersections={of=LW2 and LD12}] (intersection-1) circle (2.8pt);
	\fill[name intersections={of=LF1 and LW1}] (intersection-1) circle (2.8pt);
	\begin{scope}[] 
	\node[comp,fill] (T) [] {};
	\node[comp,fill] (T-1) [right=of T] {};	
	\draw (T) -- +(235: .3) node [] {};
	\draw (T) -- +(305: .3) node [] {};
	\draw (T-1) -- +(235: .3) node [] {};
	\draw (T-1) -- +(305: .3) node [] {};		
	\path (T) -- (T-1) node[comp, fill, midway, above=of T] (B) {}; 
	\node[circle, draw, inner sep=0pt, minimum size=12pt] (D) [above=of B] {$2$}; 		
	\path (B) edge node [order top left,xshift=2pt,yshift=-4pt] {} (T);
	\path (D) edge node [order top left,xshift=2pt,yshift=-4pt] {} (B);
	\path (T-1) edge node [order top left,xshift=2pt,yshift=-4pt] {} (B);
	\end{scope}
	\begin{scope}[xshift=-1.2cm,yshift=-2.1cm] 
	\node[] (T) [] {}; 
	\node[] (T-1) [right=of T] {}; 
	\path (T) -- (T-1) node[comp, fill, midway, above=of T] (B) {}; 
	\draw (B) -- +(235: .3) node [] {};
	\draw (B) -- +(305: .3) node [] {};
	\draw (B) -- +(270: .37) node [] {};
	\draw (B) -- +(320: .55) node [] {};	
	\node[circle, draw, inner sep=0pt, minimum size=12pt] (D) [above=of B] {$2$}; 
	\path (B) edge node [order top left,xshift=2pt,yshift=-4pt] {} (T);
	\path (B) edge node [order top left,xshift=2pt,yshift=-4pt] {} (T-1);
	\path (D) edge node [order top left,xshift=2pt,yshift=-4pt] {} (B);
	\end{scope}	
	\begin{scope}[yshift=-4cm] 
	\node[comp,fill] (T) [] {};
	\node[] (T-1) [right=of T] {}; 		
	\draw (T) -- +(235: .3) node [] {};	
	\draw (T) -- +(305: .3) node [] {};		
	\path (T) -- (T-1) node[comp, fill, midway, above=of T] (B) {}; 
	\draw (B) -- +(270: .37) node [] {};
	\node[circle, draw, inner sep=0pt, minimum size=12pt] (D) [above=of B] {$2$};		
	\path (B) edge node [order top left,xshift=2pt,yshift=-4pt] {} (T);
	\path (D) edge node [order top left,xshift=2pt,yshift=-4pt] {} (B);
	\path (T-1) edge node [order top left,xshift=2pt,yshift=-4pt] {} (B);
	\end{scope}		
	\begin{scope}[xshift=4.5cm] 
	\node[comp,fill] (T) [] {};
	\node[] (T-1) [right=of T] {}; 	
	\draw (T) -- +(235: .3) node [] {};
	\draw (T) -- +(305: .3) node [] {};	
	\path (T) -- (T-1) node[midway] (B) {};
	\node[circle, draw, inner sep=0pt, minimum size=12pt] (D) [above=of B] {$2$};
	\draw (D) -- + (300: .6) node (B-1) [comp,fill] {};
	\draw (B-1) -- +(235: .3) node [] {};
	\draw (B-1) -- +(305: .3) node [] {};					
	\path (D) edge node [order top left,xshift=2pt,yshift=-4pt] {} (T);
	\end{scope}	
	\begin{scope}[xshift=4.5cm,yshift=-4cm] 
	\node[comp,fill] (T) [] {};
	\node[] (T-1) [left=of T] {}; 	
	\draw (T) -- +(235: .3) node [] {};
	\draw (T) -- +(305: .3) node [] {};	
	\path (T) -- (T-1) node[midway] (B) {};
	\node[circle, draw, inner sep=0pt, minimum size=12pt] (D) [above=of B] {$2$};
	\draw (D) -- + (240: .6) node (B-1) [comp,fill] {};
	\draw (B-1) -- +(235: .3) node [] {};
	\draw (B-1) -- +(305: .3) node [] {};					
	\path (D) edge node [order top left,xshift=2pt,yshift=-4pt] {} (T);
	\end{scope}	
	\begin{scope}[xshift=1.3cm,yshift=-.8cm] 
	\begin{scope}[xshift=4.5cm,yshift=-2cm] 
	\node[comp,fill] (T) [] {};
	\node[comp,fill] (T-1) [right=of T] {}; 		
	\draw (T) -- +(235: .3) node [] {};
	\draw (T) -- +(305: .3) node [] {};
	\draw (T-1) -- +(235: .3) node [] {};
	\draw (T-1) -- +(305: .3) node [] {};		
	\path (T) -- (T-1) node[circle, draw, inner sep=0pt, minimum size=12pt, midway, above=of T] (B) {$2$};	
	\path (B) edge node [order top left,xshift=2pt,yshift=-4pt] {} (T);
	\path (T-1) edge node [order top left,xshift=2pt,yshift=-4pt] {} (B);
	\path (T-1) -- +(35: .6) node [] {${~\overset{\widehat{\pi}}{\longleftarrow}~}$};
	\end{scope}	
	\begin{scope}[xshift=6.1cm,yshift=-2cm] 
	\node[comp,fill] (T) [] {};
	\node[comp,fill] (T-1) [right=of T] {}; 		
	\draw (T) -- +(235: .3) node [] {};
	\draw (T) -- +(305: .3) node [] {};
	\draw (T-1) -- +(235: .3) node [] {};
	\draw (T-1) -- +(305: .3) node [] {};		
	\path (T) -- (T-1) node[circle, draw, inner sep=0pt, minimum size=12pt, midway, above=of T] (B) {$3$};	
	\path (B) edge node [order top left,xshift=2pt,yshift=-4pt] {} (T);
	\path (T-1) edge node [order top left,xshift=2pt,yshift=-4pt] {} (B);
	\end{scope}		
	\end{scope}
	\begin{scope}[xshift=8.65cm,yshift=.5cm] 
	\begin{scope} 
	\node[comp,fill] (T) [] {};
	\node[comp,fill] (T-1) [right=of T] {}; 		
	\draw (T) -- +(235: .3) node [] {};
	\draw (T) -- +(305: .3) node [] {};
	\draw (T-1) -- +(235: .3) node [] {};
	\draw (T-1) -- +(305: .3) node [] {};		
	\path (T) -- (T-1) node[circle, draw, inner sep=0pt, minimum size=12pt, midway, above=of T] (B) {$2$};	
	\path (B) edge node [order top left,xshift=2pt,yshift=-4pt] {} (T);
	\path (T-1) edge node [order top left,xshift=2pt,yshift=-4pt] {} (B);
	\path (T-1) -- +(35: .6) node [] {${~\overset{\widehat{\pi}}{\longleftarrow}~}$};
	\end{scope}	
	\begin{scope} [xshift=1.6cm] 
	\node[comp,fill] (T) [] {};
	\node[comp,fill] (T-1) [right=of T] {}; 		
	\draw (T) -- +(235: .3) node [] {};
	\draw (T) -- +(305: .3) node [] {};
	\draw (T-1) -- +(235: .3) node [] {};
	\draw (T-1) -- +(305: .3) node [] {};		
	\path (T) node[circle, draw, inner sep=0pt, minimum size=12pt, above=of T] (B) {$2$};	
	\path (B) edge node [] {} (T);
	\path (T-1) node[circle, draw, inner sep=0pt, minimum size=12pt, above=of T-1] (B-1) {$2$};	
	\path (T-1) edge node [] {} (B-1);
	\path (T) edge node [] {} (B-1);
	\path (T-1) edge node [] {} (B-1);
	\path (B) edge node [] {} (T-1);
	\end{scope}		
	\end{scope}		
	\begin{scope} [xshift=9.3cm,yshift=-4.7cm] 
	\node[comp,fill] (T) [] {};	
	\draw (T) -- +(235: .3) node [] {};
	\draw (T) -- +(305: .3) node [] {};	
	\path (T) node[circle, draw, inner sep=0pt, minimum size=12pt, above=of T] (B) {$2$};
	\draw (B) -- +(0: .4) node [] {};
	\draw (B) -- +(315: .4) node [] {};		
	\path (B) edge node [] {} (T);
	\end{scope}				
	\draw [-latex] (3,-5.5) -- (3,-6.5); 
	\draw[name path=V] (-.8,-7.5) -- +(0:9.5)  node (V) [right] {$V$};
	\draw[name path=W1] (-.8,-7.38) -- +(-15:9.8)  node (W1) [right] {$W_1$};	
	\draw[name path=W2] (1,-7.1) -- +(-23.5:6.8)  node (W2) [above,midway] {$W_2$};	
	\fill[name intersections={of=V and W2}] (intersection-1) circle (2.8pt);
	\fill[name intersections={of=V and W1}] (intersection-1) circle (2.8pt);
	\fill[name intersections={of=W1 and W2}] (intersection-1) circle (2.8pt);
	\end{tikzpicture}
	\ees
        \caption{The modified Hitchin base for $g=2$. The map $b$ to $\AA_2^\circ$ contracts
        the divisors $F_i$ and $D_{ij}$}
        \label{fig:MHB}
\end{figure}
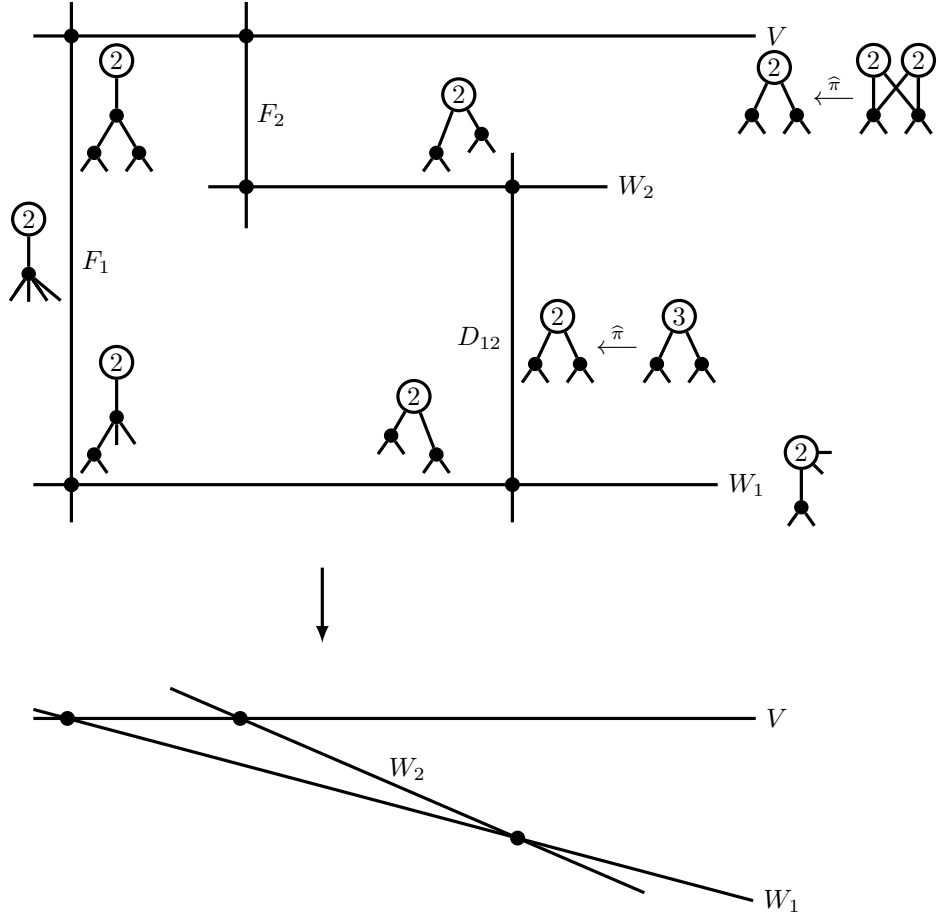

A similar constraint applies to a fourfold zero. There, the 4 marked points
on the lower level of a rational tail cannot move in an unconstrained
way (Check on the double cover: the top level is square of an abelian
differential, hence its double cover splits into two components. The double
cover of bottom level is an elliptic curve whose $j$-invariant is determined
by the position of the four simple zeros of the quadratic differential.
This is unconstrained, since the $j$-invariant is determined by the
absolute periods, which are in the $(-1)$-eigenspace. However the
point of attachment has to be chosen so that the residue of the
differential of type $(2,2,2,2,-4,-4)$ on the elliptic curve is
zero by the global residue condition.) Consequently, they can move in a
codimension one subvariety~$F_i$ of this~$\moduli[0,5]$, where
again $i=1,\ldots,6$ indexes the Weierstrass points. This subvariety
intersects the Weierstrass line~$W_i$ in the $3$-level graph, with
type $(-8,2,1,1)$ on middle level and type $(-6,1,1)$ on bottom.
It intersects the Veronese curve (the reducible locus) in the
$3$-level graph given by a cherry with a long stalk see Figure~\ref{fig:MHB}.
\par
To explain the singularity of~$F_i$ we compute that with labeled
points and the pole set as $x_5 = \infty$ the equation of the
locus of differentials
\bes
q \= (z-x_1)(z-x_2)(z-x_3)(z-x_4)dz^2
\ees
with vanishing $2$-residue at~$x_5$ is precisely the union of the
three affine-invariant loci
\bes
x_i + x_j \= x_k + x_\ell
\ees
for all $\{i,j,k,l\} = \{1,2,3,4\}$. The action of the symmetric
group permutes these three loci, and the stabilizer of, say, $x_1+x_2 = x_3+x_4$
is generated by the involutions~$(12)$ and $(34)$. This results in a
quotient singularity at the origin.
\end{example}
\par
\begin{remark} \textrm{
The dependence of $\BB_{X_\st}$ on~$X_\st$ for $g(X_\st) =3$ and similarly for
higher genus can easily be seen as follows. The modified Hitchin base
contains a 'boundary' divisor for each two-level graph. Among those graphs
is the 'compact type' graph~$\Gamma_8$ with two vertices and one level, and with
all marked points on bottom level. In this case the top level twisted
differential belongs to the stratum $\mu = (8)$. This stratum has
projectivized dimension~$4$. The forgetful map of this stratum
to $\moduli[3]$ is not dominant and for~$X_\st$ outside the range of
the forgetful map, the modified Hitchin base does not have a divisor
corresponding to $\Gamma_8$.
}\end{remark}

\section{Universal compactified Jacobians}
\label{sec:UnivJac}

In this section we recall the notions leading to a compactification
of the universal Jacobian over the moduli space of pointed stable
curves. Theorem~\ref{thm:compJ} is essentially contained in the literature,
see the references in the proof. Fix $g \geq 2$ and $n>0$, and let $d \in \ZZ$
denote the degree of the line bundles.
\par
First, we let the \emph{unrigidified universal compactified Jacobian}
$\tildecJ[g,n][d]$ be the algebraic stack parameterizing stable curves
$\calX \to B$ of genus~$g$ with~$n$ marked points together with a
family~$\calF$ of  torsion-free rank one-sheaves
on~$\calX$  of relative degree~$d$.  The restriction of $\tildecJ[g,n][d]$
to $\moduli[g,n]$ is the universal degree~$d$ Picard variety.
The multiplicative group $\GG_m$ acts by rescaling~$\calF$ in every
fiber. We denote by $\tildecJ[g,n][d] /\!\!/ \GG_m $ the
rigidification with respect to this group action, the \emph{(rigidified)
universal compactified Jacobian}. Finally we let $\barcJ[g,n][d,P] \subset 
\tildecJ[g,n][d] /\!\!/ \GG_m$ be the substack of $P$-semistable sheaves
with respect to a polarization $P$ (see below).  It comes with a forgetful
morphism $\barcJ[g,n][d,P] \to \barmoduli[g,n]$. The pointed canonical
polarization $P_{\can}$ defined in~\eqref{eq:canpol} is a universal polarization
on $\barmoduli[g,n]$ invariant under permutations of the numbering of the
marked points.
\par
\begin{theorem} \label{thm:compJ} 
The universal $P$-compactified Jacobian
$\barcJ[g,n][d,P]$ is a universally closed smooth (Artin) stack. If the
polarization~$P$ is non-degenerate it is a proper Deligne-Mumford stack. 
\par
The pointed canonical polarization $P_{\can}$ of degree~$d$ is non-degenerate
if and only if $\gcd(d-g+1,2g-2+n) = 1$. 
\end{theorem}
\par
Nodes of a curve~$X$ where~$\calF$ is not locally free will be
called \emph{Neveu-Schwarz nodes} for brevity.\footnote{This is inspired
by physics terminology, see e.g.\ \cite{JKV01}. The nodes where
the bundle is locally free are called Ramond nodes.}
The stabilizers of points in $\barcJ[g,n][d]$ may be infinite,
in fact precisely if the set of Neveu-Schwarz nodes separates the
stable curve. This is possible for degenerate polarization's only.
Consequently, $\barcJ[g,n][d]$ is not a Deligne-Mumford
stack and not separated (and thus not proper) near these points. 
\par
\medskip
\paragraph{\textbf{Torsion free rank one sheaves}} Let $X$ be a nodal
curve throughout in this section. Recall that a coherent sheaf~$\calF$ on~$X$ is
of {\em rank~$1$}, if it is of rank~$1$ at every generic point of~$X$.
It is {\em pure} if the support of every non-zero subsheaf is equal
to the support of~$\calF$. It is {\em torsion-free} if it is pure and
the support is equal to the whole curve~$X$. The {\em degree} of a
torsion-free sheaf~$\calF$ of rank one a curve~$X$ is defined by
$\deg(\calF) = \chi(\calF) - \chi(X)$.
\par
A torsion-free rank one sheaf~$\calF$
on a smooth curve is simply a line bundle. On a nodal curve~$X$ such a
sheaf is a line bundle away from the nodes. At a node~$P$ the stalk is either
isomorphic to~$\calO_{X,P}$ or to the maximal ideal~$\frakm_P$.
\par
If~$X = \cup_i X_i$ is reducible with say~$s$ components, we associate with
a torsion-free rank one sheaf~$\calF$ its multi-degree, i.e.
\bes
\deg(\calF) \= (\deg \calF_{X_1}, \ldots, \deg \calF_{X_s})\,,
\ees
where $\calF_{X_j}$ is the maximal torsion free quotient of $\calF|_{X_j}$.
If $f: \calX \to B$ is a family of curves, then a family of
torsion-free rank one sheaves~$\calF$ is a $B$-flat sheaf~$\calF$ on~$\calX$
whose fibers are rank one and torsion-free.
\par
\medskip
\paragraph{\textbf{Polarizations}} A \emph{numerical
$d$-polarization} on the stable curve~$X = \cup_{i=1}^s X_i$ is a
tuple $\phi = (\phi_1,\ldots, \phi_s)$ of rational numbers with total
degree $|\phi| = \sum \phi_i = d$. We can alternatively view it as a function
on the vertices of the dual graph~$\Gamma$ of~$X$. A \emph{universal
numerical $d$-polarization~$\phi$} for $\barmoduli[g,n]$
is a collection of numerical $d$-polarizations~$\phi_\Gamma$  for all dual graphs
of curves in $\barmoduli[g,n]$ subject to the compatibility condition that for
a contraction $c: \Gamma_1 \to \Gamma_2$ of dual graphs
\bes
\phi_{\Gamma_2}(v) \= \sum_{w \in c^{-1}(v)} \phi_{\Gamma_1}(w)\,. 
\ees
We often drop the index~$\Gamma$, if clear from the context.
In order to specify a universal polarization it is often convenient to use
a flat vector bundle $P$ of some rank~$r$ and degree $r(d-g+1)>0$
on the universal family~$\calX$. To such a vector bundle~$P$ we associate
the tuple
\bes
\phi \= \phi(P) \= \Bigl(\frac{\deg P|_{X_1}}{r}
+ \frac{\deg(\omega_{X_1})}{2}, \ldots, \frac{\deg P|_{X_s}}{r}
+ \frac{\deg(\omega_{X_s})}{2}\Bigr)\,.
\ees
In fact every numerical polarization can be given by a vector bundle,
see \cite[Remark~4.6]{KPStab}. Of special interest in the sequel is the
following polarization. Let
\bes
\omega_{\calX}(\bfz)\,:=\,\omega_{\calX}\bigl(\sum_{i=1}^{n} z_i\bigr) 
\ees
be the canonical bundle twisted by the universal sections~$z_i$. Then 
\be \label{eq:canpol}
{P}_{\can} \= {\omega}_{\calX}({\bfz})^{\otimes {d}- {g}+1}
\oplus \calO_{{\calX}}^{2{g}-3+n}
\ee
defines a universal $d$-polarization, the \emph{pointed canonical
  polarization $\phi_\can = \phi(P_\can)$}. In fact, there is a class
of polarizations
(discussed e.g.\ \cite[Paragraph 4.4.2]{melo16}) by varying the coefficients
of $z_i$ in the twisted  canonical bundle and we focus on the case of
coefficients~$+1$.
\par
\medskip
\paragraph{\textbf{Stability}}
The torsion-free rank one sheaf~$\calF$ is called \emph{$\phi$-semistable}
if for every non-empty proper subcurve $Y \subset X$ the 'basic inequality'
\be \label{eq:phisstab}
\deg(\calF_Y) \geq \sum_{X_i \subset Y} \phi_i - \frac{|Y \cap Y^c|}{2}
\ee
holds, where~$Y^c$ is the complement of~$Y$ in~$X$. The sheaf~$\calF$
is called \emph{$\phi$-stable}, if the above inequality is strict for all $Y$.
\par
The polarization~$\phi$ is called \emph{non-degenerate} if the
right hand side of~\eqref{eq:phisstab} does not assume integral
values for any dual graph~$\Gamma$ and for any proper subcurve~$Y$
of a curve~$X$ with dual graph~$\Gamma$. In particular, there
are no strictly semistable sheaves with respect to a non-degenerate
polarization.
\par
It will be convenient to give an equivalent formulation in terms of Euler
characteristics of the restriction of~$\calF$. Namely since
\bes \chi(\calF_{Y})\= \deg\calF_{Y} + 1-g_{Y} \=\deg\calF_{Y}
-\frac12 \bigl(\deg\omega_{{X}}|_{Y}-\vert {Y} \cap {Y}^c\vert \bigr)
\ees
we conclude that~\eqref{eq:phisstab} is equivalent to
\be \label{eq:phistabEC}
\chi(\calF_{Y}) \,\geq \, \frac{\deg{P}_{Y}}{\rk{P}}\,.
\ee
\par
\medskip
\paragraph{\bf Balanced line bundles on quasi-stable curves}
An alternative and essentially equivalent viewpoint uses balanced
line bundles on quasi-stable curves. Here a semistable pointed curve~$(X,\bfz)$
is called \emph{quasi-stable} if all ('exceptional') destabilizing
components are rational curves without marked points and if these
are disjoint and not contained in rational tails. A line bundle~$\calL$ is
called \emph{semibalanced} if (in the unpointed case) the \emph{basic inequality}
\be \label{eq:basicineq}
\Bigl| \deg(\calL|_Z) - \sum_{X_i \subset Z} \phi^\can_i \Bigr| \,\leq\,
\frac{|Y \cap Y^c|}{2}
\ee
holds. It is called \emph{balanced} if moreover $\deg(L|_E) = 1$ for
every exceptional component of~$X$. (See \cite{melo11} for the generalization
to the pointed case, where rational tails have to be taken into account
correctly so as to obtain a notion compatible with the morphism induced
by forgetting a marked point.)
The equivalence of the two notions is clarified by \cite[Proposition~5.4
and Proposition~6.2]{EstPac}, see also \cite[Remark~5.14]{KPStab}: The
pushforward of a balanced line bundle under the stabilization map is a
torsion-free rank-$1$ sheaf and the balancing condition~\eqref{eq:basicineq}
translates into~\eqref{eq:phisstab}.
(On the level of coarse moduli spaces this is the main content of
Pandharipande's compactification \cite{PandCompJ}.)
\par
The set of all universal numerical $d$-polarizations fits into a space~$V_{g,n}^d$
analyzed in detail in~\cite{KPStab}. Each $\phi \in V_{g,n}^d$ defines
a substack $\barcJ[g,n][d,\phi] \subset  \tildecJ[g,n][d] /\!\!/ \GG_m$
and loc.~cit. studies the dependence of the geometry of the substack
on~$\phi$.

\par
\medskip
\paragraph{\bf The versal deformations of sheaves on nodal curves}
To justify smoothness of $\barcJ[g,n][d]$ and to prepare for the next
section we recall the local deformation theory of sheaves on nodes
from \cite{CMKV}. Let $\wR = \CC[[x,y,t]]/(xy-t)$ be the complete local
algebra of a node and $\calX = \Spec \wR$ the family of curves degenerating
to a node. The miniversal deformation ring of a locally free
sheaf on a node is simply the base ring $\CC[[t]]$ that parameterizes
the degeneration to the node. The miniversal deformation ring of the
torsion-free rank one sheaf $\frakm \subset \wR_0 = \CC[[x,y]]/(xy)$ is
\be \label{eq:deforingS}
S = \CC[[a,b,t]]/(ab-t),
\ee
see \cite[Lemma~3.13]{CMKV}. 
In fact, the miniversal deformation of~$\frakm$ on the pullback family
$\calX_S = \calX \times_{\Spec \CC[[t]]} \Spec S$ is
\be \label{eq:calI}
\calI = \langle x-a, y-b \rangle \,.
\ee
It is also convenient to present $\calI$ as $\CC[[x,y,a,b]]/(xy-ab)$-module as 
\be \label{eq:praescalI}
\calI \= \langle s_1, s_2 \mid xs_1=-as_2, \, ys_2=-bs_1 \rangle\,.
\ee
\begin{proof}[Proof of Theorem~\ref{thm:compJ}]
The construction of the stack as GIT-quotient is given in \cite{melo09},
building on the earlier work of Caporaso \cite{CapoCompJ}. The
isomorphism to the functor described here is given in
\cite[Theorem~6.3]{EstPac}. The source cite above works in fact without
marked points, but \cite[Remark~5.14]{KPStab} explains how to fix this.
\par
The miniversal deformations ring of any torsion-free rank-one
sheaf at any node is smooth, combining~\eqref{eq:deforingS} and the obvious
locally free case.
The full deformation space is a product of the deformations at
the node and a power series ring, since the forgetful map from
deformations of the pair $(X,F)$ to the product of local deformations
at the node is formally smooth, see \cite[Section~6]{CMKV}.
\par
The stack $\tildecJ[g,n][d] /\!\!/ \GG_m$ satisfies the valuative
criterion for properness by \cite[Theorem~32]{Esteves}. Since $\barcJ[g,n][d]$
is a closed substack and moreover quasi-compact, it is
universally closed.
\par
The non-degeneracy of the canonical polarization under the gcd hypothesis
follows from \cite{CapoNeron}. In this case the complement of the set of
Neveu-Schwarz nodes is connected (for otherwise applying the definition
of semistability to the two components violates non-degeneracy) 
and the automorphism group of the rigidified functor is trivial
(\cite[Theorem~A]{CMKV}). Hence the stack is Deligne-Mumford in this case.
\par
For the second statement we only need to check
the pointed canonical variant of the criterion for the polarization
to be non-degenerate. Suppose the polarization has this property. Then
for any subcurve~$Y$, say with $|Y \cap Y^c| =k$ we must have
\ba \label{eq:degintegral}
\ZZ \not\ni & \phantom{\=} \, \frac{\deg(\omega_Y)}2  + \frac{\deg P|_Y}{2g-2+n}
-\frac{k}2 \\
& \= (g_Y - 1) + \frac{(d-g+1)(2g_Y-2+k+n_Y)}{2g-2+n} \\
\ea
so that in particular the second summand is not integral. We may probe
this for any 'banana curve', i.e.\ with two irreducible components connected
by $k=2$ nodes. We may take any $g_Y \leq g-2$ and~$n_Y=n$ or $g_Y = g-1$
with $n_Y = n-1$ to get a pointed stable curve. This implies that $2g_Y-2+k+n_Y$
attains all integers less or equal to $2g-2+n-1$. Hence if $\gcd(d-g+1,2g-2+n)
= \ell>1$
we may arrange for $2g_Y-2+k+n_Y = (2g-2+n)/\ell$ and get a contradiction.
\par
Conversely, if $\gcd(d-g+1,2g-2+n) =1$, \eqref{eq:degintegral} follows
since \bes 2g_Y-2+k+n_Y < 2g-2+n\ees for any subcurve of a pointed stable curve.
\end{proof}
\par

\section{The spectral correspondence} \label{sec:spectralcorr}

In this section we define an extension of the notion of spectral data on smooth
Hitchin fibers to spectral data on semistable curves. Our notion of spectral
data builds on the universal compactified Jacobian over the modified
Hitchin base constructed in the previous sections. When specifying
Higgs-related data we write '$\GL(2,\CC)^\circ$' as shorthand for 'trace-free
$\GL(2,\CC)$'.
\par 
We start by defining multi-scale spectral data point-wise and in families.
Then we analyze the pushforward of multi-scale spectral data. This motivates
the definition of multi-scale Higgs pairs and allows to formulate spectral
correspondences: A pointwise statement in Theorem~\ref{thm:Higgspairs} and
a version for families in Theorem~\ref{thm:Higgspairsfam}. These theorems
will be proven in the following section.  
\par

\subsection{Multi-scale spectral data}
Let $(X,\bfz)$ be a pointed stable curve. 
\par
\begin{definition} \label{def:SDGL2tf}
A \emph{$\GL(2,\CC)^\circ$-multi-scale spectral datum}
$(\wh{\Gamma}, \bfq, \bftau, \calF)$ of degree $\wh{d}$ on~$(X,\bfz)$
is a quadratic multi-scale differential $(\wh{\Gamma},\bfq,\bftau)$
together with a torsion-free rank-one sheaf~$\calF$ of degree $\wh{d}$
on $\wh{X}$.
\par Let $\wh{P}$ be a polarization on $\wh{X}$.
A $\GL(2,\CC)^\circ$-multi-scale spectral datum is called
\emph{$\wh{P}$-semistable} if $\calF$ is $\wh{P}$-semistable. 
\end{definition} 
\par
Semistable multi-scale spectral data are points of the fiber product
\be
\SD_g \=  \overline{\calQ}_{g,[4g-4]} \times_{\barmoduli[\wh{g},{[4g-4]}]} \barcJ[\wh{g},4g-4][\wh{d},\wh{P}]\,.
\ee

We comment on changing the polarization and on changing the Lie group below.
As a consequence of Theorem~\ref{thm:kLMS} and Theorem~\ref{thm:compJ} we record:
\par
\begin{proposition} \label{prop:barHiggsStack} Let $\wh{P}$ be a universal
polarization on $\barmoduli[\wh{g},\{4g-4\}]$. The space of multi-scale spectral
data $\SD_g$ is a smooth Artin stack that admits an action of $\CC^*$. The quotient
$\PP \SD_{g} =  \SD_{g}/\CC^*$ is universally closed. It is a proper Deligne-Mumford
stack if $\wh{P}$ is non-degenerate. 
\end{proposition}
\par
From the definition of $\SD_g$ it is clear how to define a family of
multi-scale spectral data. We record it in order to fix notation.
\par
\begin{definition} A \textit{germ of a family of multi-scale
$\GL(2,\CC)^\circ$-spectral data} over a family $\calX \rightarrow S$ is
a tuple $(\wh{\Gamma},\bfq,\bftau,\calF)$ of a family of multi-scale quadratic
differentials $(\wh{\Gamma}, \bfq, \bftau)$ on $\calX$ and flat family of
torsion-free sheaves $\calF$ on $\wh{\calX}$.
\end{definition}
One obtains the notion of a family of multi-scale $\GL(2,\CC)^\circ$-spectral data by patching together these germs. We refer to \cite[Section 7]{LMS} for further details on the 'sheafification process'.

\subsection{The local form of the pushforward} \label{sec:localHiggs}
We analyze the local form of the push-forward of a family of multi-scale spectral data $(\wh{\Gamma},\bfq,\bftau,\calF)$ at the nodes. We restrict the double covering $\wh{\pi}: \wh\calX \to \calX$ to a
neighborhood~$\calX$ of a deformation of a node $e$ and the ($\sigma$-invariant)
preimage. We may pretend that this is a node joining level~$0$ and $-1$
and label all objects accordingly. There are four cases to consider depending on whether $\kappa_e$ is odd or even and whether $\calF$ is locally free in the special fiber or not. When $\kappa_e$ is odd the node has a single preimage $\wh{e}$ fixed by the involution $\sigma$. When $\kappa_e$ is even the node has two preimages $\wh{e}_1,\wh{e}_2$ interchanged by $\sigma$.
\par 
We will assume the multi-scale abelian differential $\bflambda$ on
$\wh{\calX}$ is given in the normal form of \cite[Theorem 4.3]{LMS} 
\bes \lambda_0=(x^{\kappa_{\wh{e}}} + r)\frac{dx}{x}, \quad \lambda_{-1}=-(y^{-\kappa_{\wh{e}}}+\frac{r}{t^{\alpha\kappa_{\wh{e}}}}) \frac{dy}{y}
\ees over a deformation of the node of the form $\{xy=t^\alpha\}$. Here the
residuum $r \in \CC[[t]]$ can be assumed to be divisible by $t^{\alpha\kappa_e}$.
The abelian differential on the canonical double cover $\wh{X}$ is
$\sigma$-anti-symmetric.
Hence the residuum is zero at nodes $\wh{e}$ fixed by $\sigma$ (see also
\cite[Theorem 3.1]{kdiff}). 
\par
\medskip
\paragraph{\textbf{$\kappa_e$ odd and $\calF$ is locally free}}
The covering is given by the inclusion of rings
\ba \label{eq:covfix}
\wh{\pi}^\sharp\,:\,\, &R:=\mathbb{C}[[u,v,t]]/(uv-t^{2\alpha})
\,\rightarrow\,\hat{R}= \mathbb{C}[[x,y,t]]/(xy-t^{\alpha}), \\
 &	u \mapsto x^2,\qquad v \mapsto y^2, \qquad t \mapsto t
\ea
for some $\alpha$ determined by the family of curves, and
$\calF \cong \wh{R}$, so that
\be \label{eq:LF/fixframe}
\tag{LF/fix} 
\calE := \wh{\pi}_* \wh{R} \= R \oplus \langle x,y \rangle_R\,.
\ee
\par
The divisor $X_1$ is given by $u=0$ in these coordinates, so~$u$
is invertible on ${\calX\setminus{X_1}}$ and so $\calE|_{\calX\setminus{X_1}}=
1+\langle x\rangle_R$ (since $y = x\cdot t^a/u$).
Similarly $\calE|_{\calX\setminus{X_0}}=
1+\langle y\rangle_R$. Multiplication with $\lambda_i$ on $\calF$ induces Higgs fields 
\[ \Phi_0:  \calE|_{\calX\setminus{X_1}} \rightarrow  \calE \otimes \omega_{\calX}
|_{\calX\setminus{X_1}}, \qquad \Phi_1:  \calE|_{\calX\setminus{X_0}} \rightarrow
\calE \otimes \omega_{\calX}|_{\calX\setminus{X_0}}
\] given in this basis by
\ba \label{eq:Phi01fixed}
\Phi_0\=\wh{\pi}_*{\lambda_0}|_{\calX\setminus{X_1}}: &
	\qquad 1 \mapsto u^{\frac{\kappa_{\wh{e}}-1}{2}}x \otimes \frac{d u}{u},
	\quad &&x \mapsto u^{\frac{\kappa_{\wh{e}}+1}{2}} \otimes \frac{d u}{u},  
	\\
        \Phi_1\=\wh{\pi}_*{\lambda_1}|_{\calX\setminus{X_0}}: &
	\qquad 1 \mapsto v^{-\frac{\kappa_{\wh{e}}+1}{2}}y \otimes \frac{d v}{v}, 
	\quad &&y \mapsto v^{-\frac{\kappa_{\wh{e}}-1}{2}} \otimes \frac{d v}{v}.
\ea
\par
\medskip
\paragraph{\textbf{$\kappa_e$ odd and $\calF$ is Neveu-Schwarz}}
Using the same covering~\eqref{eq:covfix} the sheaf $\calF$ is now
locally the pullback of the miniversal deformation~$\calI$ from~\eqref{eq:calI}
with~$t$ replaced by~$t^\alpha$. Consequently, 
\be
\tag{NS/fix} \label{eq:NS/fixframe}
\calE \,:=\, \wh{\pi}_* \calI \=
\langle x-a,y-b,x(x-a),y(y-b)\rangle_R\,.
\ee
For $u \neq 0$ (resp. $v \neq 0$) this generating set simplifies to 
\[ \calE|_{\calX\setminus{X_1}} \= R(x-a)  \oplus R x(x-a)  \quad \text{resp. }
\calE|_{\calX\setminus{X_0}} \= R (y-b) \oplus R y(y-b).
\] The two Higgs fields are given by 
\begin{align*}
	\Phi_0: \quad &(x-a) \mapsto  u^{\frac{\kappa_{\wh{e}}-1}{2}}x(x-a) \otimes \frac{d u}{u}, \quad x(x-a) \mapsto u^{\frac{\kappa_{\wh{e}}+1}{2}}(x-a) \otimes \frac{d u}{u} \\
	\Phi_1: \quad &(y-b) \mapsto  v^{-\frac{\kappa_{\wh{e}}+1}{2}}y(y-b) \otimes \frac{d v}{v}, \quad y(y-b) \mapsto v^{-\frac{\kappa_{\wh{e}}-1}{2}} (y-b) \otimes \frac{d v}{v}\,.
\end{align*}
\par
\medskip
\paragraph{\textbf{$\kappa_e$ even and $\calF$ is locally free}}
The covering is given by the inclusion
\ba \label{eq:covswap}
\wh{\pi}^\sharp\ = (\wh{\pi}_1^\sharp,\wh{\pi}_2^\sharp)
\,:\,\, &R:= \frac{\mathbb{C}[[u,v,t]]}{uv-t^{\alpha}}
\,\rightarrow\,\widehat{R}= \frac{\mathbb{C}[[x_1,y_1,t]]}{x_1y_1-t^{\alpha}} \oplus
\frac{\mathbb{C}[[x_2,y_2,t]]}{x_2y_2-t^{\alpha}} \=:\, \wh{R}_1 \oplus \wh{R}_2, \\
 &	u \mapsto (x_1,x_2), \qquad v \mapsto (y_1,y_2), \qquad t \mapsto t
\ea
of rings for some $\alpha$. Note the index $j=1,2$ of $\wh{R}_j$ and $i=0,-1$ of the levels of the special fiber $X_i$. Let $\calF$ be free module of rank 1 over $\wh{R}$ and $\calE=\wh{\pi}_* \calF$. Then 
\be
\tag{LF/swap}
\calE \= \wh{\pi}_*\calF = \wh{\pi}_{1*}\calF|_{\wh{U}_1} \oplus
\wh{\pi}_{2*}\calF|_{\wh{U}_2}\,,
\ee
where $\wh{U}_j=\Spec (\wh{R}_j)$ with $j=1,2$. The multi-scale abelian
differential~$\bflambda$ is $\sigma$-antisymmetric and hence given by
level-wise abelian differentials 
\bas
\lambda_{0} &\= \left(\Bigl(x_1^\frac{\kappa_e}{2} + r\Bigr) \frac{d x_1}{x_1},\,
-\Bigl(x_2^\frac{\kappa_e}{2}+r\Bigr) \frac{d x_2}{x_2}\right), \\ 
\lambda_{-1} &\=\left(-\Bigl(y_1^{-\frac{\kappa_e}{2}} + \frac{r}{t^{\alpha\kappa_{\wh{e}}}}\Bigr)
\frac{d y_1}{y_1},\,\Bigl(y_2^{-\frac{\kappa_e}{2}}+\frac{r}{t^{\alpha\kappa_{\wh{e}}}}\Bigr)
\frac{d y_2}{y_2}\right)\,.
\eas
There exist unique abelian differentials $\eta_{i}$ on $\calX\setminus {X_{-1-i}}$,
such that 
\bes 
\wh{\pi}^*\eta_0 \= (x_1^\frac{\kappa_e}{2}+r) \frac{d x_1}{x_1} \quad \text{resp.}
\quad \wh{\pi}^*\eta_{-1} \=-(y_1^{-\frac{\kappa_e}{2}}+\frac{r}{t^{\alpha\kappa_{\wh{e}}}})
\frac{d y_1}{y_1}\,.
\ees
The Higgs fields are diagonal with respect to the given splitting of $\calE$.
They are given by
\be \label{eq:Higgsswap}
\Phi_i \= \wh{\pi}_*{\lambda_i}|_{\calX\setminus{X_{-1-i}}}
\ = \begin{pmatrix} \eta_i & 0 \\ 0 & - \eta_i \end{pmatrix}
\ee
for $i=0,-1$ with respect to such a frame. 
\par
\medskip
\paragraph{\textbf{$\kappa_e$ even and $\calF$ is Neveu-Schwarz}}
There are two cases to consider depending on whether $\calF$ is Neveu-Schwarz at one or at both nodes interchanged by $\sigma$. However, they work quite similarly. First consider the case, where $\calF$ is Neveu-Schwarz at the node $\wh{e}_1$ and locally free at $\wh{e}_2$. To take into account the equisingular deformation of $\calF_{\wh{e}_1}\cong\frakm_{\wh{e}_1}$ we have to work over the miniversal deformation ring $S$ recorded in Section \ref{sec:UnivJac} \eqref{eq:deforingS}. To compute the fiber product with the covering \eqref{eq:covswap} let $R=\CC[u,v,a,b,t]/(uv-t^\alpha,ab-t^\alpha)$ and 
\[ \wh{R}_j=\CC[x_j,y_j,a,b,t]/(x_jy_j-t^\alpha,ab-t^\alpha) \]
for $j=1,2$. Then the covering over $S$ is given by
\begin{align*} 
	\wh{\pi}^\sharp:	R \rightarrow \wh{R}_1 \oplus \wh{R}_2, \quad 
	u \mapsto (x_1,x_2), v \mapsto (y_1,y_2).
\end{align*} Let $\calI_1=\langle x_1-a,y_1-b \rangle \subset \wh{R}_1$ be the miniversal deformation of $\frakm_{\wh{e}_1}$. Then the family of torsion-free sheaves $\calF$ is given by
\be \label{eq:calIswap}
\calF \= \calI_1 \oplus \wh{R}_2 \=
\langle x_1-a,y_1-b\rangle \oplus \wh{R}_2\,.
\ee
Consequently, 
\be
\tag{NS/swap}
\calE := \wh{\pi}_* \calF \= \langle u-a, v-b \rangle \oplus
R \subset R^2  \,.
\ee
is a deformation of $\frakm_{e} \oplus R_e$. The local form of the Higgs field is given by the same matrix
as in~\eqref{eq:Higgsswap}.
\par 
In the case, of $\calF$ being Neveu-Schwarz at both nodes $\wh{e}_1,\wh{e}_2$ we have to consider the fiber product of the covering \eqref{eq:covswap} with the direct sum of two miniversal deformation rings $S_1 \oplus S_2$. The remainder of the construction can be easily adapted. In this case, the family of torsion-free sheaves $\calE=\wh{\pi}_*\calF$ is a deformation of $\frakm_e \oplus \frakm_e$ at the node. 
\par
\medskip
\paragraph{\textbf{Horizontal Nodes}}
At horizontal nodes the canonical covering $\wh{\pi}: \wh{X} \rightarrow X$ is unbranched and the local models look the same as in the case of $\kappa_e \equiv 0 \mod 2$. Instead of considering two irreducible components on two different levels $i=0,-1$ we have two irreducible components on one level $i$, that we index by $a,b$, such that on the $a$-component we have $x$-coordinates and on the $b$-component we have $y$-coordinates. The only difference to vertical nodes is that the component-wise abelian differentials glue to a section of $\omega_{\wh{\calX}}$ resp. $\omega_{\wh{X}}$. Choosing local coordinates as in \eqref{eq:covswap} they are given by
\[ \lambda_a \= \Bigl(r \frac{d x_1}{x_1},-r \frac{d x_2}{x_2}\Bigr), \qquad
\lambda_b \=\Bigl(-r\frac{d y_1}{y_1},r \frac{d y_2}{y_2}\Bigr), 
\]
with residue $r \in \CC$. This yields a gluing of the component-wise Higgs fields 
\[ \Phi_a \= \begin{pmatrix} r \frac{d u}{u} & 0 \\ 0 & -r \frac{d u}{u}
\end{pmatrix}, \quad
\Phi_b \= \begin{pmatrix} -r \frac{d v}{v} & 0 \\ 0 & r \frac{d v}{v}
\end{pmatrix}
\] to a Higgs field $\Phi: \calE \rightarrow \calE \otimes \omega$.

\subsection{Multi-scale $\GL(2,\CC)^\circ$-Higgs pairs}
The previous computations motivate the following definition.

\begin{definition} \label{def:Higgspair}
A \emph{multi-scale $\GL(2,\CC)^\circ$-Higgs pair} of degree $d$ on a $4g-4$-pointed
stable curve $(X,\bfz)$ is a tuple $(\wh{\Gamma}, \calE,\boldsymbol{\Phi}, \bftau)$
consisting of a level graph structure on a double covering $\wh{\pi}: \wh{\Gamma}
\to \Gamma$ of the dual graph of~$X$, a torsion-free rank 2 sheaf~$\calE$ on~$X$, Higgs
fields $\boldsymbol{\Phi}=(\Phi_0,\dots,\Phi_L)$ if~$\Gamma$
has~$L$ levels below zero, and prong-matchings~$\bftau$ for the
collection~$\bfq = (q_i)$ of differentials with simple zeros at the marked points,
subject to the following conditions.
\begin{itemize}
\item[i)] If we denote by $\calE_i$ the restriction $\calE |_{X_i}$ mod torsion,
then  $\Phi$ induces
\be \label{eq:HiggsEipointwise}
\Phi_i: \calE_i \rightarrow \calE_i \otimes M_i\,,
\ee
where~$M_i$ are the twists of the canonical bundle defined in~\eqref{eq:Micoverbundle},
such that $q_i = \det(\Phi_i)$\,.
\item[ii)] The traces of the Higgs fields vanish, i.e.\  $\Tr(\Phi_i)=0$.
\item[iii)] The collection~$\bfq = (q_i)$ together with~$\bftau$ is a
multi-scale differential of type $\mu = (1^{4g-4})$ compatible with the
level graph~$\Gamma$.
\item[iv)] For each node $e$ there exists an analytic neighborhood $U$, such that $(\calE,\boldsymbol{\Phi})|_U$ is given by the specialization to $t=0$ of one of the local forms of Section~\ref{sec:localHiggs}.
\end{itemize}
\par
Two Higgs pairs are called \emph{equivalent}, if they are
in the same orbit of the level rotation torus, acting by rescaling the
Higgs fields $\Phi_i$ for $i<0$ and simultaneously on~$\bftau$.
\end{definition}
\par
Note that the condition on~$q_i$ to form a multi-scale differential implies
in particular that all the Higgs fields~$\Phi_i$ are non-zero.
\par
\begin{remark} At a node $e$ with $\kappa_e$ even, including the case of
  horizontal nodes, the condition~(iv) is
equivalent to the existence of a splitting of~$\calE$ into rank 1 subsheaves
in an analytic neighborhoods of $e$, such that the Higgs fields~$\Phi_{\ell(e^+)}$
and~$\Phi_{\ell(e^-)}$ are diagonal with respect to the splitting. Here, $\ell(e^\pm)$
denote the respective levels of the preimages of the node~$e$.  However,
for $\kappa_e \equiv 1 \mod 2$ we are missing an abstract reformulation of the
condition~(iv). 
\end{remark}
\par
\begin{definition} \label{def:polarsing_bundle}
We define stability on multi-scale Higgs pairs as follows.
\begin{itemize}
\item[i)] A \emph{polarization~$P$} of a multi-scale $\GL(2,\CC)$-Higgs pair of
degree~$d$ is a locally free sheaf $P$ on $X$, such that $\deg(P)=\rk(P)(d-2g+2)$. 
\item[ii)] Let $P$ be a polarization of multi-scale $\GL(2,\CC)$-Higgs pairs
of degree~$d$ on~$X$. Then a multi-scale $\GL(2,\CC)$-Higgs pair is
called \emph{semistable}, if for every subsheaf
$\calG \hookrightarrow \calE$ that is $\bfPhi$-stable (i.e.\ such that there is a
factorization $\Phi_i: \calG|_{X_i} \to \calG|_{X_i} \otimes M_i$),
we have 
\be \chi(\calG) \,\leq \,
\frac{\sum_{v=1}^s\rk(\calG|_{X_v})\deg(P|_{X_v})}{2\rk P}\,. \label{eq:GL(2)stab}
\ee
The sum runs over all irreducible components $X_v \subset X$ for $v=1,\dots,s$. 
\end{itemize}
\end{definition}
The same convention, an index $v=1,\dots,s$ refering to vertices $v \in \Gamma$
and thus to irreducible components of~$X$, is used frequently in the sequel.
\par
\begin{example} Let $\calL$ a polarizing line bundle on $X$. Then 
\[ P\= \calL^{\otimes d-2g+2} \oplus \calO_X^{\deg \calL -1}
\] defines a polarization of multi-scale $\GL(2,\CC)$-Higgs pairs of degree~$d$
on~$X$. We recover the Simpson's $p$-stability condition for Higgs pair with
respect to $\calL$ 
\bes \frac{\chi(\calG)}{\sum_{v=1}^s\rk(\calG|_{X_v})\deg(\calL|_{X_v})} \,\leq\,
\frac{\chi(\calE)}{2\deg (\calL)}
\ees (cf. \cite[Definition 2.2]{BBN}). A universal choice for~$\calL$
over $\overline{\calQ}_{g,4g-4}(1^{4g-4})$ is the pullback of the pointed canonical bundle $\omega_X(\bfz)$ along $\overline{\calQ}_{g,4g-4}(1^{4g-4}) \to \barmoduli[g,4g-4]$. It is invariant under permuting the marked points and hence descends to the quotient stack by the action of the symmetric group.
\end{example} 
\par
\begin{lemma}\label{lem:pullback_of_polar} If $P$ is polarization of multi-scale $\GL(2,\CC)$-Higgs pairs of degree $d$, then $\wh{P}=\wh{\pi}^*P \oplus \calO_{\wh{X}}^{\rk(P)}$ is polarization of degree $\wh{d}=d+2g-2$ in the sense of Section~\ref{sec:UnivJac}.
\end{lemma}
\par
 We refer to $\wh{P}$ as the \emph{induced polarization} on $\wh{X}$.
\par
\begin{proof} We have $\rk(\wh{P})=2\rk(P)$ and 
\be \deg(\wh{P})=2\rk(P)(d-2g-2) \=\rk(\wh{P})(\wh{d}-4g+4). 
\ee Hence $\wh{P}$ satisfies the condition for a polarization of torsion-free
rank~$1$ sheaves of degree $\wh{d}=d+2g-2$.
\end{proof}
\par
\par
\begin{remark} 
It is not clear how to define a pushforward of polarizations in terms of the
bundles~$\wh{P}$ and~$P$. However, the stability condition associated to a
polarization~$\wh{P}$ of degree $\wh{d}$ on $\wh{X}$ is uniquely defined by
the \textit{slopes} (cf.\ \eqref{eq:phistabEC})
\bes \hat{s}_v \,:=\,\frac{\deg \wh{P}|_{\wh{X}_v}}{\rk \wh{P}} \,. 
\ees
We can define a pushforward of a stability condition in terms of these slopes by
\bes s_v^\flat \=\frac{1}{\sharp\wh{\pi}^{-1}(v) } \sum_{\hat{v} \in \wh{\pi}^{-1}(v)}
\hat{s}_{\hat{v}}\,.
\ees for all irreducible components $X_v$. This is compatible with the notion
of pullback of a polarization in the previous lemma.
\end{remark} 
\par 
Now, the pointwise
version of the correspondence is the following generalization of
Theorem~\ref{intro:BNR}.
\par
\begin{theorem} \label{thm:Higgspairs} Let~$X$ be a pointed stable curve. Let $P$ be a polarization of multi-scale $\GL(2,\CC)$-Higgs pairs of degree $d$ on~$X$ and $\wh{P}$ the induced polarization on $\wh{X}$.
Then there is a bijection between $\wh{P}$-(semi)stable $\GL(2,\CC)^\circ$-multi-scale spectral
data on~$\wh{X}$ and equivalence classes of $P$-(semi)stable multi-scale $\GL(2,\CC)^\circ$-Higgs pairs
on~$X$. 
\end{theorem}
\par \medskip
Our BNR-correspondence is in fact functorial, i.e.\ works in families.
To state this we globalize Definition~\ref{def:Higgspair}.
Since multi-scale differentials are defined by gluing (equivalence classes of)
germs, we will state the correspondence at the level of germs. For a germ of a
family of stable curves $f: {\calX} \to S$ we define
$\calX_{i^\complement} = \calX \setminus \cup_{j \neq i} X_j$ the complement
of the curves in the special fiber that are not at level~$i$. 
\par
There is one major difference in the form of the Higgs field between 
Definition~\ref{def:Higgspair} and Definition~\ref{def:germHiggspair} already
present in the construction of multi-scale differentials: In the pointwise
situation, and more generally for equisingular deformations,  the nodes are given
by sections of the family. They can be used for twisting line bundles, and
hence for the definition of the $M_i$ in~\eqref{eq:HiggsEipointwise}. For
families, say with smooth generic fiber, the nodes are of higher codimension.
There, alternatively, we define the Higgs field in~\eqref{eq:HiggsEifamilies}
away from the nodes. Extension across codimension two allows to recover the
local structure near the node and so the two definitions are compatible
under pullback.
\par
We emphasize that the subsequent definition treats only the case of families
with smooth generic fiber. We leave it to the reader to state the obvious
generalization to the general case where some nodes remain equisingular
while others are smoothened. 
\par
\begin{definition} \label{def:germHiggspair} 
A \emph{germ of a family of multi-scale $\GL(2,\CC)^\circ$-Higgs pairs} on a
germ $f: {\calX} \to S$ of a family of $4g-4$-pointed stable curves with
smooth generic fiber and special fiber~$(X,\bfz)$
is a tuple $(\wh{\Gamma}, \calE, \boldsymbol{\Phi}, \bftau)$ consisting 
of a level graph structure on a double covering $\wh{\pi}: \wh{\Gamma}
\to \Gamma$ of the dual graph of~$X$, a flat family of torsion-free sheaves~$\calE$ of rank $2$ on~$\calX$,
Higgs fields $\boldsymbol{\Phi}=(\Phi_0,\dots,\Phi_L)$ if~$\Gamma$
has~$L$ levels below zero, a collection~$\bfq = (q_i)$ of differentials
$q_i = \det(\Phi_i)$ with simple zeros at the marked points and prong-matchings~$\bftau$
for~$\bfq$  at the persistent nodes, defined with the following conditions.
\begin{itemize}
\item[i)] There are $\calO_{\calX_{i^\complement}}$-module homomorphisms
\be \label{eq:HiggsEifamilies}
\Phi_i: \calE |_{\calX_{i^\complement}} \rightarrow \calE |_{\calX_{i^\complement}}
\otimes \omega_{\calX/S}|_{\calX_{i^\complement}} \,,
\ee
such that $q_i = \det(\Phi_i)$ 
\item[ii)] The traces of the Higgs fields vanish, i.e.\  $\Tr(\Phi_i)=0$.
\item[iii)] The collection~$\bfq = (q_i)$ together with~$\bftau$ is a
multi-scale differential of type $\mu = (1^{4g-4})$ compatible with the
level graph~$\Gamma$.
\item[iv)] At every node $e$ the pair $(\calE,\bfPhi)$ is given by a base change of the local forms described in \ref{sec:localHiggs}.
\end{itemize}
\par
Two germs of families of Higgs pairs are called \emph{equivalent},
if they are in the same orbit of a section over~$S$ of the level rotation torus,
acting by rescaling the Higgs fields $\Phi_i$ for $i<0$ and simultaneously
on~$\bftau$.
\par
A germ of a family of Higgs pairs is called \emph{semistable}
if fiberwise the stability condition~\eqref{eq:GL(2)stab} holds.
\end{definition}
\par
We have prepared the statement of the correspondence in families:
\par
\begin{theorem} \label{thm:Higgspairsfam} Let~$f: \calX \to S$ be a germ of a
family of pointed stable curves with smooth generic fiber and let $P$ be a polarization
on $\calX$. Then there is a bijection between germs of $\wh{P}$-(semi)stable
$\GL(2,\CC)^\circ$-multi-scale spectral data on~$\calX$ and equivalence classes
of germs of $P$-(semi)stable multi-scale $\GL(2,\CC)^\circ$-Higgs
pairs on~$\calX$.
\end{theorem}

\subsection{The push-forward correspondence}

One direction of Theorem~\ref{thm:Higgspairs} without addressing stability yet
is the following proposition:
\par
\begin{proposition} Let $(\wh{\pi}: \wh{X} \to X, \bfq, \bftau, \calF)$
be a \emph{$\GL(2,\CC)^\circ$-multi-scale spectral datum} on~$X$ or on a germ of
family $f: \calX \to S$. Then  $\calE = \wh{\pi}_* \calF$ 
together with the Higgs fields $\Phi_i = \wh{\pi}_*(\,\,\cdot\, \lambda_i)$
defines a \emph{multi-scale $\GL(2,\CC)^\circ$-Higgs pair} $(\wh{\Gamma},
\calE,\boldsymbol{\Phi}, \bftau)$ on $X$ resp.\ on~$\calX$.
\end{proposition}
\par
\begin{proof} Since $\wh{\pi}|_{X_i}$ is a double covering given by $(M_i,q_i)$
as defined in~\eqref{eq:Miqicover} the classical BNR-correspondence restricted
to the~$X_i$ implies the claim on the trace and the determinant of the~$\Phi_i$.
In the case of a smooth generic fiber it is obvious that the Higgs field is a
map as required by~\eqref{eq:HiggsEifamilies}. In the equisingular
case~\eqref{eq:HiggsEipointwise}
the condition on the range (being $\calE \otimes M_i$) is a consequence of
the local descriptions of the $\Phi_i$ in Section~\ref{sec:localHiggs}, which
in turn follows from the orders of zeros and poles of $\lambda_i$ imposed
by being compatible with $\wh{\Gamma}$.
\par
The level rotation torus acts simultaneously on the abelian differentials
$\lambda_i$ and the prong-matching. This induces an action on
$\Phi_i = \wh{\pi}_*(\,\,\cdot\, \lambda_i)$ and thus on its determinant.
This is obviously compatible with the action on~$q_i = \lambda_i^2$ and
shows that push-forward is well-defined on equivalence classes.
\par
The claim about stability will be proven separately in
Proposition~\ref{prop:corr_stabl}.
\end{proof}

\section{Proof of the spectral correspondence}

We will first define an inverse map to the pushforward correspondence
considered above. In the case of a family of curves $\calX \rightarrow S$
with generically smooth fiber over a smooth base scheme $S$ this is an
easy task: We can apply the classical case as in the proof of
Theorem~\ref{thm:Hitchinfiber} over 
$\wh{\calX}\setminus \wh{N}$ where $\wh{N}$ is the nodal locus in the
fibers of the family and then push forward along the inclusion $j: \wh{\calX}\setminus
\wh{N} \rightarrow \wh{\calX}$. A flat family of torsion-free sheaves is
reflexive (see \cite[page 191]{NagSesh}) and hence determined by its values on
a codimension 2 subset (\cite[Proposition~1.4]{Hart80}).
\par
However, this argument does not work in the equi-singular case as here the nodes
have codimension~$1$. Here we need to give an explicit construction. To emphasize
the functoriality of the construction we instead give an explicit construction
using the special frames constructed in Subsection \ref{sec:localHiggs} in both cases.
\par
The second task in this section is to verify stability conditions in
this correspondence.

\subsection{The pullback correspondence over a generically smooth family of curves}

Let $f: {\calX} \to S$ be a germ of a family of $4g-4$-pointed stable curves.
Let $(\wh{\Gamma}, \calE, \boldsymbol{\Phi}, \bftau)$ be a germ of a family
of multi-scale $\GL(2,\CC)^\circ$-Higgs pairs on $\calX$. First we want to
recover a locally free sheaf $\calF_i$ on $\calX_{i^\complement}$. To do so we
apply the classical pullback correspondence from Theorem~\ref{thm:Hitchinfiber}
to  $(\calE, \bfPhi)$. For each level~$i$ let 
\bes \wh{\calB}_i \=\frac12 \div(\lambda_i) \in \mathsf{Div}^+(\wh{\calX}_{i^\complement})
\ees
be the zero divisor of the family of abelian differentials on
$\calX_{i^\complement}$. We define 
\bes \calF_i\,:=\,\ker(\wh{\pi}^*_i\Phi_i-\lambda_i \Id_{\wh{\pi}^*\calE}) \otimes \calO(\wh{\calB}_i)\,,
\ees
which is a locally free sheaf on $\wh{\calX}_{i^\complement}$. These locally free sheaves
naturally glue on the smooth fibers of $\wh{\calX} \to S$, where the
differentials $\lambda_i$ differ by an invertible function on $S$. This defines a locally free sheaf $\calF'$ on $\wh{\calX} \setminus \wh{N}$. We will now describe how to extend this locally free sheaf over the nodes $\wh{N}$ with respect to local frames. This will define the torsion-free rank 1 sheaf $\calF$ on $\wh{\calX}$.
\par
In the local considerations we can restrict to two levels $0,-1$ meeting in one node~$e$ in the special fiber $X$. (The case of a horizontal node is dealt with in the same
way, necessarily falling into the subcases 'swapped node'. There, in the following local computations, the indices $i=0$ and $i=1$ that usually refer to the level should be used to enumerate the two branches of the node.) We have to consider several special cases in parallel to Section~\ref{sec:localHiggs}. As we mentioned above, whenever the family of $\wh{\calX}$ is normal the sheaf constructed below is equal to the reflexive sheaf $j_*\calF'$. In the following for all modules over a ring $\wh{R}$ that are defined as tensor product we act by $\wh{R}$ on the left factor.
\par \medskip
\subsubsection{\textbf{Fixed node and $\calE$ has a free summand locally near the node}}
This is to say that we start with $\calE = R \oplus \langle x,y \rangle_R$
as in~(LF/fix) and the Higgs fields as in~\eqref{eq:Phi01fixed}. In this case
\bes
\ker\left( \wh{\pi}^*\Phi_0 -\lambda_0 \Id\right)
\= \langle 1 \otimes 1 +\frac{1}{x} \otimes x \rangle \quad \subset \quad
\calO(\wh{X}_1) \otimes \wh{\pi}^*\calE|_{\wh{\calX}\setminus{\wh{X}_1}}=\wh{R} \otimes \wh{R}|_{\wh{\calX}\setminus{\wh{X}_1}}
\ees
and similarly
\bes
\ker\left( \wh{\pi}^*\Phi_1 -\lambda_1 \Id\right)
\= \langle 1 \otimes 1 +\frac{1}{y} \otimes y \rangle \quad \subset \quad
\calO(\wh{X}_0) \otimes\wh{\pi}^*\calE|_{\wh{\calX}\setminus{\wh{X}_0}}=\wh{R} \otimes \wh{R}|_{\wh{\calX}\setminus{\wh{X}_0}}
\ees
These two sheaves glue to a locally free sheaf $\calF'$ on the complement
$\calX \setminus \{\wh{e}\}$ of the node. More precisely, the generating sections glue on the set $\{t\neq0\}$ by
\[ 1 \otimes 1 +\frac{t^\alpha y^2}{t^\alpha y^2}\frac{1}{x} \otimes x \= 1 \otimes 1 + \frac{1}{y} \otimes y\,.
\] Hence $\calF'$ is free of rank 1 over $\mathcal{X}\setminus\{\wh{e}\}$
and consequently $\calF'$ extends over the node as a free rank~$1$ module
$\wh{R}$ over~$\calX$. This defines~$\calF$ at~$\wh{e}$.
\par
\medskip
\subsubsection{\textbf{Fixed node and $\calE$ has no free summand}} \label{sssec:fix_node_free_summand}
We now start with~$\calE$ as in~(NS/fix). Now 
\bes
\ker\left( \wh{\pi}^*\Phi_0 -\lambda_0 \Id\right)
= \langle 1 \otimes (x-a) + \tfrac{1}{x} \otimes x(x-a) \rangle  \subset 
\calO(\wh{X}_1) \otimes \wh{\pi}^*\calE|_{\wh{\calX}\setminus{\wh{X}_1}}=\wh{R} \otimes \wh{R}|_{\wh{\calX}\setminus{\wh{X}_1}}
\ees
and similarly
\bes
\ker\left( \wh{\pi}^*\Phi_1 -\lambda_1 \Id\right)
= \langle 1 \otimes (y-b) + \tfrac{1}{y} \otimes y(y-b) \rangle \subset 
\calO(\wh{X}_0) \otimes \wh{\pi}^*\calE|_{\wh{\calX}\setminus{\wh{X}_0}}=\wh{R} \otimes \wh{R}|_{\wh{\calX}\setminus{\wh{X}_0}}
\ees
The two generators $e_0, e_1$ satisfy the relations $y e_0 = -a e_1$ and $xe_1=-be_0$ for $t \neq0$, where they are both defined. In particular, they glue to a locally free sheaf $\calF'$ on $\calX \setminus \{\wh{e}\}$. Notice that there is an isomorphism 
\[ \calI|_{\wh{\calX} \setminus \{ \wh{e} \}} \rightarrow \calF',
\quad s_1 \mapsto e_0, s_2 \mapsto e_1 .
\]
Hence, we can extend the module $\calF'$ over $\wh{e}$ by $\calI$. This defines $\calF$ at $\wh{e}$. 
\par
\medskip
\subsubsection{\textbf{Swapped node and $\calE$ is locally free}}
With $\calE$ written in a frame as in~(LF/swap) and the Higgs field
from~\eqref{eq:Higgsswap} the kernel $\ker\left( \wh{\pi}^*\Phi_0
-\lambda_0 \Id\right)$ is generated (on $\calX \setminus X_1$) by the
line bundle generated by the first coordinate on $U_1$ and by
the second coordinate on $U_2$. The same holds for the kernel
$\ker\left( \wh{\pi}^*\Phi_1 -\lambda_1 \Id\right)$. Consequently these
coordinate line bundles extend $\calF'$ over both nodes $\wh{e}_1,\wh{e}_2$ to a line bundle~$\calF$ on $\calX$.
\par
\medskip
\subsubsection{\textbf{Swapped node and $\calE$ is not locally free}}
We again only treat the case, where $\calE_e \cong \frakm_e \oplus R$. The case of $\calE_e\cong \frakm_e \oplus \frakm_e$ works along the same lines. With $\calE$ written in a frame as in~(NS/swap) and the Higgs field
still as in ~\eqref{eq:Higgsswap}, we define $\calF$ on the complement
of $X_1$ and $X_0$ with the help of the eigenspaces
\[ \ker( \wh{\pi}^*\Phi_{i=0}-\lambda_{i=0} \Id_{\wh{\pi}^*\calE} |_{x_j \neq 0} ), \quad
\text{resp.} \quad \ker( \wh{\pi}^*\Phi_{i=1}-\lambda_{i=1} \Id_{\wh{\pi}^*\calE} |_{y_j \neq 0} ).
\] Both these kernel are locally free generated by $s_1=(x_1-a,1)$ and
$s_2=(y_1-b,1)$ respectively. These glue over $t \neq 0$ by 
\[ s_1 \= \Big(-\frac{a}{y_1},1\Big) \, s_2 \quad \Leftrightarrow
\quad s_2\= \Big(-\frac{b}{x_1},1\Big)\, s_1
\] to a locally free sheaf~$\calF'$ defined away from the nodes $\wh{e}_1,\wh{e}_2$. The restriction to the first coordinate $U_1$ can be extended by $\calI_1$ over the node. The restriction to the second coordinate $U_2$ extends by the free module $R_2$. This defines $\calF$.
\par

\subsection{The pullback correspondence in the equisingular case}
For notational simplicity we consider the pullback correspondence for multi-scale
$\GL(2,\CC)^\circ$-Higgs pairs $(\calE,\bfPhi)$ on a single stable curve $X$. In
this case the levelwise abelian differentials $\lambda_i$ can be interpreted as
sections of $\wh{\pi}_i^*M_i$ having only simple zeros. Denote by $\wh{B}_i=\div(\lambda_i)$ their divisors. On each level we define the locally free sheaf $\calF_i=\ker(\wh{\pi}_i\Phi_i-\lambda_i\Id_{\wh{\pi}_i^{*}\calE}) \otimes \calO(\wh{B}_i)$. We have to glue this level-wise locally free sheaves to define a torsion-free sheaf $\calF$ on $\wh{X}$. This gluing will again be defined with respect to the special frames constructed in \ref{sec:localHiggs}. It will become apparent that the construction given here agrees with the construction of the previous section restricted to the special fiber. 
\par \medskip
\subsubsection{\textbf{Swapped nodes}} The argument in the {\bf swapped node
cases} carries over to equisingular families without any change. Let us give
some details in the {\bf swapped node, Neveu-Schwarz}-case. To define~$\calF$
we choose a frame of~$\calE$ diagonalizing~$\Phi$. Then 
\[ \wh{\pi}^*(\calE,\Phi)|_{U_1} \=  \left(
\calI_1 \oplus \sigma^*\wh{R}_2, 
\begin{pmatrix} \lambda|_{U_1} & 0 \\ 0 & - \lambda|_{U_1} 
\end{pmatrix} 	\right).
\]
The generator $s_1$ of $\calI_1$ generates the eigenspace
$\ker( \wh{\pi}^*\Phi_{i=0}-\lambda_{i=0} \Id_{\wh{\pi}^*\calE} |_{y_j = 0} )|_{U_1}$
freely.
Similarly, the generator $s_2$ of $\calI_1$ freely generates (still on
the open set $U_1$) the eigenspace $\ker( \wh{\pi}^*\Phi_{i=1}-\lambda_{i=1}
\Id_{\wh{\pi}^*\calE} |_{x_j = 0} )|_{U_1}$. These two generators glue back
inside $\wh{\pi}^*\calE$ to $\calI_1$. In the same way, one recovers $\wh{R}_2$
on $U_2$. This defines $\calF$ at $\wh{e}_1,\wh{e}_2$.
\par \medskip
\subsubsection{\textbf{Fixed node, $\calE$ has a free summand}}
We are in the situation of~\eqref{eq:covfix} specialized to $t=0$. The
level-wise locally free rank 1 sheaves are given by
\[ \calF_0 \= \ker(\wh{\pi}^*\Phi_0- \lambda_0 \Id_{\wh{\pi}^*\calE})
\otimes \calO(\wh{B}_0), \quad \calF_1\= \ker(\wh{\pi}^*\Phi_1- \lambda_1
\Id_{\wh{\pi}^*E}) \otimes \calO(\wh{B}_1)
\]
In local coordinates the eigen-sheaves are generated by
\[ s_0=\tfrac{1}{x} \otimes x + 1 \otimes 1 \in \Gamma(X_0,\wh{\pi}^*
\calE|_{X_0}\otimes \calO(\wh{e}_0)), \ s_1= \tfrac{1}{y} \otimes y + 1 \otimes 1 \in \Gamma(X_1,\wh{\pi}^*\calE|_{X_1}\otimes \calO(\wh{e}_1)).
\]
We define a local generator for $\calF$ by gluing these two generators.  
\par
There is finite number of choices of coordinates $x,y,u,v$, such the covering has above form. They differ by multiplication with constants in $\CC^\times$. This does not affect the generators $s_0,s_1$. However, we work with respect to a special frame of $(\calE,\bf\Phi)$ here and this choice does affect the definition of the generator of $\calF$. The following lemma proves that the choices of special frames are in one-to-one correspondence to choices of frames of the locally free sheaf $\calF$. In particular, $\calF$ is well-defined.

\begin{lemma}\label{prop:uni_special_frame} Fix local coordinates $x,y,u,v$,
such that the covering is given by~\eqref{eq:covfix} with $t=0$.
The choices of special frames in (\ref{eq:LF/fixframe}) such that the Higgs field
is given by~\eqref{eq:Phi01fixed} are in one-to-one correspondence with the choices
of frames $\calF \cong \hat{R}$ under the spectral correspondence.
\end{lemma}
\begin{proof} It is easy to see using the Higgs field that a special frame as in (\ref{eq:LF/fixframe}) is uniquely determined by choosing a generator of the locally free part $s_1$.  Let $\wh{1}$ denote the background generator for $\calF$. Choosing the generator $\phi \cdot \wh{1}$ with $\phi \in \wh{R}^\times$ results in the new generator $t_1s_1$ with $t_1=\phi_1 + \phi_2 x + \phi_3 y$ with $\phi_1$ the even part of $\phi$ and $\phi_2x, \phi_3y$ the odd part in the $x$ respectively $y$ coordinate. 
For the converse, choose a generator of the locally free part $t_1s_1$ with $t_1=e+fx+gy$ with respect to the background frame for $e \in R^\times$, $f \in \CC[[u]]^\times$ and $g \in \CC[[v]]^\times$. It is easy to compute that the new frame of $\calF$ under the spectral correspondence will be $(\pi^\sharp(e) + \pi^\sharp(f)x + \pi^\sharp(g)y) \cdot \wh{1}$. Every element of $\wh{R}^\times$ can be described in this way. 
\end{proof}
\par
\subsubsection{\textbf{Fixed node and $\calE$ has no free summand}}
We are in situation of \eqref{eq:NS/fixframe} specialized to $t=0$. We think of the pullback $\wh{\pi}^*\calE$ as the left $\wh{R}$-module $ \wh{R} \otimes_R \wh{R}$. In this representation, the pullback of the Higgs field $\Phi_i$ acts by multiplying with $\lambda_i$ from the right. To recover $\calI$ we first compute the eigensheaf of the Higgs fields. As in the previous case, this will not quite define $\calI$ due to the twisting by the pushforward at a branch point (cf. proof of Theorem \ref{thm:Hitchinfiber}). 
Let 
\begin{align*} 
	&e_1:=x \otimes (x-a) + 1 \otimes x(x-a), \quad e_2:= b(x	 \otimes 1 + 1 \otimes x) \\
 	&e_3: = y \otimes (y-b) + 1 \otimes y(y-b), \quad  e_4:=a(y \otimes 1 + 1 \otimes y).
\end{align*} It is easy to check that 
\[ e_1,e_2\in \ker( \wh{\pi}^*\Phi_0-\lambda_0 \Id_{\wh{\pi}^*\calE_0})
\quad \text{and} \quad e_3,e_4 \in \ker( \wh{\pi}^*\Phi_1-\lambda_1
\Id_{\wh{\pi}^*\calE_1}).
\] These sections generate the restricted eigensheaves. However, this does not
immediately give us generators for the restrictions of $\calI$. More precisely,
for $b \neq 0,a=0$ the section~$e_2$ generates $\calI_0(-\wh{e}_0)$ and $e_3$
generates $\calI_1(-\wh{e}_1)$. For $a \neq 0,b=0$ the section $e_1$ generates $\calI_0(-\wh{e}_0)$ and $e_4$ generates $\calI_1(-\wh{e}_1)$. Finally, for $a=b=0$ the
section $e_1$ generates $\calI_0(-\wh{e}_0)$ and $e_3$ generates $\calI_1(-\wh{e}_1)$.
Again, with have to glue sections with simple poles at the preimages of the node
to obtain global generators $s_1,s_2$ of $\calI$: 
\begin{itemize}
\item[i)] To obtain $s_1$ we glue $\frac{1}{x}e_1$ on $\wh{X}_0$ to $-\frac{1}{y}e_4=-a(\frac{1}{y} \otimes y + 1 \otimes 1)$ on $\wh{X}_1$, 
\item[ii)] To obtain $s_2$ we glue $\frac{1}{y} e_3$ on $\wh{X}_1$ to $-\frac{1}{x}e_2=-b(\frac{1}{x} \otimes x + 1 \otimes 1)$ on $\wh{X}_0$. 
\end{itemize}

These are the restrictions of the generators used in
Section~\ref{sssec:fix_node_free_summand} to the special fiber. By definition
the resulting generators satisfy the relations
\[ ys_1\=-as_2, \qquad xs_2\=-bs_1\,.
\] We define $\mathcal{F}$ by gluing $\calF_0$ to $\calF_1$ by the $\wh{R}$-module $\calI$ generated by these sections. Again this construction depends on the choice of a special frame of $(\calE,\bfPhi)$. The choice of such a frame is equivalent to a choice of generators for $\calI$.

\begin{lemma} The choices of special frames described by \eqref{eq:NS/fixframe}
such that the Higgs fields have the prescribed form are in one-to-one correspondence
to the choices of generating sets $s_1,s_2 \in \calI$, such that $ys_1=-as_2, xs_2=-bs_1$.
\end{lemma}
\begin{proof}
First one shows that all generators of $\calI$ satisfying the relations are given by $\phi s_1,\phi s_2$ with $\phi \in \wh{R}^\times$. We can decompose $\phi=\phi_1+ \phi_2x+\phi_3y$ with $\phi_i \in R$. Then the corresponding generators of $\pi_*\calI$ are 
\begin{align*} 
	t_1'&=\phi_1(x-a)+\phi_2 x(x-a) - \phi_3ay, \quad t_3'=xt_1', \\
 	t_2'&=\phi_1(y-b)-\phi_2bx+\phi_3 y(y-b), \quad t_4'=yt_2'.
\end{align*}  Let $t_1=x-a,t_2=y-b, t_3=x(x-a), t_4=y(y-b)$ be the standard
generators of the special frame in \eqref{eq:NS/fixframe}. It is easy to see
that any generators $t_1',t_2',t_3',t_4'$ of $\calE$, such that the Higgs field
has the desired form are uniquely determined by the choice of $t_1',t_2'$.
Now an easy but tedious computation using the relations of $\calE$ shows that
all generating sets are of the form described above. 
\end{proof}
\par
\medskip

\subsection{The proof of Theorem \ref{thm:Higgspairsfam} and Theorem \ref{thm:Higgspairs}}

\begin{proof}[Proof of Theorem \ref{thm:Higgspairsfam}]
Let $f: {\calX} \to S$ be a germ of a family of $4g-4$-pointed stable curves. Let $(\wh{\pi}: \wh{\calX} \to \calX, \bfq, \bftau, \calF)$ be a germ of families of $\GL(2,\CC)^\circ$-spectral data. Let $(\calE,\bf{\Phi})$ be the associated germ of families of multi-scale $\GL(2,\CC)^\circ$-Higgs pairs defined by pushforward. For each level we recover the restrictions of $\calF$ to $\calX_{i^\complement}$ by $\ker(\wh{\pi}^*_i\Phi_i-\lambda_i \Id_{\wh{\pi}^*\calE}) \otimes \calO(\wh{\calB}_i)$. For each fiber of $f$ this is the classical spectral correspondence revisited in Theorem \ref{thm:Hitchinfiber}. These eigensheaves glue over $t \neq 0$ to a locally free sheaf $\calF'$ on $\wh{\calX} \setminus \wh{N}$, where $\wh{N}$ is the set of nodes of the singular fibers of $\wh{\calX} \to S$. Then by the local computations above $\calF\cong j_*\calF'$. More precisely, the choice of a local frame of $\calF$ at the nodes of $\wh{X}$ determines special frames of $(\calE,\bfPhi)$ at the nodes on $X$, which in turn determine special frames of $j_*\calF'$ by the local computation above. It is with respect to this special frames of $\calF$ and $j_*\calF'$ that the isomorphism extends over the nodes of $\wh{X}$. 
\par
For the converse, start with a germ of families of multi-scale $\GL(2,\CC)^\circ$-Higgs pairs $(\wh{\Gamma}, \calE, \boldsymbol{\Phi}, \bftau)$ on $\calX$. Then the pullback construction yields a torsion-free sheaf $\calF=j_*\calF'$, such that by construction there is an isomorphism of $\wh{\pi}_*(\calF,\bflambda) \cong (\calE, \bf{\Phi})$ on $\calX\setminus N$, where $N$ is the set of nodes of the singular fibers of $\calX \to S$. Again this follows from the classical spectral correspondence applied to the fibers of $f$. The choice of a special frame at each node $e \in \calX$ determines a frame of $\calF$ at the preimages of the node. This in turn determines a special frame of $\wh{\pi}_*(\calF,\bflambda)$. By the local spectral correspondences the isomorphism on $\calX\setminus N$ extends over the nodes with respect to the special frames to an isomorphism $(\calE, \bf\Phi)\cong\wh{\pi}_*(\calF,\bflambda)$. That the spectral correspondence preserves the notions of stability will be proven separately in Proposition \ref{prop:corr_stabl}.
\end{proof}

\begin{proof}[Proof of Theorem \ref{thm:Higgspairs}]
Let $(\wh{\pi}: \wh{X} \to X, \bfq, \bftau, \calF)$ be a multi-scale $\GL(2,\CC)^\circ$-spectral datum. Let $(\calE,\bf{\Phi})$ be the associated multi-scale $\GL(2,\CC)^\circ$-Higgs pairs defined by pushforward. By Theorem \ref{thm:Hitchinfiber} applied to the case of $M_i$-twisted Higgs pairs on $X_i$ we recover the restrictions $\calF_i$ as
\[ \ker(\Phi_i-\lambda_i\Id_{\pi_i^{*}E}) \otimes \calO(\wh{B}_i).
\]
The gluing of the restrictions was described locally at all nodes. This determines a torsion-free sheaf $\calF'$ on $\wh{X}$. We claim that $\calF \cong \calF'$. This can again be checked with respect to frames. Choices of frames of $\calF$ at $\wh{N}$ determine special frames of $(\calE,\bfPhi)$ at $N$. These in turn induce frames of $\calF'$ at $\wh{N}$ by the local construction. With respect to such a framing the componentwise isomorphism $\calF'_i\cong\calF_i$ extends over the nodes.
\par
For the converse, consider a multi-scale $\GL(2,\CC)^\circ$-Higgs pair
$(\calE,\mathbf{\Phi})$ on $X$. The level-wise eigen-sheaves glue to a torsion-free
sheaf $\calF$ as described above. The pushforward $\wh{\pi}_*(\calF,\bf{\lambda})$ is
a multi-scale $\GL(2,\CC)^\circ$-Higgs pair, such that by the classical spectral
correspondence the restrictions to the levels are isomorphic to the restriction of $(\calE,\bf{\Phi})$. Choices of special frames of $\calE$ at the nodes correspond to choices of frames of $\calF$ at their preimages. Hence, the local considerations above show that the levelwise isomorphisms extend with respect to special frames, i.e.\, $\wh{\pi}_*(\calF,\bf{\lambda}) \cong (\calE,\bf{\Phi})$. That the spectral correspondence preserves the notions of stability will be proven separately in Proposition \ref{prop:corr_stabl}.
\end{proof}
\par
\begin{proposition}\label{prop:corr_stabl} Let $P$ be a polarization of multi-scale $\GL(2,\CC)$-Higgs pairs of degree $d$ on $X$. A multi-scale $\GL(2,\CC)^\circ$-Higgs pair $(\calE,\bfPhi)$ is $P$-(semi)stable if and only if the associated torsion-free
sheaf $\calF$ is $\wh{P}$-(semi)stable.
\end{proposition}
\par
\begin{proof}
As preparation we show how to translate the stability condition \eqref{eq:GL(2)stab} into a stability condition on torsion free quotients of $\calE$. Let $\calE'$ be a $\bfPhi$-invariant subsheaf of $\calE$. Then $\calE'$ defines three (closed) subcurves $Y_0,Y_1,Y_2 \subset Y$, such that $\calE'|_{Y_i}$ has generic rank $i$. On all irreducible components $X_v$
where $q_v$ is not a square of an abelian differential, there are no $\Phi_{X_v}$-invariant rank one subsheaves of $\calE|_{X_v}$. Hence, for all $X_v \subset Y_1$ the meromorphic quadratic differential $q_v$ is the square of an abelian differential. Consider the exact sequence
\bes 0 \rightarrow \calE' \rightarrow \calE \rightarrow \calQ \rightarrow 0.
\ees 
Then $\calE'$ satisfies \eqref{eq:GL(2)stab} if and only if $\calQ$ satisfies the inequality 
\be \chi(\calQ)\geq \chi(\calE)-\frac{\sum_{v=1}^s \rk\calG_{X_v}\deg P_{X_v}}{2\rk P} = \frac{\deg P_{Y_1}+2\deg P_{Y_0}}{2\rk(P)}. 
\label{eq:stab_quot}
\ee
A priori the quotient $\calQ$ has torsion. We did not put any saturation condition on the subsheaf $\calE' \subset \calE$ hence the support of the torsion submodule $\mathsf{Tor}(\calQ)$ can be any finite collection of points in $X$. Recall that for a coherent sheaf $\calS$ on $X$ and a proper subcurve $U \subset X$ we denote by $\calS_U$ the torsion-free pullback to the subcurve $U \subset X$. By definition of the subcurves $Y_j$ we have $\calQ_{Y_0 \cup Y_1}=\calQ/\mathsf{Tor}(\calQ)$. The multi-scale Higgs fields $\bfPhi$ induces levelwise Higgs fields on $\calQ$ as $\calE'$ was $\bfPhi$-invariant. These Higgs fields preserve the torsion submodule. Hence the kernel of $\calE \rightarrow \calQ_{Y_0 \cup Y_1}$ is again preserved by $\bfPhi$. This kernel has Euler characteristic $\geq \chi(\calE')$. Therefore it is enough to check the inequality \eqref{eq:stab_quot} for torsion-free quotients $\calQ$ with Higgs fields that descend level-wise.
\par 
So let $\calQ$ be a torsion-free quotient of $\calE$, such that the multi-scale Higgs field $\bfPhi$ descends. We claim that under these conditions there exists a proper subcurve $\wh{W} \subset \wh{X}$, such that $\wh{\pi}_*\calF_{\wh{W}}=\calQ$.
To prove this claim first notice that pushforward along $\wh{\pi}$ commutes with torsion-free pullback to subcurves and hence $\calE_{Y_0 \cup Y_1}=(\wh{\pi}_*\calF)_{Y_0 \cup Y_1}=\wh{\pi}_*(\calF_{\wh{\pi}^{-1}(Y_0 \cup Y_1)})$. If $Y_1=\emptyset$, we have $\calQ=\calE_{Y_0}$ and this proves the claim. If $Y_1 \neq \emptyset$, there is an exact sequence
\bes 0 \rightarrow \calL \rightarrow \calE_{Y_0 \cup Y_1}=\wh{\pi}_*(\calF_{\pi^{-1}(Y_0 \cup Y_1)}) \rightarrow \calQ \rightarrow 0\,,
\ees where $\calL$ is torsion-free rank 1 sheaf supported on $Y_1$. However, by the initial
remark, for all $X_v \subset Y_1$ there exist abelian differentials $a_v$ and line bundles $\calL_1,\calL_2$, such that $q_v=-a_v^2$ and $(\calE,\bfPhi)_{X_v}=(\calL_1  \oplus \calL_2,  \mathrm{diag}(a_v,-a_v))$. Let $\wh{\pi}^{-1}(X_v)=\wh{X}_{v_1}\cup \wh{X}_{v_2}$. By definition
$\calL_j=\wh{\pi}_*\calF|_{\wh{X}_{v_j}}$ for $j=1,2$. Being $\Phi$-invariant the kernel $\calL$ restricted to $X_v$ must agree with one of the $\calL_j$. Hence, $\calQ =\wh{\pi}_*\calF_{\wh{W}}$ for a proper subcurve $\wh{W} \subset \wh{\pi}^{-1}(Y_0 \cup Y_1) \subset \wh{X}$ defined by choosing one of the preimages of $\pi^{-1}X_v$ for all $X_v \in Y_1$. This proves the claim in general. 
\par
As the Euler characteristic is invariant under pushforward along finite maps the above inequality \eqref{eq:stab_quot} is equivalent to 
\bes
\chi(\calF|_{\wh{W}})\geq \frac{\deg P_{Y_1}+2\deg P_{Y_0}}{2\rk P}=\frac{\deg \wh{P}_{\wh{W}}}{\rk \wh{P}}.
\ees That is the stability condition of \eqref{eq:phistabEC} with respect to $\wh{P}$.
\end{proof}

\subsection{Tschirnhausen modules and multi-scale Higgs pairs}
\label{ssec:Tschirnhausen}

Fix a $\GL(2,\CC)^\circ$-multi-scale spectral datum $(\wh{\pi}:
\wh{X} \to X, \bfq, \bftau, \calF)$. The pushforward of the structure sheaf
of the admissible cover $\wh{\pi}: \wh{X} \rightarrow X$ defines an
$\calO_X$-algebra $\calA:=\wh{\pi}_*\calO_{\wh{X}}$. One each irreducible
component $X_{v}$ of $X$ the covering $\wh{\pi}$ is defined by an
equation $\lambda^2-q_v$, hence it is affine.
By \cite[Exercise II.5.17]{hartshorne} the pushforward induces an equivalence
of categories between the category of quasi-coherent sheaves on $\wh{X}$ and
the category of quasi-coherent $\calA$-modules on $X$.  In particular,
$\calE=\wh{\pi}_*\calF$ has an $\calA$-module structure. In the following we
want to compare this $\calA$-module structure to a multi-scale Higgs
pair $(\calE,\boldsymbol{\Phi})$ compatible with $\bfq$ (cf.\
Definition~\ref{def:Higgspair}). 
\par
The $\calA$-module structure is decoded in the action of the Tschirnhausen module
on $\calE$. The Tschirnhausen module $\tau$ is defined through 
\[ \calA\=\wh{\pi}_*\calO_{\wh{X}}\= \calO_X \oplus \tau\,.
\]
An $\calA$-module structure on $\calE$ induces a morphism
of $\calO_X$-modules
\[ \varphi: \tau \otimes \calE \rightarrow \calE
\] satisfying the relations of the sheaf of $\calO_X$-algebras $\calA$. 
On each level~$X_i$ the covering is defined by the projective class of
differentials $[q_i] \in \PP H^0(X_i,M_i^2)$ with $M_i$ defined in~\eqref{eq:Micoverbundle}.
This identifies the restriction $\tau|_{X_i}=M_i^{-1}$ (see for
example \cite[Corollary 4.11]{Horn}). Choosing a representative~$q_i$ is equivalent to
choosing an embedding $X_i \subset \mathsf{Tot}(M_i)$. By the classical
BNR-correspondence such a choice induces a Higgs field 
\[ \Phi_i: \calE_i \rightarrow \calE_i \otimes M_i,
\] such that $\Phi_i^2+\Id_{\calE} \otimes q_i=0$.
We will show the converse: A choice of levelwise
Higgs fields compatible with the quadratic multi-scale  differential~$\bfq$
induces an $\calA$-module structure on $\calE$. A local description of $\tau$
is apparent from the construction of the pushforward in
Section~\ref{sec:localHiggs}. Thereby $\tau$ is locally free at each node
where $\kappa$ is even, and isomorphic to the maximal ideal $\frakm$ at
each node where $\kappa$ is odd.  
\par
\begin{lemma} Let $X'$ be the partial normalization of $X$ at all nodes~$e$
with $\kappa_e$ being an odd number. Let $(\wh{\Gamma}, \calE,
\boldsymbol{\Phi}, \bftau)$ be a multi-scale $\GL(2,\CC)^\circ$-Higgs pair
with respect to $\bfq$. Define $\varphi_i: M_i^{-1} \otimes \calE_i \rightarrow \calE_i$
to be the morphism of $\calO_{X_i}$-sheaves obtain form $\Phi_i$ by tensoring with
$M_i^{-1}$. A morphism of $\calO_{X'}$-modules
$\varphi_{X'}: \tau|_{X'} \otimes \calE|_{X'} \to \calE|_{X'}$, such that the
restriction to $X_i$ is given by $\varphi_i$ defines a unique morphism of
$\calO_X$-sheaves 
\[ \varphi: \tau \otimes \calE \rightarrow \calE.
\]
\end{lemma}
\begin{proof} Let $\iota:X' \to X$ be the inclusion. By the local description of $\tau$
given above we have $\iota_{*}\tau|_{X'}=\tau$. The pushforward of the restricted
morphism $\tau|_{X'} \otimes \calE|_{X'} \to \calE|_{X'}$ defines a morphism 
\[ f: \tau \otimes \iota_* \calE|_{X'} \to \iota_* \calE|_{X'}.
\] At a node $e$ with $\kappa_e$ odd $\calE$ is locally given by $R \oplus \frakm$ or $\frakm \oplus \frakm$ (see \ref{eq:LF/fixframe}, \ref{eq:NS/fixframe}). Hence there is an exact sequence
\[ 0 \rightarrow \calE \rightarrow \iota_* \calE|_{X'} \rightarrow \calS \rightarrow 0,
\] where $\calS$ is a coherent sheaf supported at the nodes $e$ with $\kappa_e$ odd with stalks of length at most 1. As $\wh{\pi}$ is affine, $\calA$ is flat and so is $\tau$. Hence we obtain a diagram
\[ \begin{tikzcd} 0 \ar[r]& \tau \otimes \calE \ar[r]\ar[dashed,d] & \tau \otimes \iota_* \calE|_{X'} \ar[r]\ar["f",d] & \tau \otimes \calS \ar[r] \ar["0",d] & 0 \\
0 \ar[r]& \calE \ar[r] & \iota_* \calE|_{X'} \ar[r] & \calS \ar[r] & 0
\end{tikzcd}.
\] By the local description of the Higgs field (\ref{eq:Phi01fixed}) the map $f$ descends to the quotient by zero and hence defines a map on the kernels. This is the desired
morphism $\varphi: \tau \otimes \calE \rightarrow \calE$. 
\end{proof}
\par
\begin{proposition}
A multi-scale $\GL(2,\CC)^\circ$-Higgs pair $(\wh{\Gamma}, \calE,\boldsymbol{\Phi}, \bftau)$ with respect to $\bfq$ induces the structure of an $\calA$-module on $\calE$. 
\end{proposition}
\begin{proof} Let $j: X' \rightarrow X$ be the partial normalization at the nodes $e$ with odd $\kappa_e$. Recall that $\tau_{X'}$ is an invertible sheaf. The level-wise Higgs fields $\boldsymbol{\Phi}$ can be glued to a morphism of $\calO_{X'}$-modules $\calE|_{X'} \to \calE|_{X'} \otimes \tau|_{X'}^{-1}$. To see this recall that the level-wise Higgs fields are sections $\Phi_i \in H^0(X_i,\End(\calE_i) \otimes M_i)$. Hence, there is only one choice of gluing for each level-passage. However, $\tau|_{X'}^{-1}$ is locally free at each node of $X'$. Let $U$ be a neighborhood of a node $e$ with $\kappa_e$ even, such that $\tau^{-1}|_U \cong \calO_U$. Then we can define a $\tau^{-1}|_U$-action by letting $1 \in \calO_U$ act by $\Phi_{\ell^+(e)}$ on the higher level and by $\Phi_{\ell^-(e)}$ on the lower. This defines unique morphism $\calE|_{X'} \to \calE|_{X'} \otimes \tau|_{X'}^{-1}$. By the previous lemma it can be extended to a morphism of $\calO_X$-sheaves $\varphi: \tau \otimes \calE \rightarrow \calE$. We are left with showing that this morphism induces an $\wh{\pi}_*\calO_{\wh{X}}$- module structure on $\calE$. 

We can realize $\calA= \Spec\mathsf{Sym}(\tau)/\calI$, where $\calI$ is an $\mathsf{Sym}(\tau)$-ideal sheaf. We need to show that the morphism $\tau \otimes \calE \rightarrow \calE$ satisfies the relation of $\calI$. Restricted to levels this means that the Higgs fields satisfy $\Phi_i^2+\Id_{\calE} \otimes q_i=0$, which is satisfied by construction. At a node $e$ with $\kappa$ even the Tschirnhausen module is the module of odd local functions under $\sigma: \wh{X} \to \wh{X}$. On the other hand, $\lambda_{\ell^+(e)}, \lambda_{\ell^-(e)}$ are odd  with respect to $\sigma$ and non-vanishing in a neighborhood of the node as sections of $\wh{\pi}^*M_{\ell^+(e)}, \wh{\pi}^*M_{\ell^-(e)}$. Hence $\Phi_{\ell^+(e)}, \Phi_{\ell^-(e)}$, which correspond to multiplication by $\lambda_{\ell^+(e)}, \lambda_{\ell^-(e)}$ by definition, induce an $\calA$-module structure on $\calE$ in a neighborhood of the node $e$. At a node~$e$
with $\kappa_e$ odd the module $\calA$ is given by $\calA= 1 \oplus \langle x,y \rangle$
with the relations are $x^2=u, y^2=v, xy=0$. The first two conditions are satisfied
as a consequence of $\Phi_i^2+\Id_{\calE} \otimes q_i=0$.
The third condition is satisfied by extending the levelwise Higgs fields by zero to the other levels. In summary, we showed that the morphism $\varphi'$ defines an $\calA$-module structure on $\calE$.
\end{proof}
\par
\begin{proof}[Alternative proof of Theorem \ref{thm:Higgspairs}]
Let $(\wh{\pi}: \wh{X} \to X, \bfq, \bftau, \calF)$ be a $\GL(2,\CC)^\circ$-spec\-tral datum. Then $\wh{\pi}_*\calF$ has an $\wh{\pi}_*\calO_{\wh{X}}$-module structure. Hence, it corresponds to a unique quasi-coherent sheaf on $\wh{X}$. This recovers $\calF$. 
\par
For the converse, a multi-scale Higgs pair $(\calE,\bf{\Phi})$ induces a $\wh{\pi}_*\calO_{\wh{X}}$-module structure on $\calE$ by the previous proposition. Hence it determines a unique quasi-coherent sheaf $\calF$, such that $\calE=\wh{\pi}_*\calF$. Using the special frames of~$\calE$ given in Section~\ref{sec:localHiggs} at each node we conclude that $\calF$ is torsion-free of rank 1.
\end{proof}

\section{Comparison to the original Hitchin fibers}\label{sec:Compare}

In this section we compare for a smooth curve~$X_\st$ the fibers of the
compactified Jacobian, hence the fibers of the map~$h: \SD_{X_\st} \to \BB_{X_\st}$
over the modified Hitchin base, with the fibers of the Hitchin map~$\Hit$,
i.e. we compare the vertical arrows in Proposition~\ref{intro:modHFvsHF}.
To do so we work on the level of moduli spaces instead of stacks.
For concreteness, we work in this section with the weighted
canonical polarization $\wh{P} = \wh{P}_{\can}$. For other polarizations
the differences of the fibers are similar. This polarization
is indeed induced by the polarization of the pointed stable curves by
the pointed canonical bundle. In particular, the coarse moduli spaces
for the $\wh{P}$-compactified Jacobian functor can be obtained by Simpson's
construction of moduli spaces of semistable sheaves. The points of this
moduli space correspond to Jordan-Hölder equivalence classes of semistable
torsion-free sheaves on~$\wh{X}$.
\par
The fibers of both maps, the Hitchin fibration and our variant using the
universal Jacobian, are
stratified into semi-abelian varieties, but the combinatorics and the
semi-abelian
varieties in the strata are quite different. 
\par
\medskip
\paragraph{\bf The singular fibers of the Hitchin map}
The geometry of the singular Hitchin fibers $\Hit^{-1}(X_\st,q)$, in the case
where $q$ is not a global square, was analyzed in \cite{Horn} 
using Hecke modifications.
There is a stratification
\be \label{eq:HitDecomp}
\Hit^{-1}(X_\st,q) \= \bigcup_{D} \calS_D\,,
\ee
where $D$ runs over effective divisors on~$X_\st$, such that
$\div(q) - 2D$ is also effective. There are exact sequences
\be \label{eq:SDsequence}
0 \to (\CC^*)^{r_1} \times \CC^{r_2} \to \calS_D \to \Pic^{d}(\Sigma^n) \to 0
\ee
where $\Sigma^n$ denotes the normalization of the spectral curve and $r_1,r_2,d$ can be read off from~$q$ and~$D$, see the examples below.
Splitting $\CC = \CC^* \cup \{0\}$ gives the stratification into
semi-abelian varieties. 
\par
\medskip
\paragraph{\bf The fibers of the universal Jacobian} We next describe
the fibers of the modified Hitchin map~$h$ over a multi-scale differential.
It only depends on the underlying semistable curve~$\wh{X}$.
By definition it is the fiber of the universal compactified Jacobian
$\barcJ[g,n][d]$ over~$\wh{X}$. This fiber admits a stratification
\be \label{eq:uniJacDecomp}
h^{-1}(\wh{X},\bfz,\bfq) \= \bigcup_{(\bfd, N)\atop \text{$\phi^\can$-semistable}}
\Pic^{\bfd}(\wh{X}_N)
\ee
by semi-abelian varieties, where $N$ is a subset of the nodes of~$\wh{X}$
and $\wh{X}_N$ is the partial normalization there. The semi-abelian variety
of line bundles with multi-degree~$\bfd$ is denoted by~$\Pic^{\bfd}$. Here the inclusion~ $\Pic^{\bfd}(\wh{X}_N) \rightarrow h^{-1}(\wh{X},\bfz,\bfq)$ is defined by pushforward along the partial normalization $\wh{X}_N \to \wh{X}$.
The notion of semistability is designed so that the right hand side is a
finite union. This stratification is essentially in the paper of
Caporaso \cite{CapoCompJ} or Oda-Seshadri \cite{OdSesh},
see \cite{CapChrist} for in-depth
examples and \cite{MMUV} for a modern overview.
\par
Recall from the construction of the modified Hitchin base and that if the image of $(\wh{X},\bfz,\bfq)$ under the
forgetful map~$b$ is~$q$, then the top level curve
$\wh{X}_0$ is the normalization of the spectral curve associated to $X_\st$ and the top level differential~$q_0=q$.
\par
\begin{proposition} \label{prop:maptoJac}
Suppose that $q$ is not a global square. Each of the strata
of $\Hit^{-1}(X_\st,q)$ and of $h^{-1}(\wh{X},\bfz,\bfq)$ has a map to
$\Jac^d(\wh{X}_0)$. If~$q$ has an even order zero, these maps do not
glue to a global map to $\Jac^d(\wh{X}_0)$, neither for the original
Hitchin fiber, nor for $h^{-1}(\wh{X},\bfz,\bfq)$.
\end{proposition}
\par
\begin{proof} The existence of the map on each stratum is obvious from
the definition in~\eqref{eq:uniJacDecomp} for the universal Jacobian
and from the description of the spectral data after~\eqref{eq:HitDecomp}
for the Hitchin fiber. The non-existence of a map to $\Jac^d(\wh{X}_0)$ on the universal compactified Jacobian is clear from the non-uniqueness of stable multi-degrees in~\eqref{eq:uniJacDecomp}. For the Hitchin fiber the non-existence of this map is stated in~\cite[Example~8.3]{Horn}, see also \cite[Section~5]{GotOli}. (With an analogous argument one shows that the map to $\Jac^d(\wh{X}_0)$ does not extend to the irreducible components of the universal compactified Jacobian.) 
\end{proof}
\par
We now compare the two fibers in typical loci, see Figure~\ref{cap:23together}.
\par
\subsection{Quadratic differentials with one zero of order $m=2k$:
banana curves} \label{sec:banana}

The special case $k=1$ is illustrated in Figure~\ref{cap:23together} left.
In general, there are $2k$ simple zeros of~$\bfq$ on bottom level~$X_{-1}$ of the pointed
stable curve~$X$ and the genera of the covering curves are $\wh{g}_0
= g(\wh{X}_0) =4g-3-k$ and $\wh{g}_1= g(\wh{X}_{-1}) = k-1$ respectively.
\par
\begin{proposition} \label{prop:comparebanana} 
For $\gcd(d+2g-2,6g-6)=1$, the fibers $\Hit^{-1}(X_\st,q)$ and $h^{-1}(\wh{X},\bfz,\bfq)$ have
different numbers of irreducible components for each $k \geq 1$. For $d = 2g-2 \mod 6g-6$ they have the same number of irreducible components for all $k \geq 1$.
\par
For the generic point in the singular locus of the Hitchin base, the
case of one double zero (i.e.\ $k=1$) the strata of~\eqref{eq:uniJacDecomp}
and of~\eqref{eq:HitDecomp} of the same dimension are isomorphic.
\par
For $k \geq 2$ not even the top-dimensional strata are isomorphic. 
\end{proposition}
\par
In fact, associated to a singularity of the spectral curve one defines
an associated moduli space $\mathrm{Heck}(q)$ of allowable Hecke modifications
and one obtains a map from a fiber bundle of spectral data $\calS$ 
\[ \mathrm{Heck}(q) \rightarrow \calS \rightarrow \Jac(\wh{X}_0)
\] to the $\Hit^{-1}(X_\st,q)$ by applying Hecke modifications. We have
$\mathrm{Heck}(q)\cong\PP^1$ for $k=1$  while $\mathrm{Heck}(q) \cong \PP(1,1,2)$
is a weighted projective space for $k=2$. The map to the Hitchin fiber is birational but no isomorphism. E.g.\ for $k=1$ the $\PP^1$-bundle is
twisted by gluing the zero section over $L \in \Jac^d(\wh{X}_0)$ to the $\infty$-section
over $L(p_+-q_+)$, where $p_+,q_+$ are the preimages of the singularity in $\wh{X}_0=\Sigma^n$.
\par
As preparation for the proof we list the possible multi-degrees~$\bfd$
appearing in~\eqref{eq:uniJacDecomp}, see \cite[Example 7.3]{CapoCompJ}.
If the polarization is non-degenerate (see Theorem~\ref{thm:compJ}
when this holds for the weighted canonical polarization)
the possible multi-degrees of the line bundles in $h^{-1}(\wh{X},\bfz,\bfq)$
are a shift of $(1,0),(0,1)$, otherwise a shift of $(2,0),(1,1),(0,2)$. In the first case, both multi-degrees are stable and correspond to irreducible components of $h^{-1}(\wh{X},\bfz,\bfq)$. In the second, only middle multi-degree is stable and corresponds to an irreducible component. The other two are strictly
semistable. By the identification of numerical stability conditions with the GIT-stability conditions \cite{Alexeev} (see also \cite[Fact 2.8]{CMKV}) these multi-degrees correspond to strata in the boundary of this irreducible component. So in this case, $h^{-1}(\wh{X},\bfz,\bfq)$ is irreducible. 
Concretely, for the canonical pointed
numerical polarization the torsion-free pullbacks to the components must
satisfy the condition 
\begin{align*}
\deg(\calF_{\wh{X}_0}) \geq \wh{d}+\frac{k-2}{3}+\wh{d}\frac{1-2k}{6g-6}, \qquad \deg(\calF_{\wh{X}_{-1}}) \geq \wh{d}\frac{2k-1}{6g-6}-\frac{k+4}{3}\,.
\end{align*}
\par
\medskip
\begin{proof} The irreducibility of the Hitchin fibers with irreducible spectral
curve and at least one zero of odd order is proven in \cite[Corollary~8.5]{Horn}.
In fact the stratum~$\calS_0$ is dense. The count above shows that the fibers
in the compactified Jacobian are reducible for $\gcd(\wh{d},6g-6)=1$ and always
irreducible for $\wh{d} = 4g-4 \mod 6g-6$ for all $k \geq 1$. This proves
the first claim by noting that $\wh{d}=d+2g-2$.
\par
In the $\Hit$-fibers the rank of the abelian part is the genus of the
normalization of the spectral curve, which is the same as the top level
curve~$\wh{X}_0$. In the $h$-fibers the rank of the abelian part is
$g(\wh{X}_0) +g(\wh{X}_1)$. Since for $k \geq 2$ the bottom level curve has
positive genus, this proves the last statement. For $k=1$ the closed strata
are in both cases the Jacobians of~$\wh{X}_0$ and the open strata are
$\CC^*$-extension of this Jacobian. This follows from~\eqref{eq:uniJacDecomp}
and since $r_1=1$ and $r_2=0$ in~\eqref{eq:SDsequence} in this case
by \cite[Theorem~6.2]{Horn}.
\end{proof}
\par

\subsection{Quadratic differentials with one zero of order $m=2k+1$:
compact type} \label{sec:ctcurves}

The special case $k=1$ is illustrated in Figure~\ref{cap:23together} right.
In general, there are $2k+1$ simple zeros of~$\bfq$ on bottom level~$X_{-1}$ of the pointed
stable curve~$X$. The genera of the covering curves are $\wh{g}_0
= g(\wh{X}_0) =4g-3-k$ and $\wh{g}_1= g(\wh{X}_{-1}) = k$.
\par
\begin{proposition} \label{prop:comparect}
The fibers $\Hit^{-1}(X_\st,q)$ and $h^{-1}(\wh{X},\bfz,\bfq)$ are irreducible.
\par
Neither the fibers nor their strata are isomorphic for any $k \geq 1$.
\end{proposition}
\par
In this case, the pointed canonical stability conditions are
\[ \deg(\calF_{\wh{X}_0}) \geq \wh{d}- \wh{d}\frac{k}{3g-3} + \frac{k}{3}, \quad \deg(\calF_{\wh{X}_{-1}}) \geq \wh{d}\frac{k}{3g-3} - \frac{k}{3} -1.
\] 
\par
\begin{proof}
The irreducibility of $\Hit^{-1}(X_\st,q)$ is shown in \cite[Corollary~7.10]{Horn},
see also \cite{GotOli}.
The irreducibility of $h^{-1}(\wh{X},\bfz,\bfq)$ is discussed in \cite[Example~7.1]{CapoCompJ}. There are cases, where there are two strictly semistable
multi-degrees $(d_1+1,d_2)$ and $(d_1,d_2+1)$. However, up to Jordan-Hölder
equivalence they can be represented by the pushforward of a locally free sheaf of multi-degree $(d_1,d_2)$ along the normalization map of $\wh{X}$ (or equivalently by the locally free sheaf on a semistable model with one rational bridge). Hence, the irreducibility still holds.
\par
The strata of the Hitchin fiber in~\eqref{eq:HitDecomp} are semi-abelian
varieties with abelian part $\Jac(\wh{X}_0)$. The
stratification~\eqref{eq:uniJacDecomp} is reduced to a single stratum isomorphic
to $\Jac(\wh{X}_0) \times \Jac(\wh{X}_1)$. Since $g(\wh{X}_{-1}) = k>0$
the remaining claims follow.
\end{proof}
\par
In fact, in this case the fibers of $\Hit^{-1}(X_\st,q)$ over $\Jac^d(\wh{X}_0)$ are isomorphic to
$\mathrm{Heck}(q)\cong\PP^1$ for $k=1$ while $\mathrm{Heck}(q) \cong \PP(1,1,2)$
is a weighted projective space for $k=2$.
\par

\subsection{Square of abelian differential with simple zeros}
In this case the spectral curve has two irreducible components interchanged by $\sigma$. We compare the fibers in
the generic stratum where~$\alpha$ with $\alpha^2 = q$ has only simple zeros.
The level graph and that of the covering are given in Figure~\ref{cap:qtosquare}.
To discuss a concrete example we consider the case $g=3$ given there.
The double cover $\wh{X}$ has two irreducible components on level $0$ of
genus $3$, denoted by $\wh{X}_{01}$,$\wh{X}_{02}$ and four components of
genus $0$ on the lower level $\wh{X}_{11}, \dots,  \wh{X}_{14}$. We have 
\bas \deg( \omega_{\wh{X}}(\wh{\bfz})|_{\wh{X}_{01}})
&\= \deg( \omega_{\wh{X}}(\wh{\bfz})|_{\wh{X}_{02}}) && \=8. \\
\deg( \omega_{\wh{X}}(\wh{\bfz})|_{\wh{X}_{1j}}) &\= 2 \qquad  
&& \text{for all $i=1,\ldots,4$}
\eas
Computing the stability conditions with respect to the $\wh{P}$-polarization of degree $\wh{d}=2g-2$ we obtain: Let $Y \subset X$ be a connected subsurface. If $Y$ contains one component on top level and $0 \leq i \leq 4$ components on bottom level, then $\deg(\calF_Y) \geq \frac{2-i}{3}$. If $Y$ contains the top level and $1 \leq i \leq 3$ components on bottom level, then $\deg(\calF_Y) \geq \frac{2(i-1)}{3}$. If $Y=X_{1j}$, then $\deg(\calF_Y) \geq -\frac{4}{3}$. 
We list all possible multi-degrees up to permutation on the first $2$ or last $4$ entries.
\begin{align*} &(2,2,0,0,0,0),\ (3,1,0,0,0,0),\ (2,3,-1,0,0,0,0),\ (3,3,-1,-1,0,0), \\ & (3,4,-1,-1,-1,0), \ (4,4,-1,-1,-1,-1), \ (2,4,-1,-1,0,0).
\end{align*}
All but the last multi-degree are stable and correspond to irreducible
components of the $\wh{P}$-compactified Jacobian. The last multi-degree is
strictly semistable. Let $Y \subset \wh{X}$ be a proper subcurve containing one irreducible component on top level and two irreducible components of the bottom level corresponding such that the multi-degree of the restriction to $Y$ is $(2,-1,-1)$. Then $\deg(\calF_Y)=0$, which yields equality in the stability condition to $Y$.  As above this multi-degree corresponds to a stratum in the boundary of the compactified Jacobian instead of an irreducible component.
\par
In this case the classical Hitchin fiber is also reducible. Let us shortly
explain how to see that. Let $(E,\Phi) \in \Hit^{-1}(X_\st,\alpha^2)$. Define
the eigen-line bundles 
\[ L_1\= \ker(\Phi-i\alpha \Id_E), \qquad L_2= \ker(\Phi+i\alpha\Id_E)\,..
\]
The semistability of $(E,\Phi)$ is equivalent to $\deg L_i \leq 0$. By \cite{GotOli} there is an open subset of the Hitchin fiber, where the Higgs field is non-vanishing, i.e.\ for all $x \in X: \Phi(x) \neq 0$. 
\begin{lemma} There is a morphism of coherent sheaves 
\[ E^\vee \rightarrow L_1 ^\vee \oplus L_2 ^\vee,
\] that is an isomorphism away from $Z(\alpha)$, such that the dualized Higgs field $\Phi^\vee$ is the pullback of the diagonal Higgs field $\mathsf{diag}(i\alpha,-i\alpha)$ on $L_1 ^\vee \oplus L_2 ^\vee$ along this map. Furthermore, for $\Phi$ non-vanishing we have $L_1 \otimes L_2 = \omega_X^{-1}$. 
\end{lemma}
\begin{proof}
We have an exact sequence of coherent sheaves
\[ 0 \rightarrow L_1 \oplus L_2 \xrightarrow{\iota_1 + \iota_2} E \rightarrow \mathcal{T} \rightarrow 0,
\] where $\iota_1,\iota_2$ are inclusions of subbundles and $\mathcal{T}$ is a torsion sheaf supported at $Z(\alpha)$. The Higgs field $\Phi$ on $E$ induces the diagonal Higgs field on $L_1 \oplus L_2$. The dual of the first map is the desired morphism. Generically, the line subbundles $L_1, L_2\subset E$ agree of order 1 on $Z(\alpha)$. In this case, $\mathcal{T}$ has a stalk of length 1 at $p \in Z(\alpha)$. Since~$\alpha$ is an abelian differential, this yields $\det(\mathcal{T})=\omega_X$ and hence $L_1 \otimes L_2 = \omega_X^{-1}$. This genericity condition is equivalent to the Higgs field being non-vanishing (cf. \cite[Theorem 5.5]{Horn}).
\end{proof}

The morphism of coherent sheaves described in the lemma is referred to as Hecke transformation (cf. \cite[Section 4]{Horn}). 
By the last formula of the lemma and since $g(X) = 3$, the stability of $(E,\Phi)$ in the open stratum is equivalent to 
\[ -4 < \deg(L_1) < 0.
\] Parameterizing the Hecke transformation, we obtain a description of the open stratum as a 
$(\CC^\times)^3$-fiber bundle over 
\[  \Pic^{-3}(X) \cup \Pic^{-2}(X) \cup \Pic^{-1}(X).
\] The closure of the three connected components of the open stratum are the irreducible components of $\Hit^{-1}(\alpha^2)$. These correspond to multi-degrees $(1,3,0,0,0)$, $(2,2,0,0,0,0)$, $(3,1,0,0,0,0)$. The other irreducible components of the $\wh{P}$-compacti\-fied Jacobian do not correspond to the original Hitchin fiber. 

\section{The fixed determinant case}
\label{sec:fixeddet}
The goal of this section is to deduce from the multi-scale spectral correspondence
that we proved in Theorem~\ref{thm:Higgspairsfam} a multi-scale version for
$\SL(2,\CC)$-Higgs pairs or more generally  multi-scale $\GL(2,\CC)^\circ$-Higgs
pairs with a fixed determinant~$\calL$. This correspondence boils down for smooth
curves and quadratic differentials with simple zeros to the classical case recalled in
Theorem~\ref{thm:Hitchinfiber}. For this purpose we have to define compactified
Prym varieties and determinants for torsion-free sheaves which are locally free except
for the special form given in Section~\ref{sec:localHiggs} at the nodes.
\par
There are difficulties to find a good notion of determinant for vector bundles on
nodal curves, see \cite[Section 8]{NagSeshI}, \cite{SunSLn} for solutions in some
cases. Similarly, for covers of curves with singular target one has to be careful
with the definition of the norm map, see e.g.\ \cite{CarbThesis} for various attempts.
For this reason we restrict to the situation of a family of quadratic multi-scale
differentials $(\wh{\pi}: \wh{\calX} \to \calX, \bfq, \bftau)$ on a fixed
Riemann surface $X_\st$ and over special basis. We will develop in this section
the notions for the following spectral correspondence.
\par
\begin{theorem}\label{thm:sl_spectral_corr}
Let $S$ be a reduced scheme and let $(\wh{\pi}: \wh{\calX} \to \calX,
\bfq, \bftau)$ be a family of quadratic multi-scale differentials over~$S$ with a
fixed underlying  Riemann surface~$X_\st$. Let $\calL$ be a line bundle on $X_\st$. The 
spectral correspondence of Theorem \ref{thm:Higgspairsfam} restricts to a
correspondence between $\wh{P}$-semistable torsion-free sheaves~$\calF$ on~$\wh{\calX}$
satisfying the $\calL$-twisted Prym condition and $P$-semistable $\GL(2,\CC)^\circ$-Higgs
pairs with fixed determinant~$\calL$. 
\end{theorem}
\par
The condition on the underlying  Riemann surface can be stated by the requiring
that the moduli map sends $S$ to the modified Hitchin base $\BB_{X_\st}$ for a fixed
Riemann surface~$X_\st$ of genus~$g$ rather then allowing $X_0$ to vary. We say
that such families are 'in  $\BB_{X_\st}$' for brevity. Concretely,
the fibers of~$\calX$ are given by $X_\st$ augmented by rational tails. 
\par
Both the Prym condition and the determinant condition depend on a line bundle~$\calL$
on~$X_\st$. Since the fibers of~$\calX \to S$ are curves of compact type, we may
extend~$\calL$ uniquely by the trivial bundle on the rational tails to a line
bundle on $X$ that we keep calling~$\calL$.
\par
In order to generalize the Prym condition from~\eqref{eq:defHitfiber} to families
of stable curves~$\wh{\calX}$ we have to allow to twist by components of the special
fiber. If the base~$S$ is a point, components of the special fiber are no longer 
divisors and hence we have to rephrase twisting by its effect on the normalization.
The definition requires some preparation. For an admissible double cover
$\wh{\pi}: \wh{X} \to X$, we denote by $\nu: X^\nu \to X$ respectively
by $\wh{\nu}: \wh{X}^\nu \to \wh{X}$
the normalizations and by $\wh{\pi}^\nu: \wh{X}^\nu \to X^\nu$ the induced covering of
smooth curves. We denote by $\nu^T$ (and $\wh{\nu}^T$) the torsion-free pullback, i.e.\
$\nu^T(\calG)=\nu^*\calG /\mathsf{Tor}(\nu^*\calG)$ for a coherent sheaf $\calG$ on $X$.
In the following, we identify the edges~$E(\Gamma)$ of the dual graph~$\Gamma$ with the corresponding nodes of~$X$. 
For each node~$e$ of ~$X$ (respectively ~$\wh{e}$ of~$\wh{X}$) we label its
preimages by $\nu^{-1}(e)=\{e_1, e_2\}$ (respectively~$\wh{\nu}^{-1}(\wh{e})=\{\wh{e}_1,\wh{e}_2\}$). 
Since~$\Gamma$ is a tree we can write the effect of twisting, a priori a sum
over vertices of~$\Gamma$, equivalently as a sum over its edges, see the second
sum in the following definition.
\par
\begin{definition} Let $(\wh{\pi}: \wh{X} \to X, \bfq, \bftau)$ be in
$\BB_{X_\st}$. Let $\calF \in \barcJ[\wh{g},\wh{n}][\wh{d},
\wh{P}](\wh{X})$ be a  torsion-free rank one sheaf. Let $\wh{N}$
be the set of nodes, where $\calF$ is not locally free. Then $\calF$ satisfies
the $\calL$\emph{-twisted Prym condition},  if for a choice (or equivalently for all
choices) $i(e) \in \{1,2\}$ of the preimage of each node $\wh{e} \in \wh{N}$ in the
normalization there exists $m_{e} \in \ZZ$, such that
\begin{align} \label{equ:Prym_cond}
\wh{\nu}^T \Bigl(\calF \otimes \sigma^* \calF(-\wh{B})\Bigr) 
\cong
(\wh{\pi} \circ \wh{\nu})^*\calL
\Bigl(-\sum_{\wh{e} \in \wh{N}} \wh{e}_{i(e)} + \sigma({\wh{e}}_{i(e)})
+ \sum_{e \in E(\Gamma)} m_{e} \sum_{\wh{e} \in \pi^{-1}(e)} \wh{e}_1-\wh{e}_2 \Bigr)\,.
\end{align} 
\par
The Prym condition defines a subfunctor of the compactified Jacobian for
families $(\wh{\pi}: \wh{\calX} \to \calX, \bfq, \bftau) \to S$ with fixed
underlying~$X_\st$ and arbitrary base~$S$. 
We define the \emph{compactified Prym functor} as the functor that associates to
a scheme $T \to S$ the set 
\[ \barPr[\wh{\pi}][\calL](T)=\left\{ \calF \in \barcJ[\wh{g},\wh{n}][\wh{d},\wh{P}](\wh{\calX} \times_S T) \middle|
\begin{array}{l}
\calF \ \text{satisfies the } \calL \text{-twisted Prym condition} \\  \text{\eqref{equ:Prym_cond} for all geometric points in } T 
\end{array}\right\}.
\]
\end{definition}
\par 
First we justify the claim that this condition is the specialization of the Prym
condition in Theorem \ref{thm:Hitchinfiber} under degeneration of the quadratic
differential, i.e.\ that the various Prym conditions are consistent in families.
\par
\begin{proposition} \label{prop:specialisation} Let $(\wh{\pi}: \wh{\calX}
\to \calX, \bfq, \bftau)$ be a germ of families of quadratic multi-scale differentials in $\BB_{X_\st}$ over a DVR $S$ with generically smooth fiber. Let $\eta$ denote the generic point of $S$. A torsion-free sheaf $\calF \in \barcJ[\wh{g},\wh{n}][\wh{d},\wh{P}](\wh{\calX})(S)$ satisfies the $\calL$-twisted Prym condition if and only if $\calF \otimes \sigma^* \calF(-\wh{B}) \cong \wh{\pi}^*\calL$ on $\wh{\calX}_{\eta}$. In particular, $\wh{d}=\deg \calL +2g-2$.
\end{proposition}
\par
The proposition relies on the following extension result. A crucial ingredient
is that the dual graph~$\Gamma$ of the special fiber of~$\calX$ is a tree.
\par
\begin{lemma}\label{lemma:free_extension} Let $(\wh{\pi}: \wh{\calX} \to \calX, \bfq, \bftau)$ be a germ of families of quadratic multi-scale differentials in $\BB_{X_\st}$ over a DVR $S$ with generically smooth fiber $\calX_\eta$ and special fiber $X$.
\begin{itemize}
\item[i)] Let $\calL$ be a family of line bundles on $\calX_\eta$. Then there exists an extension of $\calL$ to $\calX$ as a locally free rank 1 sheaf and all possible such extension differ by twisting with components of the special fiber $X$.
\item[ii)] Let $\calL$ be a family of line bundles on $\wh{\calX}_\eta$. Then there exists an extension of $\calL \otimes \sigma^*\calL$ as a locally free rank 1 sheaf and all possible such extension differ by twisting with components of the special fiber $\wh{X}=\wh{\calX}_s$.
\end{itemize}
\end{lemma}
\begin{proof}
  
\begin{itemize}[wide=\parindent]
\item[i)] Choose a polarization $\phi$ on $\calX$. Then there exists a semistable extension~$\calF$ of $\calL$ to $\calX$. Assume that $\calF$ is not locally free. Let $e \in X$ be a node where this happens. Then we can obtain this torsion-free rank one sheaf from a locally free rank 1 sheaf $\calG$ on a quasi-stable model of $X$, where the node $e$ is replaced by a rational bridge, such that $\calG$ has degree 1 on this rational bridge. By assumption the node $e$ separates $X=Y_1 \cup Y_2$ into disjoint subcurves. Denote by $d_i=\deg(\calF_{Y_i})$. Now it easy to see $\calG(-Y_1)$ has degree $d_1+1$ on $Y_1$, $d_2$ on $Y_2$ and degree $0$ on the rational bridge. Hence, $\calG(-Y_1)$ defines extension of $\calL$ that is locally free at the node~$e$. We also did not create any new nodes, where $\calG(-Y_1)$ is not locally free. Repeating this procedure at every other node, where $\calG$ is not locally free yields a locally free extension. On the other hand, if we have two extensions of $\calF_1,\calF_2$ of $\calL$ as locally free sheaves, then $\calF_1 \otimes \calF_2^{-1}$ is a locally free sheaf that is trivial, when restricted to $\calX_\eta$. Hence it is isomorphic $\calO(\sum n_i X_i)$ where the sum
is over all irreducible components $X_i \subset X$.
\item[ii)] The proof works similarly although~$\wh{\Gamma}$ is not a tree, thanks
to~$\sigma$-invariance of the line bundle to be extended: Take an extension $\calF$ of $\calL$ to $\calX$ as a torsion-free rank 1 sheaf. Let $\wh{e} \in \wh{X}$ be a node, such that $\calF$ is not locally free. Choose a quasi-stable model of $\wh{X}$ so that $\calF$ and $\sigma^*\calF$ are represented by locally free sheaves $\calG$ and $\sigma^*\calG$. If $\wh{e}$ is fixed by the involution $\sigma$, then $\calG \otimes \sigma^* \calG$ has degree $2$ on the rational bridge. Hence, it is equivalent to a torsion-free sheaf with degree $0$ on the rational bridge by twisting with the rational component. If $\wh{e}$ is not fixed by the involution then $\sigma^* \calF$ is not locally free at $\sigma \wh{e}$. In particular, $\calG \otimes \sigma^* G$ is a locally free sheaf that has the same degree on the two rational bridges corresponding to $\wh{e}$ and $\sigma\wh{e}$. The pair of nodes $\wh{e}, \sigma\wh{e}$ cuts $\wh{X}$ into two $\sigma$-invariant subcurves $\wh{Y}_1$ and $\wh{Y}_2$. (The preimages of the components $Y_1,Y_2$ as defined above using the node $\wh{\pi}(\wh{e})=\wh{\pi}(\sigma\wh{e})$.) $\calG \otimes \sigma^* G$ is equivalent to a sheaf with degree $0$ on both rational bridges by twisting with one of the components $\wh{Y}_1$ or $\wh{Y}_2$. Hence, $\calL \otimes \sigma^*L$ has an extension that is locally free at $\wh{e}$ and $\sigma\wh{e}$. We did not change $\calF \otimes \sigma^*\calF$ at any other node. Hence, by induction we obtain a locally free sheaf on $\calX$ extending $\calL \otimes \sigma^*\calL$. The uniqueness up to twisting by components works as before. 
\end{itemize}
\end{proof}
\par
\begin{proof}[Proof of Proposition~\ref{prop:specialisation}]
By Lemma \ref{lemma:free_extension} we can find a locally free extension $\calG$
of $\calF \otimes \sigma^*\calF$. By the proof of the lemma, the restriction of $\calG$
to the components  differs from the restriction of $\calF \otimes \sigma^*\calF$
by twisting with ${\wh{e}}_{i(e)}+\sigma^*\wh{e}_{i(e)}$ for certain preimages
determined by $i(e) \in \{1,2\}$ for all $\wh{e} \in \wh{N}$. On the other hand,
$\calL$ is a locally free extension. Hence, again by Lemma \ref{lemma:free_extension}
there exist $m_i \in \ZZ$, such that 
\[ \calG(-\wh{B}) \= \wh{\pi}^*\calL
\Bigl(\sum_{\wh{X}_i \subset \wh{X}\ \mathrm{ irred.}} m_i \wh{X}_i\Bigr).
\] Restricted to the special fiber of $\wh{\calX} \to S$ twisting by the irreducible
component $\wh{X}_i$ yields a twist of the torsion-free pullback to the normalization
$\wh{X}^\nu$ by $\sum_{ \wh{e} \in \wh{X}_i} (\wh{e}_1-\wh{e}_2)$. As $\wh{\Gamma}$ is
obtained from the tree $\Gamma$ by doubling some edges, this is equivalent to twisting
by $\sum_{e \in \Gamma} m_{e} \sum_{\wh{e} \in \pi^{-1}(e)} (\wh{e}_1-\wh{e}_2)$. All together
we obtain formula~\eqref{equ:Prym_cond}. For different choices of the preimages
${\wh{e}_{i(e)}} \in \wh{X}^{\nu}$ in \eqref{equ:Prym_cond} we can compensate by
varying the~$m_{e}$.  
\end{proof}
\par
\begin{proposition} Let $(\wh{\pi}: \wh{\calX} \to \calX, \bfq, \bftau) \to S$
be a family of quadratic multi-scale differentials in $\BB_{X_\st}$. Let $\wh{d}=\deg \calL-2g+2$.
Then $\barPr[\wh{\pi}][\calL]$ is universally closed in $\barcJ[\wh{g},\wh{n}][\wh{d},\wh{P}]$.
It is a proper Deligne-Mumford stack, if $\wh{P}$ is non-degenerate.
\end{proposition}
\par
\begin{proof} For a family of quadratic differentials on $X_\st$ with simple zeros
and thus a family of smooth curves $\wh{X} \to S$ the result follows from the properness of the Prym variety.
Otherwise let $\calF \in \barcJ[\wh{g},\wh{n}][\wh{d}](\wh{\calX})[T]$ be a family of torsion-free sheaves over a DVR $T$ that satisfies the $\calL$-twisted Prym condition on the generic point $\eta_T$ of $T$. We may assume that $\calX \times_S T$ is singular over the special point $t \in T$. Let $\wh{N}_{\text{per}}=\{ n_i: T \to \wh{\calX} \times_S T \}$ be the set of persistent nodes over $T$. We consider the partial normalization $f: \wh{\calY} \to \wh{\calX}$ at all nodes in $\wh{N}_{\text{per}}$. Then the family $\wh{\calY} \to T$ has generically smooth fiber and the Prym condition~\eqref{equ:Prym_cond} is a condition on $f^T\calF$ over the generic point. When the family of smooth curves $\wh{\calY} \to T$ develops new nodes over the special point of $t \in T$, then $\calF_t$ still satisfies the $\calL$-twisted Prym condition by Proposition~\ref{prop:specialisation}.
\par
We are left with considering the case, when $\calF$ is locally free at a persistent node $n_j(\eta_T)$ but Neveu-Schwarz at $n_j(t)$. Let $\wh{\calY}_1,\wh{\calY}_2$ be the connected components of $\wh{\calY}$ meeting in the node $n_j$. If a locally free sheaf becomes Neveu-Schwarz at $n_j(s)$, then the degree $\deg(f^T\calF_{\calY_i})$ drops by $1$ for one of the components $\calY_i$ and stays constant on the other component. (See the
alternative viewpoint using quasi-stable curves in Section~\ref{sec:UnivJac} for an explanation.)
Let $p_{n_j(s)}$ be the preimage of the node in $\calY_i$. Then the degree drop is compensated in~\eqref{equ:Prym_cond} by twisting with $p_{n_j(s)} + \sigma p_{n_j(s)}$. This shows that by choosing the correct preimage of the node $n_j(s)$ we can ensure that the formula~\eqref{equ:Prym_cond} is constant under such degeneration. In particular, $\calF_s$ still
satisfies the $\calL$-twisted Prym condition. 
\end{proof}
\par
For a multi-scale $\GL(2,\CC)^\circ$-Higgs pair $(\calE,\bfPhi)$ on a
quadratic multi-scale differential $(\wh{\pi}: \wh{X} \to X, \bfq, \bftau)$
in $\BB_{X_0}$ we partition the set~$N$ of nodes where $\calE$ is not locally free
into subsets $N=N_{o,f} \cup N_{o,n} \cup N_{e,n} \cup N_{e,nn}$ according
to ramification and the local structure as follows. Let
\begin{itemize}
\item[i)] $N_{o,f}$ be the set of nodes, where $\kappa$ is odd and it is of the form \eqref{eq:LF/fixframe}, 
\item[ii)] $N_{o,n}$ be the set of nodes, where $\kappa$ is odd and it is of the form \eqref{eq:NS/fixframe},
\item[iii)] $N_{e,n}$ be the set of nodes, where $\kappa$ is even and it is of the form \eqref{eq:calIswap} and 
\item[iv)] $N_{e,nn}$ be the set of nodes, where $\kappa$ is even and it is the pushforward of torsion-free sheaf that is Neveu-Schwarz at both nodes in $\wh{\pi}^{-1}(e)$.
\end{itemize} 
\par
\begin{definition}
Let $(\wh{\pi}: \wh{X} \to X, \bfq, \bftau)$ be a quadratic multi-scale differential
in $\BB_{X_\st}$. Let $(\calE,\bfPhi)$ be a multi-scale $\GL(2,\CC)^\circ$-Higgs pair on $X$.
Then $(\calE,\bfPhi)$ has \emph{fixed determinant $\calL$} if and only if for a choice
(or equivalently all choices) $i(e) \in \{1,2\}$ of the preimage of a node
$e \in N$ there exists $m_{e} \in \ZZ$, such that
\begin{align}\label{equ:det:cond} \det(\nu^T(\calE))
\Bigl(\sum_{e \in N_{o,f} \cup N_{e,n}} e_{i(e)} + \sum_{e \in N_{o,n}\cup N_{e,nn}} 2 e_{i(e)} \Bigr)
\=\nu^*\calL\Bigl(\sum_{e \in E} m_{e} (e_1-e_2)\Bigr)\,.
\end{align}
\par
We call $(\calE,\bfPhi)$ a \emph{multi-scale $\SL(2,\CC)$-Higgs pair} if it has fixed
determinant equal to~$\calO_{X_0}$.
\end{definition}
\par
Again we want to show that this condition is a specialization of the fixed
determinant condition for smooth curves.
\par
\begin{proposition}
Let $(\calE, \boldsymbol{\Phi})$ be a multi-scale $\GL(2,\CC)^\circ$-Higgs pair
on a germ of families of quadratic multi-scale differentials $(\wh{\pi}: \wh{\calX} \to \calX, \bfq, \bftau) \to S$ in $\BB_{X_\st}$ over a DVR $S$ with generically smooth fiber. If $(\calE,\bfPhi)$ has fixed determinant $\calL$ on the generic fiber $\calX_{\eta}$ then it satisfies \eqref{equ:det:cond} on the special fiber. In particular, $\deg\calE=\deg{\calL}$.
\end{proposition}
\par
\begin{proof}
By Lemma \ref{lemma:free_extension} we can find a locally free extension $\calG$
of $\det(\calE)$ from the generic fiber $\calX_\eta$ to $\calX$. Since $\calL$ is another
locally free extension there exist $m_i \in \ZZ$, such that 
\[ \calG \= \calL\Bigl(\sum_{X_i \subset X\ \mathrm{ irred.}} m_i X_i\Bigr)\,.
\] By restriction to the special fiber we obtain
\[ \nu^*\calG \= \nu^*\calL\Bigl(\sum_{e \in E(\Gamma)} m_e (e_1-e_2)\Bigr)\,.
\] To relate $\nu^*\calG$ to $\det(\nu^*\calE)$ we compare both of them to $\bigwedge^2
\calE$. This sheaf is not torsion-free in general. We recover $\calG$ by first taking
the reflexive hull of $\bigwedge^2 \calE|_{\calX \setminus E(\Gamma)}$ and then in the 'second step' applying 
Lemma~\ref{lemma:free_extension} to the resulting torsion-free rank 1 sheaf.  For a
node~$e$ let~$X_{e_1}$ and~$X_{e_2}$ be the connected components of the normalization
of~$X$ at~$e$, so that $e_1 \in X_{e_1}$ and $e_2 \in X_{e_2}$. For every node~$e$, where
the reflexive hull is Neveu-Schwarz the 'second step' involves the choice of one of the connected
components $X_{e_i}$ by which we twist to obtain a locally free extension. For the
comparison with $\det(\nu^T\calE)$ it is necessary to remember these choices.  
\begin{itemize}
\item If $e \in N_{o,f}$, then $\calE$ is locally in a neighborhood $U_e$ of $e$
isomorphic to $R \oplus \langle x,y \rangle$. Then $\bigwedge^2\calE\cong \langle
1 \otimes x, 1\otimes y \rangle$ is Neveu-Schwarz at $e$. Let $\nu^{-1}U=U_{e_1}
\sqcup U_{e_2}$ be the connected components of the preimage in the normalization. Assume
we have twisted in the 'second step' by~$X_{e_1}$. Using the isomorphism $\langle x,y
\rangle \cong \langle u,t \rangle$ defined by multiplying by $x$ we see that
$\calG_{U_{e_1}}=\det(\calE_{U_{e_1}})(e_1)$ and $\calG_{U_{e_2}}=\det(\calE_{U_{e_2}})$.
Alternatively, in case we twisted with $X_{e_2}$ in the 'second step', we use the
isomorphism $\langle x,y \rangle \cong \langle t,v \rangle$ and recover
$\calG_{U_{e_1}}=\det(\calE_{U_{e_1}})$ and $\calG_{U_{e_2}}=\det(\calE_{U_{e_2}})(e_2)$.
\item Let $e \in N_{o,n}$ and assume the degree of $\calE$ has dropped at $X_{e_2}$.
Then locally at~$e$ the rank two torsion-free sheaf~$\calE$ is isomorphic
to $\langle 1_u,v \rangle \oplus \langle x,y \rangle \cong \langle 1_u,v \rangle
\oplus \langle u,t \rangle$ (or $\langle 1_u,v\rangle \oplus \langle t,v\rangle$). Here $1_u,1_v$ denote the unit in the
coordinate rings $\calO_{U_{e_1}} \cong \CC[[u]]$ resp.\ $\calO_{U_{e_1}} \cong \CC[[v]]$.
When choosing to twist by~$X_{e_1}$ (resp.\ by~ $X_{e_2}$) in the 'second step', we
obtain $\calG_{U_{e_1}}=\det(\calE_{U_{e_1}})(e_1)$ and $\calG_{U_{e_2}}=\det(\calE_{U_{e_2}})(e_2)$
(resp.\ we obtain $\calG_{X_{e_1}}=\det(\calE_{X_{e_1}})$ and
$\calG_{X_{e_2}}=\det(\calE_{X_{e_2}})(2e_2)$).
However up to twisting by $e_2-e_1$ there is no difference between these two choices. 
\item For the cases of $e \in N_{e,n}$ and $e \in N_{e,nn}$ the argument works exactly
as in the proof of Lemma~\ref{lemma:free_extension}. Here $\calE$ can be identified
with $\calF \oplus \sigma^*\calF$ locally at~$e$. 
\end{itemize}
Summing up we obtain the formula in \eqref{equ:det:cond}.
\end{proof}
\par
\begin{proof}[Proof of Theorem \ref{thm:sl_spectral_corr}]
We formulated  the Prym condition and the fixed determinant condition fiberwise. Hence
it is enough to relate the two notions on a single quadratic multi-scale differential
$(\wh{\pi}: \wh{X} \to X, \bfq, \bftau)$ in $\BB_{X_\st}$. For a torsion-free sheaf~$\calF$ on~$\wh{X}$
pushforward and pullback to the normalization commute, i.e.\ $\nu^T\wh{\pi}_* \calF
= \wh{\pi}^{\nu}_* \wh{\nu}^T \calF$. The result essentially follows from
formula~\eqref{eq:detEpullback} since the map $\wh{\pi}^{\nu*}$ is still injective
because restricted to a connected component of $\wh{X}^\nu$ this map is not unbranched.
To give the details it suffices to check the correspondence  for $\calF$ a locally
free sheaf on $\wh{X}$ because both formulas~\eqref{equ:Prym_cond} and~\eqref{equ:det:cond}
continue to hold under degeneration of $\calF$ to a torsion-free sheaf. Notice that for $\calF$ a locally free sheaf $N_{o,n}=N_{e,n}=N_{e,nn}=\varnothing$.
We convert the fixed determinant condition into the Prym condition as follows.
\begin{align*} &\det(\nu^T(\calE))\Bigl(\sum_{e \in N_{o,f}} e_{i(e)}\Bigr)&&\=
\nu^*\calL\Bigl(\sum_{e \in E(\Gamma)} m_{e} (e_1-e_2)\Bigr) \\
\Leftrightarrow \quad &\Bigl(\wh{\pi}^{\nu*}\det(\nu^T\wh{\pi}_* \calF)\Bigr)
\Bigl(\sum_{e \in N_{o,f}} 2\wh{e}_{i(e)}\Bigr)
&&\= (\wh{\pi}^{\nu*}\nu^*\calL)\Bigl(\sum_{e \in E(\Gamma)} m_{e}\sum_{\wh{e} \in \pi^{-1}(e)} (\wh{e}_1-\wh{e}_2)\Bigr) \\ 
\Leftrightarrow \quad &(\nu^T \calF \otimes \sigma^* \nu^T \calF)
\Bigl(-\nu^*\wh{B}+\sum_{e \in N_{o,f}}&&2\wh{e}_{i(e)} - \wh{e}_1-\wh{e}_2\Bigr) \\
& &&\= (\wh{\nu}^*\wh{\pi}^*\calL)\Bigl(\sum_{e \in E(\Gamma)} m_{e}\sum_{\wh{e} \in \pi^{-1}(e)} (\wh{e}_1-\wh{e}_2)\Bigr) 
\end{align*} Note that for $e \in N_{o,f}$ the preimages in the normalization $\wh{e}_i$
are branch points of~$\wh{\pi}^\nu$. This explains the extra twist on the left hand
side. Finally, $ 2\wh{e}_{i(e)} -\wh{e}_1-\wh{e}_2$ is of the form $\pm(\wh{e}_1-\wh{e}_2)$.
Bringing these divisors to the right, we see that there exists $m_{e}' \in \ZZ$, such
that condition~\eqref{equ:Prym_cond} is satisfied. This proves the theorem. 
\end{proof}

\newpage
\paragraph{\textbf{Summary of Notation:}} \
\par \medskip
\begin{tabular}{l|l}
$X_\st$, $(X,\bfz)$ & $X_\st \in \overline{\calM}_g$ and $(X,\bfz) \in \barmoduli[g,n]$
base stable (pointed) curve \\
$\pi: \Sigma \rightarrow X$ & spectral cover of $(X,q) \in \calQ_{g}^+(\mu)$ \\
$\wh{\pi}: \wh{X} \rightarrow X$ & canonical double cover,
with marked points $\wh{\bfz}$\\
$\wh{g}$ & $\wh{g}=4g-3$ genus of $\wh{X}$ \\
$e$ resp. $\wh{e}$ & node of pointed stable curve $X$ respectively $\wh{X}$ \\
$N$ resp. $\wh{N}$ & set of nodes of $X$ respectively $\wh{X}$ \\
$B$ resp. $\wh{B}$ & branch divisor resp. ramification divisor of $\wh{\pi}$\\
$\sigma: \wh{X} \rightarrow \wh{X}$ & involution on $\wh{X}$ interchanging the sheets \\ 
$\nu: X^\nu \rightarrow X$ & Normalisation of $X$ \\
$\wh{\nu}: \wh{X}^\nu \rightarrow \wh{X}$ & Normalisation of $\wh{X}$ \\
$\omega(\bfz)$, $\wh{\omega}(\wh{\bfz})$ & twisted dualizing sheaf
$\omega_X(\sum_{i=1}^{n} z_i)$
resp. $\omega_{\wh{X}}(\sum_{i=1}^{n} \wh{z}_i)$ \\
 $(X,q)$, $(X,q,\bfz)$ & twisted quadratic differential
in $\overline{\calQ}_{g}^+(\mu)$ resp.\ in  $\overline{\calQ}_{g,n}^+(\mu)$ \\
$X_i$ or $\wh{X}_i$ & $i$-level of $X$ or $\wh{X}$\\
$\calX_{i^\complement}$ & $\calX \setminus \cup_{j \neq i} X_j$\\
$\bftau$ & prong-matching \\
$L, \calL$ & locally free sheaf of rank 1 \\
$\calF$ & coherent sheaf, often on $\wh{X}$ \\
$\calE$ & vector bundle or special rank two bundle, often on ${X}$ \\
$M$ & locally free sheaf of rank 1, twist of Higgs field \\
$\Lambda$ & locally free sheaf of rank 1, fixed determinant   \\

\end{tabular}

\printbibliography

\end{document}